\documentclass[11pt]{amsart}
\usepackage{amssymb,amsmath,amsthm,amsfonts,amsopn,url,color,hyperref,enumerate,mathtools,microtype,MnSymbol,bm,array,comment,wrapfig,caption,subcaption}
\usepackage[top=1in, bottom=1in, left=1in, right=1in]{geometry}
\usepackage[normalem]{ulem}
\usepackage[all]{xy}
\input xy
\xyoption{all}
\usepackage{amscd}
\usepackage{tikz}
\usepackage{float}
\usepackage{graphicx}
\usepackage{multicol}

\newtheorem{thm}{Theorem}[section]

\newtheorem{prop}[thm]{Proposition}

\theoremstyle{definition}
\newtheorem{defn}[thm]{Definition}

\newtheorem{ex}[thm]{Example}

\numberwithin{equation}{section}
\numberwithin{figure}{section}
\numberwithin{table}{section}

\newcommand{\df}[2]{ \left(#1, #2 \right)! }

\newcommand\<{\langle}
\renewcommand\>{\rangle}
\newcommand\CC{\ensuremath{\mathbb C}}
\newcommand\II{\ensuremath{\mathbb I}}
\newcommand\NN{\ensuremath{\mathbb N}}

\newcommand\RR{\ensuremath{\mathbb R}}
\newcommand\ZZ{\ensuremath{\mathbb Z}}

\newcommand{\bba}{\ensuremath{\bm{a}}}
\newcommand{\bbc}{\ensuremath{\bm{c}}}
\newcommand{\bbd}{\ensuremath{\bm{d}}}

\newcommand{\bbm}{\ensuremath{\bm{m}}}
\newcommand{\bbt}{\ensuremath{\bm{t}}}
\newcommand{\bbx}{\ensuremath{\bm{x}}}
\newcommand{\boldzero}{\ensuremath{\bm{0}}}

\title[An illustrated view of differential operators]{An illustrated view of differential operators of a reduced quotient of an affine semigroup ring}
\author[Berkesch]{Christine Berkesch}
\address{Department of Mathematics \\ University of Minnesota\\ Minneapolis, MN}
\email{cberkesc@umn.edu}
\author[Chan]{C-Y. Jean Chan}
\address{Department of Mathematics\\ Central Michigan University \\ Mt. Pleasant, MI}
\email{chan1cj@cmich.edu}
\author[Klein]{Patricia Klein}
\address{Department of Mathematics \\ University of Minnesota\\ Minneapolis, MN}
\email{klein847@umn.edu}
\author[Matusevich]{Laura Felicia Matusevich}
\address{Department of Mathematics \\ Texas A \& M \\ College Station, TX}
\email{laura@math.tamu.edu }
\author[Page]{Janet Page}
\address{Department of Mathematics\\ University of Michigan \\ Ann Arbor, MI}
\email{jrpage@umich.edu}
\author[Vassilev]{Janet Vassilev}
\address{Department of Mathematics and Statistics \\ University of New Mexico \\ Albuquerque, NM}
\email{jvassil@math.unm.edu }

\thanks{CB was partially supported by NSF DMS 2001101. 
LFM was partially supported by the Simons Foundation.
}
\begin{document}
\begin{abstract}
Through examples, we illustrate how to compute differential operators on a quotient of an affine semigroup ring by a radical monomial ideal, when working over an algebraically closed field of characteristic 0. 
\end{abstract}
\maketitle

\setcounter{section}{1}
\section*{Introduction}
\label{sec:intro}

In this paper, we provide illustrative examples and visualizations of some differential operators on the quotient of an affine semigroup ring by a radical monomial ideal.  
These examples motivate our work in \cite{BCKMPVJournal}.
For a finitely-generated algebra $R$ over a field, let $D(R)$ denote the ring of differential operators of $R$ and use $\ast$ to denote an action of a differential operator. 
A foundational result in this area relates the ring of differential operators of an arbitrary quotient of a polynomial ring by an ideal $J$ to differential operators on that polynomial ring and the idealizer of $J$. 

\begin{prop}
\cite[Propostion~1.6]{SmStDO} 
\label{prop:SmStDO-Prop1.6}
Let $S$ be the coordinate ring of a nonsingular affine variety over an algebraically closed field of characteristic $0$, and let $J$ be an $S$-ideal. Then there is an isomorphism 
\vspace{-2mm}
\begin{equation}
D\left(\frac{S}{J}\right) \cong \frac{\II( J )}{ JD(S)  }, 
\label{eq:regularDO}
\end{equation}
where $\II(J):=\{\delta \in D(R) \mid \delta \ast J \subset J \}$, which is called the \emph{idealizer} of $J$.
\qedhere
\end{prop}

The phrasing of \cite[Proposition~1.6]{SmStDO} given above differs from the original by making use of the equivalence of the conditions, with notation as above, of $\delta J D(S) \subset JD(S)$ and $\delta \ast J \subset J$ \cite[Lemma~2.3.1]{Sai-Tr-DASR}. 
Traves \cite{TrDM} uses Proposition~\ref{prop:SmStDO-Prop1.6} to give concrete descriptions of rings of differential operators of 
Stanley--Reisner rings, and Saito and Traves \cite{Sai-Tr-DASR} use the same to compute rings of differential operators of affine semigroup rings. 
For a non-regular affine semigroup ring $R_A$ over an algebraically closed field of characteristic $0$, Proposition~\ref{prop:SmStDO-Prop1.6} fails, even for a radical monomial ideal $J$ in $R_A$. 
However, the differential operators on $R_A$ that induce maps on $R_A/J$ are precisely those in $\II(J)$. 
Further, there is an embedding of rings 
\begin{align*}
\label{eq:key-embedding}
\frac{\II(J)}{D(R_A,J)}\hookrightarrow 
D\left(\frac{R_A}{J}\right),
\end{align*}
so the description of $D(R_A)$ by Saito and Traves can still be used to compute this subring of $D(R_A/J)$. 

One primary goal of this article is to visually illustrate the computation of $\II(J)/D(R_A,J)$ when $R_A$ is an affine semigroup ring and $J$ is a radical monomial ideal in $J$; in \cite{BCKMPVJournal}, we provide an explicit formula for this computation. 
In particular, the pictures we provide explain how to compute differential operators of the form 
\[
D(I,J):=\{ \delta \in D(R) \mid \delta \ast I \subseteq J\},
\]
where $I$ and $J$ are subsets of the ring $R$, and $\ast$ denotes an action by a differential operator. 
These sets were originally instrumental in both the work of Smith and Stafford \cite{SmStDO} and of Musson \cite[Section 1]{Mussonzd}, and they are essential building blocks of our  computations, as well.  

Our second major goal in this paper is to compare $JD(R_A)$ and $D(R_A,J)$ for an affine semigroup ring $R_A$ and radical monomial ideal $J$ of $R_A$.  
Towards this end, the first set of examples we consider consists of quotients of the coordinate rings of rational normal curves. 
These quotients are all isomorphic to $\CC[x,y]/\<x y\>$, 
 which is handled in~\cite{TrDM}. 
From the standpoint of comparing $JD(R_A)$ and $D(R_A,J)$, what we will see is that 
$JD(R_A)=D(R_A,J)$ for the rational normal curves of degrees $1$ and $2$ but fail to coincide in all degrees larger than $2$.  
For rational normal curves of degree at most $2$, 
$J$ is principally generated and so is $JD(R_A)$.
In this case, it is straightforward to see (by definition) that $D(R_A,J)$ is also principal and is isomorphic to $JD(R_A)$. 
However, for degree $n>2$, $J$ has $n-1$ generators. 
We will see that this number of generators greatly impacts $\II(J)/JD(R_{A})$ but not $\II(J)/D(R_{A},J)$. 

In \cite{BCKMPVJournal}, 
we consider cases where  
\begin{equation}\label{eq:JDDJequal}
JD(R_A) = D(R_A, J),
\end{equation} 
for the particular ideal $J$ generated by all monomials corresponding to the interior of the semigroup, i.e., $J = \omega_{R_A}$. 
When $R_A$ is Gorenstein and normal,  ~\eqref{eq:JDDJequal} holds since $\omega_{R_A}$ is principal. We show that the converse is also true; that is,~\eqref{eq:JDDJequal} is equivalent to $R_A$ being a Gorenstein ring. 

To provide intuition beyond the two-dimensional case, we also offer a three-dimensional normal affine semigroup ring $R_A$ for which we compute $\II(J)/D(R_A,J)$ for two choices of $J$. 
Then, we return to the two-dimensional setting to consider a scored but not normal example, as well as a non-scored example, computing differential operators for quotients of both rings. 

\subsection*{Outline}
Section~\ref{sec:back} fixes notation to be used throughout the article and describes the main result of~\cite{Sai-Tr-DASR}.  
Sections~\ref{sec:2D},~\ref{sec:3D}, and~\ref{sec:higherD} describe 
$\II(J)/D(R_A,J)$ 
for rational normal cones of degrees 2, 3, and higher, respectively. 
Section \ref{sec:hdim} considers a three-dimensional normal affine semigroup ring modulo two different choices of radical monomial ideal $J$,  
and Section \ref{sec:nonnormal} computes $\II(J)/D(R_A,J)$ for two different non-normal two-dimensional affine semigroup rings, where $J = \omega_{R_A}$ is the radical monomial ideal corresponding to the interior of the semigroup $\NN A$. 

\subsection*{Acknowledgements}
We are grateful to 
the organizers, Karen Smith, Sandra Spiroff, Irena Swanson, and Emily Witt, of \emph{Workshop: Women in Commutative Algebra (WICA)}, at which this work began. Additionally, we are thankful to BIRS for their hospitality in hosting this meeting. 
We also owe thanks to Jack Jeffries, who pointed us to an example that improved our work.  We truly appreciate the comments of the referee which added extra insight into our work and greatly improved the paper.

\section{Background and Notation}
\label{sec:back}

In this section, we fix notation and conventions to be assumed throughout the article.  Although the results we discuss in this paper hold over any algebraically closed field of characteristic 0, we will use in this illustrated view the field of complex numbers for the sake of simplicity. Having fixed notation, we will then state and discuss \cite[Theorem 3.2.2] {Sai-Tr-DASR}.

\begin{defn}
\label{def:semigp-ring}
Let $A$ be a $k \times \ell$ matrix with entries in $\ZZ$. 
Let  $\NN A$ denote the semigroup in $\ZZ^k$ that is generated by the columns of $A$. 
The \emph{affine semigroup ring} determined by $A$ is 
\[
R_A=\CC[\NN A]= \bigoplus_{{\bba} \in \NN A} \CC\cdot  {\bbt}^{\bba},
\]
where ${\bbt}^{\bba}=t_1^{a_1}t_2^{a_2} \cdots t_k^{a_k}$ for $\bba= (a_1,a_2, \ldots, a_k)\in\NN A$. 
Throughout this article, we assume that the group generated by the columns of $A$ is the full ambient lattice, so $\ZZ A = \ZZ^k$, and also that the real positive cone over $A$, $\RR_{\geq 0} A$, is pointed, meaning that it is strongly convex.  
\end{defn}

\begin{defn}
\label{def:normal-scored}
A semigroup $\NN A$ is \emph{normal} if $\NN A = \RR_{\geq 0} A \cap\ZZ A$. 
A semigroup $\NN A$ is \emph{scored} if the difference $(\RR_{\geq 0}A \cap \ZZ A) \setminus \NN A$ consists of a finite union of hyperplane sections of $\RR_{\geq 0} A \cap \ZZ A$ that are all parallel to facets of the cone $\RR_{\geq 0} A$. 
An affine semigroup ring $\CC[\NN A]$ is said to be \emph{normal}, or \emph{scored}, if $\NN A$ is normal, or scored, respectively. 
Note that normal semigroups are scored. 
\end{defn}

When we write a facet $\sigma$ of $\RR_{\geq 0} A$ (or $A$ or $\NN A$),
we will always mean the integer points in the corresponding facet of $\RR_{\geq 0}A$. 
When $A$ is normal, this is the same as the semigroup generated by the columns of $A$ that lie in the corresponding facet of $\RR_{\geq 0} A$.

Throughout, we use $\<-\>$ to indicate ideals in the commutative rings $R_A$ or the polynomial ring $\CC[\theta]=\CC[\theta_1, \ldots, \theta_k]$. 
It will be clear from the context if the ideals are in $R_A$ 
or in $\CC[\theta]$.

Notice that the $\ZZ^k$-graded prime ideals in $\RR_A$ are in one-to-one correspondence with the faces of $A$ (or $\RR_{\geq 0} A$), as a face $\tau$ of $A$ corresponds to the multigraded prime $R_A$-ideal 
\[
P_\tau := \left\<\bbt^{\bbd}\, \big\vert \,  \bbd\in\NN A\setminus \tau \right\>. 
\]
In this paper, we compute rings of the form $\II(J)/D(R_A,J)$, where $J$ is a radical monomial ideal in $R_A$; as such, $J$ is always as intersection of primes of the form $P_{\tau}$, for various faces $\tau$ of $A$.

We are mainly following the description
of Saito and Traves \cite{Sai-Tr-DASR}, 
although Musson and Van den Bergh describe the ring of differential operators of a toric ring  $\CC[\NN A]$  first in \cite{Mussontori}, \cite{Mussontv} and \cite{MuVDBTor}  when viewed as a subring of the ring of differential operators of the Laurent polynomials, i.e., $D(\CC[\ZZ^k])=\CC\{t_1^{\pm 1}, \ldots, t_k^{\pm 1}, \partial_1, \ldots, \partial_k\}/\sim$, 
where $\partial_i$ denotes the differential operator $\frac{\partial}{\partial t_i}$ 
and $\sim$ denotes the usual relations on the free associative algebra $\CC\{t_1^{\pm 1}, \ldots, t_k^{\pm 1}, \partial_1, \ldots, \partial_k\}$ that describe the behavior of differential operators.
To explain 
the differential operators of  $\CC[\NN A]$ as presented by Saito and Traves,
note first that $D(\CC[\ZZ^k])$ carries a $\ZZ^k$-grading, where $e_i$ denotes the $i$-th column of the identity matrix $I$ and $\deg(t_i) = e_i = -\deg(\partial_i)$.   Note also that if $\bba_i$ is a column of $A$, then $\deg(\bbt^{\bba_i})=\bba_i$.
Set $\theta_j=t_j\partial_j$ for $1 \leq j \leq k$, and 
set 
\[
\Omega({\bbd}) 
:= \{\bba \in \NN A \mid \bba +{\bbd} \notin \NN A\} 
= \NN A \setminus (-{\bbd}+ \NN A). 
\]
We note also that for any $\bbt^{\bbm} \in R_A$ and $f(\theta) \in \CC[\theta]$, $f(\theta) \ast \bbt^{\bbm} = f(\bbm) \bbt^{\bbm}$.
The \emph{idealizer of $\Omega ({\bbd})$} is defined to be
\[ 
\II( \Omega({\bbd}) ) 
:= \<f (\theta) \in \CC[\theta]=\CC[\theta_1, \theta_2 \ldots, \theta_k] 
\mid f(\bba)=0 \text{ for all } \bba\in \Omega({\bbd})\>, 
\]
which is viewed as an ideal in the ring $\CC[\theta]$, where $\theta_i = t_i\partial_i\in D(\CC[\ZZ^d])$ is of degree $\boldzero$. 
We will soon see that $\II( \Omega({\bbd}))$ consists of $f(\theta)$ such that $\bbt^{\bbd} f(\theta) \in D(R_A)$. 

To compute $\II(\Omega({\bbd}))$ for a normal semigroup ring, consider a facet $\sigma$ of $A$, recalling that by this we mean all lattice points on the corresponding facet of the cone $\RR_{\geq 0} A$.
The {\em primitive integral support function} $h_{\sigma}$ is a uniquely determined linear form on $\RR^k$ such that: 
\begin{enumerate}
    \item $h_{\sigma}(\RR_{\geq 0}A) \geq 0$,
    \item $h_{\sigma}(\RR \sigma)=0$, and 
    \item $h_{\sigma}(\ZZ^k)=\ZZ.$
\end{enumerate}
We write $ \sigma_1, \sigma_2, \ldots, \sigma_m$ for the facets of $A$, 
so that for the remainder of the paper, we set 
\[
h_j := h_j(\theta) = h_{\sigma_j}(\theta). 
\] 

For a non-negative integer $n$, we will use the following notation to denote this {\em descending factorial}: 
\[
\df{\alpha}{n} 
:= \prod\limits_{i=0}^{n}(\alpha - i) 
= \alpha ( \alpha -1) \cdots(\alpha - n),
\]
where $\alpha$ is a function, which could be  constant or already evaluated. 
For example, 
\[
\df{h_j}{h_j({-\bbd})-1)}
= \prod\limits_{i=0}^{h_j({-\bbd})-1}(h_j-i)
\]
will be a common expression throughout this article, as it is a factor in the generator of the idealizer $\II(\Omega(\bbd))$.   To streamline our presentation, we make the convention that $\df{\alpha}{n} = 1 $ for all $n <0$.   

\begin{thm}
\cite[Theorem 3.2.2] {Sai-Tr-DASR} 
\label{thm:ST-dasr-3.2.2}
If $R_A$ is a pointed, normal affine semigroup ring with $\ZZ A = \ZZ^k$, then 
\begin{align}
\nonumber
D(R_A) 
&= 
\bigoplus\limits_{ {\bbd} \in \ZZ^d} 
\bbt^{\bbd}\cdot 
\II(\Omega({\bbd})), 
\\ 
\label{eq:idealizer-normal}
&=
\bigoplus\limits_{ {\bbd} \in \ZZ^d} 
\bbt^{\bbd}\cdot 
\left\<\prod\limits_{h_i({\bbd})<0}\df{h_i}{h_i(-{\bbd})-1}\right\>,
\end{align}
where the outer product runs over those $j=1, 2, \dots, m$ with $h_{j}({\bbd})<0$. 
\end{thm}

We now turn to an example that illustrates Theorem~\ref{thm:ST-dasr-3.2.2}, as well as our notation. 

\begin{ex}
\label{ex:exampleDO}
Let $R = \CC[s, st, st^2, st^3]$, where $\deg(s) = 1 = \deg(t)$.
Let $\partial_1$ and $\partial_2$ denote the partial derivatives with respect to $s$ and $t$, and let $\theta_1 = s \partial_1$ and $\theta_2 = t \partial_2$. 
By Theorem~\ref{thm:ST-dasr-3.2.2}, 
$D(R)_{\bbd}$, for some $\bbd = (d_1,d_2)  \in \ZZ^2$, is generated by an element in the form of 
$s^{d_1}t^{d_2}\cdot  f(\theta_1, \theta_2)$, 
where $ f(\theta_1, \theta_2)$ is a polynomial in $\theta_1$ and $\theta_2$. 

To understand the philosophy from \cite{Sai-Tr-DASR}, fix $\bbd = (d_1,d_2) \in \ZZ^2$, and consider  a differential operator on 
$R=k[s, st, st^2, st^3]$ of the form $\delta=s^{d_1}t^{d_2} f(\theta_1, \theta_2)$. By definition, $\delta \ast (s^{m_1}t^{m_2}) \in R$
for any monomial $s^{m_1}t^{m_2} \in R$.  We can check that
\[
\delta \ast s^{m_1}t^{m_2} 
= f(\bbm) s^{d_1+m_1}t^{d_2+m_2}. 
\]
In particular, if $s^{d_1 + m_1} t^{d_2 + m_2} \notin R$ then we must have that $\delta \ast s^{m_1}t^{m_2} = 0$ so that $f(\bbm) = 0$.
In summary, $\delta \ast (s^{m_1}t^{m_2})$ may be nonzero if and only if $s^{d_1 + m_1} t^{d_2 + m_2} \in R$.

These multidegrees $\bbm =(m_1, m_2)$ where $\delta \ast s^{m_1}t^{m_2}$ vanishes are exactly the points in $\Omega({\bbd})$.  We illustrate examples of $\Omega({\bbd})$ for specific $\bbd$ in Figure~\ref{fig:exampleDO-idealizer} by the integral points that are on the dotted lines in the first quadrant, including those on the positive horizontal axis. 

\begin{center}
\begin{figure}[b]
\begin{tabular}{ccc} 
\begin{tikzpicture}[scale=0.6]
 \draw (-3, 0) -- ( 3, 0) ; 
 \draw (0, -3) -- (0,6) ;
 \draw[blue, thick] (0,0) -- (1,3) ; 
 \draw[blue, thick] (0,0) -- (3,0) ;
 \filldraw (-1,2) circle (2pt);
 \filldraw (1,-2) circle (2pt);
 \coordinate [label=left:$\bbd$] (A) at (-1.15,2.15);
 \coordinate [label=right:$-\bbd$] (B) at (1.1,-2.1);
 \draw[thick] (1,-2) -- (2,1) ;  
 \draw[densely dotted, thick] (1.333, 0) -- (2, 2) ; 
 \draw[densely dotted, thick] (1, 0) -- (2, 3) ; 
 \draw[densely dotted, thick] (.666, 0) -- (2, 4) ;  
 \draw[densely dotted, thick] (.333, 0) -- (2, 5) ;  
 \draw[densely dotted, thick] (1, 3) -- (2, 6) ; 
\end{tikzpicture}

& 
\hspace{.8in}
&
\begin{tikzpicture}[scale=0.6]
 \draw (-3, 0) -- ( 2, 0) ; 
 \draw (0, -3) -- (0,6) ;
 \draw[blue, thick] (0,0) -- (1,3) ; 
 \draw[blue, thick] (0,0) -- (2.5,0) ;
 \filldraw (-2,-1) circle (2pt);
 \filldraw (2,1) circle (2pt);
 \coordinate [label=left:$\bbd$] (A) at (-2.15,-1.15);
 \coordinate [label=right:$-\bbd$] (B) at (2.1,1.3);
 \draw[thick] (2,1) -- (3,4) ;  
 \draw[thick] (2,1) -- (4,1) ; 
 \draw[densely dotted, thick] (1.333, 0) -- (3, 5) ;  
 \draw[densely dotted, thick] (1, 0) -- (3, 6) ; 
 \draw[densely dotted, thick] (.666, 0) -- (2.666, 6) ;   
 \draw[densely dotted, thick] (.333, 0) -- (2.333, 6) ;   
 \draw[densely dotted, thick] (1, 3) -- (2, 6)  ; 
 \draw[densely dotted, thick] (2.5, 0) -- (4,0) ;
\end{tikzpicture}
\end{tabular}
\caption{Half-lines of vanishing for various $\bbd$; here, $\bbd = (-1,2)$ and $(-2,-1)$. }
\label{fig:exampleDO-idealizer}
\end{figure}
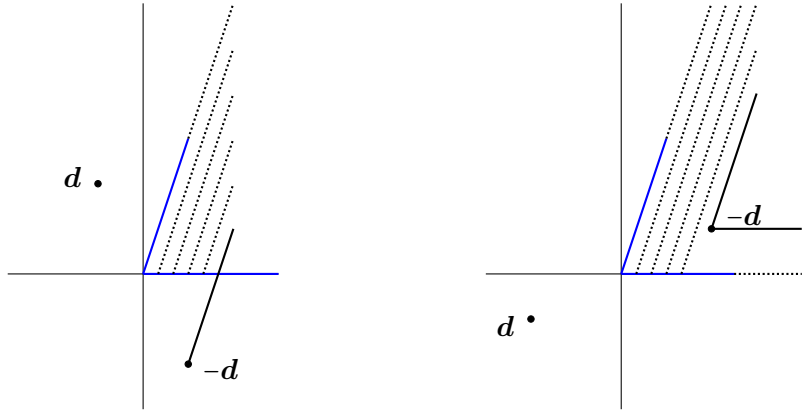
\end{center}

To explicitly compute the polynomial $f = f(\theta_1,\theta_2)$ for a fixed $\bbd\in\ZZ^2$, we will
use the primitive integral support functions for the two facets of $A$, which are $h_1=\theta_2$ and $h_2 = 3\theta_1 -\theta_2$. Both of these must divide $f$, along with all of the linear forms representing the dotted lines in Figure~\ref{fig:exampleDO-idealizer}. 
Specifically, $\II(\Omega(-1,2))$ is generated by 
\[
f(\theta_1, \theta_2) 
=   h_2 (h_2 -1) ( h_2 -2) ( h_2 -3) (h_2 - 4) 
=  (h_2, 4)!,
\]
and $\II(\Omega(-2,-1))$ is generated by 
\[
f(\theta_1, \theta_2) 
=  h_1 h_2 (h_2 -1) ( h_2 -2) ( h_2 -3) (h_2 - 4)  
=   (h_1, 0)! (h_2, 4)!.
\]
For $\bbd = (-1, 2)$, we have $h_2( -\bbd ) = 5$ while $h_1(- \bbd) = - 2 <0$, so there are no linear forms involving $h_1$ in this generator. 
Similarly for $\bbd = (-2, -1)$, 
$h_1( - \bbd) = 1$,
$h_2( -\bbd) = 5$. In particular,  
$(h_1, h_1(-\bbd) -1)! = (h_1, 0)! =h_1$ and $(h_2, h_2(-\bbd) -1)! = (h_2, 4)! =(h_2-4)(h_2-3)(h_2-2)(h_2-1)h_2$. 
\end{ex} 

As shown in Example~\ref{ex:exampleDO}, $\II(\Omega({\bbd}))$ is an ideal in $\CC[\theta]$, and any polynomial $f(\theta) \in \CC[\theta]$ has multidegree $\boldzero$; for any monomial ${\bbt}^{\bbm} \in R_A$, 
$f(\theta)\ast {\bbt}^{\bbm} = f(\bbm)\bbt^{\bbm}$, which belongs to $R_A$.  
So to determine if the $\bbd$-th graded piece of $D(R_A)$ applied to a monomial lands in $R_A$ or an $R_A$-ideal $J$, it is enough to test the membership of any monomial  $\bbt^{\bbd+\bbm}$ and then adjust the $\theta$-portion of the differential operator appropriately.  
Membership failure for $R_A$ only occurs if $\bbm$ lies in $\NN A$ but outside of the cone $-\bbd+\RR_{\geq 0}A$, whereas for $J \neq R_A$, failure is more likely to occur.  Since the monomials whose exponents lie on a face of $-\bbd+\RR_{\geq 0}A$ will lie on the corresponding face $\tau$ of $A$, if an associated prime of $J$ is a prime ideal associated to a face that contains $\tau$, then membership failure will also occur for $\bbm$ on this face.
In either case, there are only finitely many potential linear forms to be determined, and they are of the form 
$h_{\sigma}(\theta) - h_{\sigma}(-{\bbd})$ where $h_{\sigma}$ is the primary support function of a facet $\sigma$. 

In \cite{Er_SRDops,TrDM,Tripp}, Eriksson, Traves, and Tripp separately computed the ring of differential operators of a Stanley--Reisner ring over an arbitrary field, i.e., the quotient of any polynomial ring over a field by a squarefree monomial ideal.  
We include here the ring of differential operators of an ordinary double point $R=\CC[x,y]/\<x y\>$ using the viewpoint presented by the above authors, as it exhibits behavior akin to our computations in this article. 
In~\cite{Mussonzd}, Musson also considered the differential operators on an ordinary double point.

\begin{ex}
\label{ex:Traves}
The ring of differential operators on $\CC[x,y]$ is the Weyl algebra $W=\CC\{x,y,\partial_x, \partial_y\}$, which is the free associative algebra generated by $x,y,\partial_x,\partial_y$, subject to the relations: 
\[
\{ xy-yx, 
\partial_x\partial_y-\partial_y\partial_x,
\partial_x y-y\partial_x, 
\partial_y x-x\partial_y, 
\partial_x x-x\partial_x +1, 
\partial_y y-y\partial_y +1
\}.
\]
The Weyl algebra $W$ is a graded ring with $\deg(x)=\deg(y)=1=-\deg(\partial_x)=-\deg(\partial_y)$. 

For the ordinary double point ring  $R=\CC[x,y]/\<x y\>$, 
Traves shows in \cite[Theorem 3.5]{TrDM} that the idealizer of  $\<x y\>$ in $W$ is also graded and generated by 
\[
\{1, x^m, y^n, x^my^n, 
x^m\partial_x^i, y^n\partial_y^j, 
x^my^n \partial_x^i\partial_y^j 
\mid i \geq m>0, j \geq n > 0\}.
\]
Notice in particular that $x^m\partial_x^m$ and $y^n\partial_y^n$ both have degree 0. 
Setting $\theta_x=x\partial_x$ and $\theta_y=y\partial_y$, 
from the Weyl algebra relations, it follows that 
\[
\partial_x^i=x^{-i}\cdot \prod\limits_{\ell=0}^{i-1} (\theta_x-\ell)=x^{-i}\df{\theta_x}{i-1}
\quad\text{and}\quad 
\partial_y^j=y^{-j}\cdot \prod\limits_{\ell=0}^{j-1} (\theta_y-\ell)=\df{\theta_y}{j-1}.
\]
Hence $x^m\partial_x^i=x^{m-i}\cdot\df{\theta_x}{i-1}$ and $y^n\partial_y^j=y^{n-j}\cdot\df{\theta_y}{j-1}$, 
and $x^my^n\partial_x^i\partial_y^j$ has multidegree $(m-i,n-j)$ in $W$. 
In fact, 
\[
W = \bigoplus\limits_{(m,n) \in \ZZ^2}
x^my^n \cdot
\< \df{\theta_x}{m-1}\df{\theta_y}{n-1}\>,
\] 
which is a presentation of the Weyl algebra using the Saito--Traves approach from Theorem~\ref{thm:ST-dasr-3.2.2}.
From this viewpoint, Traves showed that  
\[
\II(\<x y\>)_{(m,n)} = 
\begin{cases}
x^my^n\cdot \CC[\theta_x,\theta_y] 
&\text{ if }m,n \geq 0,\\
x^my^n\cdot \< \df{\theta_x}{-m}\> 
&\text{ if } m<0, n\geq 0,\\
x^my^n \cdot \< \df{\theta_y}{-n}\> 
&\text{ if } n<0, m\geq 0,\\
x^my^n \cdot \< \df{\theta_x}{-m}\df{\theta_y}{-n}\> &\text{ if } m,n <0.
\end{cases}
\]
Further, $\<x y\>W$ can be expressed as a multigraded $W$-ideal that is contained in $\II(\<x y\>)$ as follows: 
\[
\<x y\>W_{(m,n)} = 
\begin{cases}
x^my^n\cdot \CC[\theta_x,\theta_y] &\text{ if }m,n > 0,\\
x^my^n \cdot \< \df{\theta_x}{-m} \> &\text{ if } m\leq 0, n> 0,\\
x^my^n \cdot \< \df{\theta_y}{-n}\> &\text{ if } n\leq 0, m> 0,\\
x^my^n \cdot \< \df{\theta_x}{-m}\df{\theta_y}{-n}\>  &\text{ if } m,n \leq 0.\\
\end{cases}
\]
Now, applying 
Proposition~\ref{prop:SmStDO-Prop1.6} yields a computation for the ring of differential operators for the ordinary double point $R=\CC[x,y]/\<xy\>$: 
\begin{equation} 
D\left(\frac{\CC[x,y]}{\<xy\>}\right)_{(m,n)}=\begin{cases}
0 
&\text{ if }m n \neq 0,\\
\displaystyle\frac{\CC[\theta_x,\theta_y]}{\< \theta_x \theta_y\> } 
&\text{ if } m=n=0,\\
x^m\cdot  \displaystyle\frac{\< \df{\theta_x}{-m} \>}{\< \df{\theta_x}{-m}\theta_y\> } 
&\text{ if } m \neq 0, n=0,\\
y^n \cdot \displaystyle\frac{\< \df{\theta_y}{-n} \>}{\< \df{\theta_y}{-m}\theta_x\> } 
&\text{ if } n \neq 0, m=0.\\
\end{cases} 
\label{eq:bigtriangleup}
\end{equation}
\end{ex}

We next fix some notation to be used in  the remainder of this paper. 
Whenever the dimension of the semigroup ring under consideration is $k=2$, instead of using $t_1,t_2$ as our variables, we will use $s,t$.  
Further, when considering subsets of $\RR^2$ and $\RR^3$ that contain the lattice points that describe a set of monomials in our semigroup, such as lines or planes, we will describe them with the variables $x, y,$ and $z$, 
for example, the line $y=2x-1$ in $\RR^2$ or the plane $x-z=2$ in $\RR^3$.

Consider the matrix
\[
A_n=\begin{bmatrix} 1 & 1 & 1 & \cdots & 1 \\
0 & 1 & 2 & \cdots & n \\
\end{bmatrix}.
\]  
We call $R_{A_n}=\CC [\NN A_n]=\CC[s,st,st^2, \ldots st^n]$ the ring of the {\em rational normal curve of degree $n$}, since it is the coordinate ring of the affine cone of the projective rational normal curve. 
This ring will be the subject of the next three sections. 
The two facets of $A_n$ are 
\begin{align*}
\sigma_1 
= \NN \begin{bmatrix} 1 \\ 0 \end{bmatrix}
= \{
(x,y)\in\NN^2\mid x\geq 0, y=0 
\}
\quad\text{and}\quad
\sigma_2 
= \NN\begin{bmatrix}1 \\ n\end{bmatrix}
= \{
(x,y)\in\NN^2\mid x\geq 0,  
y=n x\}. 
\end{align*}
The prime ideal associated to $\sigma_1$ is $P_1=\<st,st^2, \ldots, st^n\>$, and 
the prime ideal associated to $\sigma_2$ is $P_2=\<s, st, \ldots, st^{n-1}\>$. 
We will consider the radical ideal 
\[
J = P_1 \cap P_2 = 
\begin{cases}
\<st,st^2,\ldots, st^{n-1}\> &\text{if } n>1,\\
\<s^2t\> &\text{if } n=1.
\end{cases}
\]
Observe that $R_{A_n}/ J  \cong \CC[x, y]/\<x y\>$ for all $n>0$.  
At the end of Section~\ref{sec:higherD}, we revisit Example~\ref{ex:Traves} and present a $\CC$-algebra isomorphism between $D(\CC[x, y]/\<x y\>)$ and $D(R_{A_n}/J )$ and compare it to our calculations for $\II(J)/D(R_{A_n},J)$.

\section{Differential operators on the rational normal curve of degree 2}
\label{sec:2D}

In this section, we compute the idealizer, $\II(I)$, along with the subset of differential operators $D(R_{A_2},I)$ for the ideals $I=P_1=\<st,st^2\>$ and $I=J=\<st\>$ 
over the ring of the
rational normal curve of degree $2$, $R_{A_2}=\CC[s,st,st^2]$.
To aid our computations we include illustrations of the lattice representing the multidegrees in the plane broken down into four chambers where the operators will be determined by similar expressions. 
The facets of $A_2$ are  
\begin{align*}
\sigma_1 = 
\{(x,y) \in \NN^2 \mid x \geq 0, y=0\}
\quad\text{and}\quad
\sigma_2 = 
\{(x,y) \in \NN^2 \mid x, y \geq 0, y=2x\},
\end{align*} 
which have primitive integral support functions 
\[
h_1=\theta_2\quad \text{and}\quad h_2=2\theta_1-\theta_2.
\]
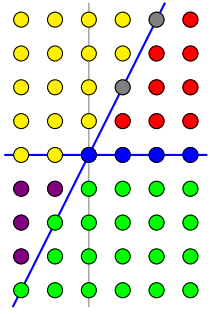
\begin{wrapfigure}[14]{l}[-5pt]{1.45in}
\vspace{-5pt}
\centering
\begin{tikzpicture}[scale=0.45]
\draw[blue, thick] (-2.5,0) -- (3.5,0);
\draw[gray] (0,-4.5) -- (0,4.5);
\draw[blue, thick] (-2.25,-4.5) -- (2.25,4.5);
\filldraw[violet] (-1,-1) circle (2pt);
\filldraw[violet] (-2,-1) circle (2pt);
\filldraw[violet] (-2,-2) circle (2pt);
\filldraw[violet] (-2,-3) circle (2pt);
\foreach \y in {-2,-3}{ \node[draw,circle,inner sep=2pt,violet,fill] at (-2,\y) {};}
\foreach \x in {-1,-2}{ \node[draw,circle,inner sep=2pt,violet,fill] at (\x,-1) {};}
\foreach \y in {2,3}{ \node[draw,circle,inner sep=2pt,red,fill] at (2,\y) {};}
\foreach \x in {1,2,3}{ \node[draw,circle,inner sep=2pt,red,fill] at (\x,1) {};}
\foreach \x in {1,2}{ \node[draw,circle,inner sep=2pt,gray,fill] at (\x,2*\x) {};}
\foreach \x in {0,1,2,3}{ \node[draw,circle,inner sep=2pt,blue,fill] at (\x,0) {};}
\foreach \y in {1,2,3,4}{ \node[draw,circle,inner sep=2pt,red,fill] at (3,\y) {};}
\foreach \x in {0,1,2,3}{
      \foreach \y in {-4,-3,...,-1}{
        \node[draw,circle,inner sep=2pt,green,fill] at (\x,\y) {};}}
      \foreach \y in {-4,-3,-2}{
        \node[draw,circle,inner sep=2pt,green,fill] at (-1,\y) {};}
\foreach \y in {-4}{ \node[draw,circle,inner sep=2pt,green,fill] at (-2,\y) {};}
\foreach \x in {-2,-1,0,1}{
      \foreach \y in {3,4}{
        \node[draw,circle,inner sep=2pt,yellow,fill] at (\x,\y) {};}}
\foreach \x in {-2,-1,0}{
      \foreach \y in {1,2}{
        \node[draw,circle,inner sep=2pt,yellow,fill] at (\x,\y) {};}}
\foreach \x in {-2,-1}{ \node[draw,circle,inner sep=2pt,yellow,fill] at (\x,0) {};}
\foreach \x in {-2,-1,...,3}{
      \foreach \y in {-4,-3,...,4}{
        \node[draw,circle,inner sep=2pt] at (\x,\y) {};}}
\end{tikzpicture}
\captionsetup{margin=0in,width=1.5in,font=small,labelsep=newline,justification=centering}
\caption{Chambers of $D(R_{A_2})$}
\label{fig:2D-lattice}
\end{wrapfigure}
Figure~\ref{fig:2D-lattice} illustrates the integer lattice, divided into four \emph{chambers} that are colored as follows: 
\hfill

\begin{itemize}
\item[\textbf{C1}:] The red multidegrees correspond to monomials in $J$, the gray multidegrees correspond to monomials in $P_1 \setminus J$ and the blue multidegrees correspond to monomials in $R_{A_2} \setminus P_1$,
\item[\textbf{C2}:] The yellow multidegrees are the $\bbd$ with $h_1({\bbd})\geq 0$ and $h_2 ({\bbd})< 0$, 
\item[\textbf{C3}:] The violet multidegrees are the $\bbd$ with both $h_1({\bbd})<0$ and $h_2({\bbd})<0$, and 
\item[\textbf{C4}:] The green multidegrees are the $\bbd$ with $h_1({\bbd})<0$ and $h_2 ({\bbd})\geq 0$. 
\end{itemize}

\smallskip 
Still following the convention $\df{h}{ n} =1$ if $n <0$, 
by Theorem~\ref{thm:ST-dasr-3.2.2}, the graded pieces of $D(R_{A_2})$ are 
\[ 
\hspace*{2in}
D(R_{A_2})_{\bbd} =
s^{d_1}t^{d_2} \cdot 
\left\<  \df{h_1}{ h_1(  {- \bbd} ) -1 }   \df{h_2}{ h_2(  {- \bbd} ) -1 } \right\>. 
\]
Broken down by chambers, this amounts to:
\[
D(R_{A_2})_{\bbd}=
\begin{cases}
s^{d_1}t^{d_2}\cdot \CC[\theta] 
&\text{if }\bbd\in\textbf{C1},\\
s^{d_1}t^{d_2} \cdot \< \df{h_2}{-2d_1+d_2-1}\> 
&\text{if }\bbd\in\textbf{C2},\\
s^{d_1}t^{d_2} \cdot \< \df{h_1}{-d_2-1} \df{h_2}{-2d_1+d_2-1}\> 
&\text{if }\bbd\in\textbf{C3},\\
s^{d_1}t^{d_2} \cdot\< \df{h_1}{-d_2-1}\> 
&\text{if }\bbd\in\textbf{C4}.\\
\end{cases}
\]  

\smallskip
\begin{ex}\label{ex:A_2:P1}
We first compute the graded pieces of the sets of differential operators $\II(P_1)$ and $D(R_{A_2},P_1)$. 
Later, since $\CC[x] \cong R_{A_2}/P_1$, we will exhibit a $\CC$-algebra isomorphism between $D(\CC[x])$ and $D(R_{A_2}/P_1)$.

To begin, recall that
\[
\II(P_1) =\{\delta \in D(R_{A_2}) \mid \delta * P_1 \subseteq P_1\} 
\quad\text{and}\quad 
D(R_{A_2},P_1) =\{\delta \in D(R_{A_2}) \mid \delta * R \subseteq P_1\}.
\] 
Now if $\bbd \in \textbf{C1}$ and $s^{m_1}t^{m_2} \in P_1$ 
or if  $\bbd \in \textbf{C1}\setminus \sigma_1$ and $s^{m_1}t^{m_2} \in R_{A_2}$, then for any $g(\theta) \in \CC[\theta]$, 
\[
s^{d_1}t^{d_2}\cdot g(\theta) * s^{m_1}t^{m_2}=g(\bbm) s^{d_1+m_1}t^{d_2+m_2} \in P_1, \text{ so} 
\]
\[
\II(P_1)_{\bbd}=D(R_{A_2})_{\bbd} \text{ for all } \bbd \in \textbf{C1}
\quad\text{and}\quad 
D(R_{A_2},P_1)_{\bbd}=D(R_{A_2})_{\bbd} \text{ for all } \bbd \in \textbf{C1}\setminus\sigma_1.
\]
\end{ex}
\begin{figure}[h]
\centering
\begin{subfigure}{1.55in}
 \centering
\begin{tikzpicture}[scale=0.6]
\draw[blue, thick] (0,0) -- (3.5,0);
\draw[gray] (0,0) -- (0,6.5);
\draw[red, thick] (0,0) -- (3.5,7);
\draw[red] (0.5,0) -- (3,5);
\draw[red]  (1,0) -- (3,4);
\draw[red] (1.5,0) -- (3,3);
\draw[red]  (2,0) -- (3,2);
\draw[red] (2.5,0) -- (3,1);
\foreach \y in {1,2}{ \node[draw,circle,inner sep=2pt,red,fill] at (1,\y) {};}
\foreach \x in {1,2,3}{ \node[draw,circle,inner sep=2pt,red,fill] at (\x,2*\x) {};}
\foreach \x in {0,1,2,3}{ \node[draw,circle,inner sep=2pt,blue,fill] at (\x,0) {};}
\foreach \y in {1,2,3}{ \node[draw,circle,inner sep=2pt,red,fill] at (2,\y) {};}
\foreach \y in {1,2,3,4,5}{ \node[draw,circle,inner sep=2pt,red,fill] at (3,\y) {};}
\end{tikzpicture}
\captionsetup{margin=0in,width=1.5in,font=small,justification=centering}
\caption{Lines parallel \\ to $\sigma_2$.}
\label{fig:2D:P_1:diagonal}
\end{subfigure}
\begin{subfigure}{1.55in}
 \centering
\begin{tikzpicture}[scale=0.6]
\draw[blue, thick] (0,0) -- (3.5,0);
\draw[gray] (0,0) -- (0,6.5);
\draw[red, thick] (0,0) -- (3.5,7);
\draw[red] (0.5,1) -- (3,1);
\draw[red]  (1,2) -- (3,2);
\draw[red] (1.5,3) -- (3,3);
\draw[red]  (2,4) -- (3,4);
\draw[red]  (2.5,5) -- (3,5);
\foreach \y in {1,2}{ \node[draw,circle,inner sep=2pt,red,fill] at (1,\y) {};}
\foreach \x in {1,2,3}{ \node[draw,circle,inner sep=2pt,red,fill] at (\x,2*\x) {};}
\foreach \x in {0,1,2,3}{ \node[draw,circle,inner sep=2pt,blue,fill] at (\x,0) {};}
\foreach \y in {1,2,3}{ \node[draw,circle,inner sep=2pt,red,fill] at (2,\y) {};}
\foreach \y in {1,2,3,4,5}{ \node[draw,circle,inner sep=2pt,red,fill] at (3,\y) {};}
\end{tikzpicture}
\captionsetup{margin=0in,width=1.65in,font=small,justification=centering}
\caption{Lines parallel \\ to $\sigma_1$.}
\label{fig:2D:P_1:horizontal}
\end{subfigure}
\caption{Lines parallel to the facets.}
\label{fig:2D:parallel}
\end{figure}
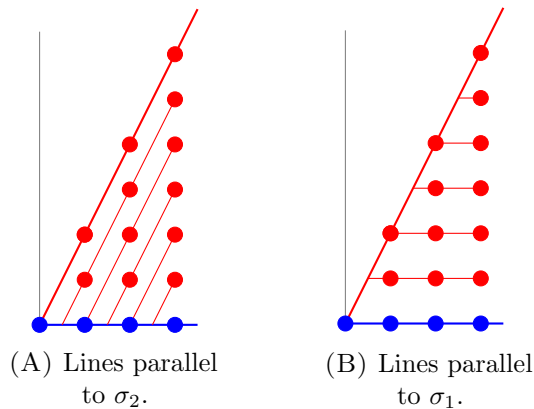

With the aid of Figure~\ref{fig:2D:parallel}, we will explain how to determine $\II(P_1)$ and $D(R_{A_2},P_1)$ in the other chambers. 
The red lattice points in Figure~\ref{fig:2D:parallel}, indicate the monomials in $P_1$.  

First, note that if 
$\bbd \in \textbf{C2}$ and $s^{m_1}t^{m_2} \in P_1$ 
or $\bbd \in \textbf{C2} \setminus (-\sigma_1)$ and $s^{m_1}t^{m_2} \in R_{A_2}$,  and  $\bbm$ lies on one of the lines $y=2x-r$ shown in Figure~\ref{fig:2D:P_1:diagonal}, 
then $s^{d_1}t^{d_2}\df{h_2}{-2d_1+d_2-1}$, the generator of $D(R_{A_2})_{\bbd}$, 
applied to such a monomial will either be a constant times a monomial represented by a red lattice point on one of the lines $y=2x-j$ for $0 \leq j \leq r$ or, when $\df{h_2(\bbm)}{-2d_1+d_2-1}=0$, it will be $0$. 
Hence, 
\[
\II(P_1)_{\bbd}=D(R_{A_2})_{\bbd} \text{ for } \bbd \in \textbf{C2}
\quad\text{and}\quad
D(R,P_1)_{\bbd}=D(R_{A_2})_{\bbd} \text{ for } \bbd \in \textbf{C2}\setminus (-\sigma_1).
\]

Now if $\bbd \in \textbf{C4}$ and $s^{m_1}t^{m_2} \in P_1$ or $\bbd \in \textbf{C4} \cup \sigma_1$ and $s^{m_1}t^{m_2} \in R_{A_2}$  and  $\bbm$ lies on one of the lines $y=r$ shown in Figure~\ref{fig:2D:P_1:horizontal}, then $s^{d_1}t^{d_2}\df{h_1}{-d_2-1}$, the generator of $D(R_{A_2})_{\bbd}$, applied to such a monomial will either be a constant times a monomial represented by a red lattice point on one of the lines $y=j$ for $ 0 \leq j \leq r$ or, when $\df{h_1(\bbm)}{-2d_1+d_2-1}=0$, it will be $0$.  
However, the monomials represented by the lattice points on $y=0$ do not lie in $P_1$.  
The monomials represented by the red lattice points on the line $y=-d_2$ are precisely the monomials whose image is a term represented by a lattice point on $\sigma_1$.  
Hence, we need to further right-multiply any operator in $D(R_{A_2})_{\bbd}$ by $\theta_2+d_2$, so that  
\begin{align*}
\II(P_1)_{\bbd}
=s^{d_1}t^{d_2}\cdot \< \df{h_1}{-d_2} \> \text{ for } \bbd \in \textbf{C4} 
\text{ \ and \ }
D(R_{A_2},P_1)
=s^{d_1}t^{d_2}\cdot \< \df{h_1}{-d_2} \> \text{ for } \bbd \in \textbf{C4} \cup \sigma_1.
\end{align*}

Using similar reasoning, we can easily see that 
\begin{align*}
\II(P_1)_{\bbd}
&=s^{d_1}t^{d_2}\cdot \< \df{h_1}{-d_2}\df{h_2}{-2d_1+d_2-1} \> \text{ for } \bbd \in \textbf{C3} 
\text{ and}\\
D(R_{A_2},P_1)_{\bbd}
&=s^{d_1}t^{d_2}\cdot \< \df{h_1}{-d_2}\df{h_2}{-2d_1+d_2-1} \> \text{ for } \bbd \in \textbf{C3} \cup (-\sigma_1).
\end{align*}

Putting these all together, the graded pieces of $\II(P_1)$ and $D(R_{A_2},P_1)$ are as follows: 
\[
\II(P_1)_{\bbd} =
\begin{cases}
s^{d_1}t^{d_2}\cdot \CC[\theta] &\text{if } {\bbd} \in \textbf{C1} = \NN A_2,\\
s^{d_1}t^{d_2}\cdot \< \df{h_2}{-2d_1+d_2-1} \> &\text{if } {\bbd} \in \textbf{C2},\\
s^{d_1}t^{d_2}\cdot \< \df{h_1}{-d_2} \df{h_2}{-2d_1+d_2-1}\> &\text{if } {\bbd} \in \textbf{C3},\\
s^{d_1}t^{d_2} \cdot \< \df{h_1}{-d_2} \> &\text{if } {\bbd} \in \textbf{C4},
\end{cases} \] 
\[
D(R_{A_2},P_1)_{\bbd} =
\begin{cases}
s^{d_1}t^{d_2} \cdot \CC[\theta] 
&\text{if } s^{d_1}t^{d_2} \in P_1,\\
s^{d_1}t^{d_2} \cdot\< \df{h_2}{-2d_1+d_2-1} \> 
&\text{if } {\bbd} \in \textbf{C2}\setminus (-\sigma_1),\\
s^{d_1}t^{d_2} \cdot\<
\df{h_1}{-d_2} \df{h_2}{-2d_1+d_2-1}\> &\text{if } {\bbd} \in \textbf{C3} \cup ( -\sigma_1),\\
s^{d_1}t^{d_2} \cdot\< \df{h_1}{-d_2} \> 
&\text{if } {\bbd} \in \textbf{C4} \cup \sigma_1.
\end{cases}
\] 
Now taking the quotient, we obtain:
\[ 
\left (\displaystyle\frac{\II(P_1)}{D(R_{A_2},P_1)} \right)_{\bbd} = 
\begin{cases}
0 
&\text{ if } {\bbd} \notin \ZZ \sigma_1, \\
s^{d_1}\cdot\dfrac{\CC[\theta]}{\< h_1\>}  
&\text{ if } {\bbd} \in \sigma_1,\\
s^{d_1}\cdot \dfrac{\< \df{h_2}{-2d_1-1} \>}{\< h_1\df{h_2}{-2d_1-1}\>  }
&\text{ if } {\bbd} \in (-\sigma_1 \setminus {\bf 0}).
\end{cases}
\]
Viewing $D(\CC [x])$ as a $\ZZ$-graded algebra over $\CC[\theta_x]$, we note that 
\[
D(\CC[x])= \sum\limits_{d=1}^{\infty} \partial_x^d \cdot \CC[\theta_x] \oplus \sum\limits_{d=0}^{\infty}x^d\CC[\theta_x]=\sum\limits_{d \in \ZZ} x^d \df{\theta_x}{-d-1} \cdot \CC[\theta_x].
\]

The map
\[
\phi\colon 
\sum\limits_{d \in \ZZ} x^d \cdot \df{\theta_x}{-d-1} \cdot \CC[\theta_x]
\rightarrow 
\sum\limits_{d \in \ZZ}s^{d} \cdot \dfrac{\<\df{h_2}{-2d-1}\>}{\<h_1\df{h_2}{-2d-1}\>},
\]
which is defined on the generators by 
\[
\phi(x^{d}\df{\theta_x}{-d-1})
=s^{d}\df{h_2}{-2d-1}+\<h_1\df{h_2}{-2d-1}\>, 
\]
although a $\CC$-vector space isomorphism, does not produce a ring isomorphism.   The graded pieces in negative degree are generated by polynomials in $\theta$ whose degrees are twice as large as large as the degrees in $\theta_x$ given in the Weyl algebra.  In fact, 
\[
D(R_{A_2}/P_1)=\begin{cases}
0 
&\text{ if } {\bbd} \notin \ZZ \sigma_1, \\
s^{d_1}\cdot\dfrac{\CC[\theta]}{\left\< \dfrac{h_1}{2}\right\>}  
&\text{ if } {\bbd} \in \sigma_1,\\
s^{d_1}\cdot \dfrac{\left\< \df{\dfrac{h_2}{2}}{-d_1-1} \right\>}{\left\< \dfrac{h_1}{2}\df{\dfrac{h_2}{2}}{-d_1-1}\right\>  }
&\text{ if } {\bbd} \in (-\sigma_1 \setminus {\bf 0}).
\end{cases}
\]
Therefore, we can produce an isomorphism of graded rings 
$\psi\colon D(\CC[x]) \rightarrow D(R_{A_2}/P_1)$ 
defined for any $m\in\NN$ by 
\[
\psi(x^{m})=s^m, 
\quad\text{and}\quad 
\psi(\partial_x^{m})=s^{-m} \cdot \df{\frac{h_2}{2}}{m-1}+\left\< \frac{h_1}{2}\df{\frac{h_2}{2}}{m-1}\right\>.
\]

\begin{ex}\label{ex:A_2:J}
We will now compute the graded pieces of the sets of differential operators $\II(J)$ and $D(R_{A_2},J)$, 
as well as the graded pieces of $JD(R_{A_2})$.
Recall that 
\[
\II(J) =\{\delta \in D(R_{A_2}) \mid \delta * J \subseteq J\} 
\quad\text{and}\quad 
D(R_{A_2},J) =\{\delta \in D(R_{A_2}) \mid \delta * R \subseteq J\}.
\] 
\end{ex}
\begin{wrapfigure}[15]{l}{1.5in}
\centering
\begin{tikzpicture}[scale=0.6]
\draw[cyan, thick] (0,0) -- (3.5,0);
\draw[gray] (0,0) -- (0,6.5);
\draw[blue, thick] (0,0) -- (3.5,7);
\draw[red] (0.5,0) -- (3,5);
\draw[cyan,thick]  (1,0) -- (3,4);
\foreach \y in {1,2}{ \node[draw,circle,inner sep=2pt,red,fill] at (1,\y) {};}
\foreach \x in {1,2,3}{ \node[draw,circle,inner sep=2pt,blue,fill] at (\x,2*\x) {};}
\foreach \x in {0,1,2,3}{ \node[draw,circle,inner sep=2pt,blue,fill] at (\x,0) {};}
\foreach \y in {1,2,3}{ \node[draw,circle,inner sep=2pt,red,fill] at (2,\y) {};}
\foreach \y in {1,2,3,4,5}{ \node[draw,circle,inner sep=2pt,red,fill] at (3,\y) {};}
\end{tikzpicture}
\captionsetup{margin=0in,width=1.5in,font=small,labelsep=newline,justification=centering}
\caption{Vanishing \\for $\bbd=(-1,0)$}
\label{fig:2D:J}
\end{wrapfigure}

We will soon give a general formula for the graded pieces of $\II(J)$ and $D(R_{A_2},J)$; 
however, for illustrative purposes, first consider the graded piece of $D(R_{A_2})$ at $(-1,0)$: 
$s^{-1}\cdot \< \df{h_2}{1} \>.$ 
Applying $s^{-1}\cdot  \df{h_2}{1} $  to a monomial whose exponent lies in the two parallel half-lines $\sigma_2$ and $y=2x-1$ in $\NN A$ will yield $0$, which certainly lives inside $J$.
However, when we let $s^{-1}\cdot \df{h_2}{1}$ act on a monomial whose exponent is a member of the  half-lines $y=2x - 2$ or $y=0$ lying inside $\NN A_2$, 
we obtain an integer multiple of a monomial whose exponent lies in one of the facets of $A_2$, $\sigma_2$ or $ \sigma_1$ respectively, and these are not in $J$.  
The remaining monomials in $J$ yield another element of $J$ when they are acted upon by any operator in $s^{-1}\cdot \< \df{h_2}{1}\>$ . 

In Figure~\ref{fig:2D:J}, the two light blue lines indicate the two half-lines representing the multidegrees of monomials in $R_{A_2}$ that, after application of an element in $D(R_{A_2})_{(-1,0)}$, fails to yield an element in $J$.  
To correct for this lack of membership in $J$ for the monomials on $y=2x-2$, every element of $D(R_{A_2})_{(-1,0)}$ should be multiplied by $(h_2-2)$;
applying $s^{-1}\cdot \df{h_2}{2}$ to these monomials 
yields $0$. 
The application of $s^{-1}\cdot \df{h_2}{2}$ to the remaining monomials in $J$  
will output a term inside $J$. 
Thus, 
\[
\II(J)_{(-1,0)}=s^{-1}\cdot \< \df{h_2}{2} \>.
\]
Similarly, for every operator $\delta\in
D(R_{A_2})_{(-1,0)}$, 
$\delta(h_2-2)h_1\ast s^{m_1}t^{m_2}=0$
for every $\bbm$ on $y=2x-2$ or in  $\sigma_1$.
Also, $\delta(h_2-2)h_1\ast s^{m_1}t^{m_2}\in J$ for all other $\bbm$ in $\NN A$, so 
\[
D(R_{A_2},J)_{(-1,0)} 
= s^{-1}\cdot \< h_1\df{h_2}{2} \>.
\]

In fact, using a similar argument applied to any $\bbd\in\textbf{C2}$ in the case of $\II(J)_{\bbd}$ and $\bbd \in \textbf{C2}\setminus (-\sigma_1)$ for $D(R_{A_2},J)_{\bbd}$, 
applying an operator in $D(R_{A_2})_{\bbd}$ to 
a monomial with exponent on the half-lines $y=2x-j$ for $0 \leq j < h_2({-\bbd})$ will give $0$; whereas, these operators applied to a monomial with multidegrees on $y=2x+h_2({\bbd})$ yields a constant multiple of a monomial with exponent in $\sigma_2$. 
Hence, right-multiplying 
$s^{d_1}t^{d_2}\cdot \df{h_2}{h_2({-\bbd})-1}$ by $h_2+h_2({\bbd})$ produces an operator that, when applied to monomials with multidegrees on the lines $y=2x+h_2({\bbd})$, becomes $0$, and 
\begin{align*}
\II(J)_{\bbd} 
&=s^{d_1}t^{d_2}\cdot \< \df{h_2}{-2d_1+d_2} \> 
\quad\text{for all $\bbd\in\textbf{C2}$}, \quad\text{and}
\\
D(R_{A_2},J)_{\bbd} 
&= s^{d_1}t^{d_2}\cdot \< \df{h_2}{-2d_1+d_2} \>
\quad\text{for all $\bbd\in\textbf{C2}\setminus (-\sigma_1)$}.
\end{align*} 
We will discuss the multidegrees  $\bbd\in\textbf{C2}\cap (-\sigma_1)$ momentarily, when we turn to \textbf{C3}. 

\begin{wrapfigure}[15]{l}[-5pt]{1.5in}
\centering
\begin{tikzpicture}[scale=0.6]
\draw[blue, thick] (0,0) -- (3.5,0);
\draw[gray] (0,0) -- (0,6.5);
\draw[cyan, thick] (0,0) -- (3.5,7);
\draw[red] (0.5,1) -- (3,1);
\draw[cyan,thick]  (1,2) -- (3,2);
\foreach \y in {1,2}{ \node[draw,circle,inner sep=2pt,red,fill] at (1,\y) {};}
\foreach \x in {1,2,3}{ \node[draw,circle,inner sep=2pt,blue,fill] at (\x,2*\x) {};}
\foreach \x in {0,1,2,3}{ \node[draw,circle,inner sep=2pt,blue,fill] at (\x,0) {};}
\foreach \y in {1,2,3}{ \node[draw,circle,inner sep=2pt,red,fill] at (2,\y) {};}
\foreach \y in {1,2,3,4,5}{ \node[draw,circle,inner sep=2pt,red,fill] at (3,\y) {};}
\end{tikzpicture}
\captionsetup{margin=0in,width=1.65in,font=small,labelsep=newline,justification=centering}
\caption{Vanishing \\ for $\bbd=(-1,-2)$}
\label{fig:2D:horizontal}
\end{wrapfigure}

Determining the graded piece of multidegree $\bbd$ for both $\II(J)$ in \textbf{C4} and $D(R_{A_2},J)$ in $\textbf{C4}\setminus (-\sigma_2)$ is quite similar to the arguments we used to determine $\II(J)_{\bbd}$ for $\bbd$ in \textbf{C2} and $D(R_{A_2},J)_{\bbd}$ for $\bbd$ in $\textbf{C2}\setminus (\sigma_1)$, respectively.  We will briefly describe $\II(J)_{(-1,-2)}$ and $D(R_{A_2},J)_{(-1,-2)}$ with the aid of Figure~\ref{fig:2D:horizontal} and then immediately describe the general case.
Recall that $D(R_{A_2})_{(-1,-2)}=s^{-1}t^{-2} \cdot \< \df{h_1}{1} \>$. 
Similar to the argument for degree $\bbd= (-1,0)$ above, the monomials corresponding to the multidegrees which lie on the two light blue lines ($y=2$ and $\sigma_2$) in Figure~\ref{fig:2D:horizontal} are the only exponents of monomials that fail to land inside $J$ after the application of 
$s^{-1}t^{-2} \cdot  \df{h_1}{1} $. 

To correct this deficiency, right-multiply by $(h_1-2)$ for $\II(J)_{(-1,-2)}$ and $(h_1-2)h_2$ for $D(R_{A_2},J)_{(-1,-2)}$. Notice that applying the operator $s^{-1}t^{-2}\cdot \df{h_1}{2}$ to a monomial corresponding to 
$\bbd\in\NN A_2$ along the half-lines $y=2$ or the operator $s^{-1}t^{-2}\cdot \df{h_1}{2}h_2$ to a monomial corresponding to $\bbd \in \NN A_2$ along $y=2$ or $y=2x$ now yields $0$, and no problems are created for the remaining monomials in $J$ or $R_{A_2}$, respectively. 
Thus, 
\begin{equation*}
\II(J)_{(-1,-2)} =s^{-1}t^{-2}\cdot\< \df{h_1}{2} \> \quad \text{and} \quad
D(R_{A_2},J)_{(-1,-2)} 
=s^{-1}t^{-2}\cdot\< \df{h_1}{2}h_2 \>.
\end{equation*}
In fact,
\begin{align*}
\II(J)_{\bbd} 
&= s^{d_1}t^{d_2}\cdot \< \df{h_1}{-d_2} \>
\quad \text{for all }{\bbd}\in\textbf{C4}, 
\quad \text{and}
\\
D(R_{A_2},J)_{\bbd} 
&= s^{d_1}t^{d_2}\cdot \< \df{h_1}{-d_2} \>
\quad\text{for all }{\bbd}\in\textbf{C4}\setminus (-\sigma_2).
\end{align*} 
We will return to $D(R_{A_2},J)_{\bbd}$ for $\bbd \in \textbf{C4} \cap (-\sigma_2)$ when we turn to \textbf{C3}.

\begin{wrapfigure}[16]{l}[-5pt]{1.5in}
\centering
\begin{tikzpicture}[scale=0.6]
\draw[blue, thick] (0,0) -- (3.5,0);
\draw[gray] (0,0) -- (0,6.5);
\draw[blue, thick] (0,0) -- (3.5,7);
\draw[cyan, thick] (0.5,1) -- (3,1);
\draw[cyan,thick]  (0.5,0) -- (3,5);
\foreach \y in {1,2}{ \node[draw,circle,inner sep=2pt,red,fill] at (1,\y) {};}
\foreach \x in {1,2,3}{ \node[draw,circle,inner sep=2pt,blue,fill] at (\x,2*\x) {};}
\foreach \x in {0,1,2,3}{ \node[draw,circle,inner sep=2pt,blue,fill] at (\x,0) {};}
\foreach \y in {1,2,3}{ \node[draw,circle,inner sep=2pt,red,fill] at (2,\y) {};}
\foreach \y in {1,2,3,4,5}{ \node[draw,circle,inner sep=2pt,red,fill] at (3,\y) {};}
\end{tikzpicture}
\captionsetup{margin=0in,width=1.65in,font=small,labelsep=newline,justification=centering}
\caption{Vanishing \\for $\bbd=(-1,-1)$}
\label{fig:2D:1-1}
\end{wrapfigure}

Determining $\II(J)_{\bbd}$ for $\bbd \in \textbf{C3}$ or $D(R_{A_2},J)_{\bbd}$ for $\bbd \in \textbf{C3} \cup (-\sigma_1) \cup  (-\sigma_2)$ again is akin to the arguments we gave above for $D(R_{A_2},J)_{(-1,0)}$ and $D(R_{A_2},J)_{(-1,-2)}$.  
With the aid of Figure~\ref{fig:2D:1-1}, we will briefly describe 
\[
\hspace*{2in}
\II(J)_{(-1,-1)}
\quad\text{and}\quad
D(R_{A_2},J)_{(-1,-1)}.
\]
This explanation can easily be extended from $\bbd=(-1,-1)$ to all multidegrees $\bbd$ in \textbf{C3} (or in $\textbf{C3}\cup(-\sigma_1)\cup(-\sigma_2)$ in the case of $D(R_{A_2},J)$). 
If the operator 
$s^{-1}t^{-1}\cdot h_1h_2$, 
which is the generator for $D(R_{A_2})_{(-1,-1)}$,  is applied to any of the monomials corresponding to multidegree $\bbd\in\NN A_2$ 
along the two light blue half-lines in Figure~\ref{fig:2D:1-1} (the portion of $y=1$ or $y=2x-1$ in $\textbf{C1}$), 
we obtain an integer multiple of a monomial with exponent in the facets $ \sigma_1$ or $\sigma_2$, respectively, which is not in $J$,
and no problems are created for the remaining monomials in $J$ (or $R_{A_2}$ for $D(R_{A_2},J)$). 

Hence, right-multiplying by $(h_1-1)(h_2-1)$ yields a new operator 
$s^{-1}t^{-1}\cdot \df{h_1}{1} \df{h_2}{1}$ that will send to $0$ all monomials with multidegrees $\bbd\in\NN A_2$ along the half-lines $y=2x-1$ and $y=1$.  
No problems are created for the remaining monomials in $R_{A_2}$ and we obtain  
\[
\II(J)_{(-1,-1)} 
= D(R_{A_2},J)_{(-1,-1)} 
= s^{-1}t^{-1}\cdot \< \df{h_1}{1} \df{h_2}{1} \>.
\] 
In fact,
\begin{align*}
\II(J)_{\bbd} = D(R_{A_2},J)_{\bbd}
&=s^{d_1}t^{d_2}\cdot \< \df{h_1}{-d_2)} \df{h_2}{-2d_1+d_2}\>
\quad\text{ for all }\bbd\in
\textbf{C3}, 
\end{align*}
and 
\begin{align*}
D(R_{A_2},J)_{\bbd} 
&= s^{d_1}t^{d_2}\cdot \< \df{h_1}{-d_2)} \df{h_2}{-2d_1+d_2}\>
\quad\text{ for all }\bbd\in
(-\sigma_1)\cup(-\sigma_2).
\end{align*}
Hence, the graded pieces of $\II(J)$ and $D(R_{A_2},J)$ are as follows:
\[
\II(J)_{\bbd} =
\begin{cases}
s^{d_1}t^{d_2}\cdot \CC[\theta] &\text{if } {\bbd} \in \textbf{C1} = \NN A_2,\\
s^{d_1}t^{d_2}\cdot \< \df{h_2}{-2d_1+d_2} \> &\text{if } {\bbd} \in \textbf{C2},\\
s^{d_1}t^{d_2}\cdot \< \df{h_1}{-d_2} \df{h_2}{-2d_1+d_2}\> &\text{if } {\bbd} \in \textbf{C3},\\
s^{d_1}t^{d_2} \cdot \< \df{h_1}{-d_2} \> &\text{if } {\bbd} \in \textbf{C4},
\end{cases} 
\] 
\[
D(R_{A_2},J)_{\bbd} =
\begin{cases}
s^{d_1}t^{d_2} \cdot \CC[\theta] 
&\text{if } s^{d_1}t^{d_2} \in J,\\
s^{d_1}t^{d_2} \cdot\< \df{h_2}{-2d_1+d_2} \> 
&\text{if } {\bbd} \in (\textbf{C2}\setminus (-\sigma_1)) \cup (\sigma_2  \setminus \{\boldzero\}),\\
s^{d_1}t^{d_2} \cdot\<
\df{h_1}{-d_2} \df{h_2}{-2d_1+d_2}\> &\text{if } {\bbd} \in \textbf{C3} \cup ( -\sigma_1)\cup (-\sigma_2),\\
s^{d_1}t^{d_2} \cdot\< \df{h_1}{-d_2} \> 
&\text{if } {\bbd} \in (\textbf{C4}\setminus (-\sigma_2)) \cup (\sigma_1 \setminus \{\boldzero\}).
\end{cases}
\] 
Now taking the quotient, we obtain:
\[ 
\left (\displaystyle\frac{\II(J)}{D(R_{A_2},J)} \right)_{\bbd} = 
\begin{cases}
0 
&\text{ if } {\bbd} \notin (\ZZ \sigma_1 \cup \ZZ \sigma_2), \\
s^{d_1}t^{d_2}\cdot \dfrac{\CC[\theta]}{\< h_1h_2\>}
&\text{ if } {\bbd=\boldzero}, \\
s^{d_1}t^{d_2}\cdot\dfrac{\CC[\theta]}{\< h_1\>}  
&\text{ if } {\bbd} \in \sigma_1 \setminus \{\boldzero\},\\
s^{d_1}t^{d_2}\cdot\dfrac{\CC[\theta]}{\< h_2\>  }
&\text{ if } {\bbd} \in \sigma_2  \setminus \{\boldzero\},\\
s^{d_1}t^{d_2}\cdot \dfrac{\< \df{h_2}{-2d_1} \>}{\< h_1\df{h_2}{-2d_1}\>  }
&\text{ if } {\bbd} \in (-\sigma_1), \\
s^{d_1}t^{d_2}\cdot \dfrac{\< \df{h_1}{-d_2}\>}{\< h_2 \df{h_1}{-d_2} \>  }
&\text{ if } {\bbd} \in (-\sigma_2).
\end{cases}
\]

As both $J$ and $D(R_{A_2})$ are graded,  we can similarly determine $JD(R_{A_2})$.
Our goal for the remainder of the section is to compute the graded pieces of $JD(R_{A_2})$, in order to observe that, in this case,  $D(R_{A_2},J)=JD(R_{A_2})$.

To begin, note that for all $s^{d_1}t^{d_2} \in \ZZ^2$, 
\[
D(R_{A_2})_{(d_1-1,d_2-1)} 
= s^{d_1-1}t^{d_2-1} \cdot\< \df{h_1}{-d_2}\df{h_2}{-2d_1+d_2}\>.
\]
For $s^{m_1}t^{m_2} \in J$ the graded piece at the multidegree $(d_1-m_1,d_2-m_2)$ will be 
\[
s^{d_1-m_1}t^{d_2-m_2}\cdot
\< \df{h_1}{-(d_2-m_2)-1}\df{h_2}{-2(d_1-m_1)+d_2-m_2-1}\>.
\]  
Note that since $m_1$ and $m_2 $ are both positive,  
\[
s^{m_1}t^{m_2}s^{d_1-m_1}t^{d_2-m_2} = s^{d_1}t^{d_2}, 
\] 
and
\[
\< \df{h_1}{-(d_2-c_2)-1}\df{h_2}{-2(d_1-c_1)+d_2-c_2-1}\> 
\subseteq \< \df{h_1}{-d_2}\df{h_2}{-2d_1+d_2}\>.
\]  
Thus, to determine $(JD(R_{A_2}))_{\bbd}$ for any $\bbd\in\ZZ^2$, it is enough to (left) multiply $D(R_{A_2})_{(d_1-1,d_2-1)}$ by $st$.  
Hence, from our previous computation, it follows that
\[
(JD(R_{A_2}))_{\bbd} = 
\begin{cases}
s^{d_1}t^{d_2}\cdot \CC[\theta] 
&\text{ if } s^{d_1}t^{d_2} \in J,\\
s^{d_1}t^{d_2}\cdot \< \df{h_2}{-2d_1+d_2} \> 
&\text{ if }{\bbd} \in (\textbf{C2}\setminus (-\sigma_1)) \cup (\sigma_1 \setminus \{\boldzero\}),\\
s^{d_1}t^{d_2}\cdot \< \df{h_1}{-d_2}\df{h_2}{-2d_1+d_2} \> 
&\text{ if }{\bbd} \in \textbf{C3}\cup (-\sigma_1) \cup (-\sigma_2),\\
s^{d_1}t^{d_2}\cdot \< \df{h_1}{-d_2}\> 
&\text{ if }{\bbd} \in (\textbf{C4} \setminus (-\sigma_2)) \cup (\sigma_1 \setminus \{\boldzero\}).
\end{cases}
\]
Combining all cases together, it follows that for an arbitrary $\bbd\in\ZZ^2$, there is an equality 
\[ 
(JD(R_{A_2}))_{\bbd} 
= s^{d_1}t^{d_2}\cdot  \left\< \df{h_1}{h_1( - {\bbd} )}  \df{h_2}{h_2( - {\bbd} )} \right\>
= D(R_{A_2},J)_{\bbd}. 
\]

\section{Differential operators on the rational normal curve of degree 3}
\label{sec:3D}

In this section, we determine some
subsets of the ring of differential operators 
for the ring of the rational normal curve of degree $3$ determined by its interior ideal $J$.  
We contrast these computations with the degree $2$ case that was determined in Section~\ref{sec:2D}. 
Although the description follows the same reasoning as the degree $2$ setting, we ultimately get quite different behavior for the operators that end up in $JD(R_{A_3})$. 

\begin{wrapfigure}[19]{l}{1.60in}
\centering
\vspace{5pt}
\begin{tikzpicture}[scale=0.50]
\draw[blue, thick] (-2.5,0) -- (3.5,0);
\draw[gray] (0,-6.5) -- (0,6.5);
\draw[blue, thick] (-2.33,-7) -- (2.33,7);
\foreach \y in {1,2}{ \node[draw,circle,inner sep=2pt,red,fill] at (1,\y) {};}
\foreach \x in {1,2}{ \node[draw,circle,inner sep=2pt,blue,fill] at (\x,3*\x) {};}
\foreach \x in {0,1,2,3}{ \node[draw,circle,inner sep=2pt,blue,fill] at (\x,0) {};}
\foreach \y in {1,2,...,7}{ \node[draw,circle,inner sep=2pt,red,fill] at (3,\y) {};}
\foreach \y in {1,2,...,5}{ \node[draw,circle,inner sep=2pt,red,fill] at (2,\y) {};}
\foreach \x in {0,1,2,3}{
      \foreach \y in {-7,-6,...,-1}{
        \node[draw,circle,inner sep=2pt,green,fill] at (\x,\y) {};}}
      \foreach \y in {-7,-6,...,-3}{
        \node[draw,circle,inner sep=2pt,green,fill] at (-1,\y) {};}
      \foreach \y in {-7,-6}{
        \node[draw,circle,inner sep=2pt,green,fill] at (-2,\y) {};}
\foreach \x in {-2,-1,0,1}{
      \foreach \y in {4,5,6}{
        \node[draw,circle,inner sep=2pt,yellow,fill] at (\x,\y) {};}}
\foreach \x in {-2,-1,0}{
      \foreach \y in {1,2,3}{
        \node[draw,circle,inner sep=2pt,yellow,fill] at (\x,\y) {};}}
     \foreach \x in {-2,-1,...,2}{ \node[draw,circle,inner sep=2pt,yellow,fill] at (\x,7) {};}
\foreach \x in {-1,-2}{ \node[draw,circle,inner sep=2pt,yellow,fill] at (\x,0) {};}
\foreach \y in {-5,-4,...,-1}{ \node[draw,circle,inner sep=2pt,violet,fill] at (-2,\y) {};}
\foreach \y in {-2,-1}{ \node[draw,circle,inner sep=2pt,violet,fill] at (-1,\y) {};}
\foreach \x in {-2,-1,...,3}{
      \foreach \y in {-7,-6,...,7}{
        \node[draw,circle,inner sep=2pt] at (\x,\y) {};}}
\end{tikzpicture}
\captionsetup{margin=0in,width=1.5in,font=small,labelsep=newline,justification=centering}
\caption{Chambers of $D(R_{A_3})$}
\label{fig:3D-lattice}
\end{wrapfigure}
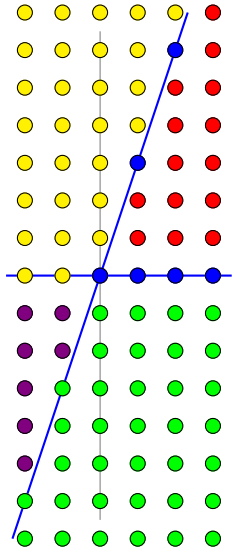

The ring of the rational normal curve of degree $3$ is $R_{A_3}=\CC[s,st,st^2,st^3]$, 
and we will compute $\II(J)/D(R_{A_3},J)$, 
where $J=\<st,st^2\>$ is a radical ideal. 
The facets of $A_3$ are
\begin{align*}
\sigma_1 = 
\{(x,y) \in \NN^2 \mid x \geq 0, y=0\}
\quad\text{and}\quad
\sigma_2 = 
\{(x,y) \in \NN^2 \mid x, y \geq 0, y=3x\},
\end{align*} 
which have primitive integral support functions 
\[
h_1=\theta_2 \quad \text{and}\quad h_2=3\theta_1-\theta_2.
\] 
Figure~\ref{fig:3D-lattice} illustrates the integer lattice, divided into four chambers that are colored as follows: 
\hfill

\vspace*{-7mm}
\medskip
\begin{itemize}
\item[\textbf{C1}:] The red multidegrees correspond to monomials in $J$, and the blue multidegrees correspond to monomials in $R_{A_2} \setminus J$,
\item[\textbf{C2}:] The yellow multidegrees are the $\bbd$ with $h_1({\bbd})\geq 0$ and $h_2 ({\bbd})< 0$, 
\item[\textbf{C3}:] The violet multidegrees are the $\bbd$ with both $h_1({\bbd})<0$ and $h_2({\bbd})<0$, and 
\item[\textbf{C4}:] The green multidegrees are the $\bbd$ with $h_1({\bbd})<0$ and $h_2 ({\bbd})\geq 0$. 
\end{itemize}
\medskip

Still following the convention $\df{h}{n}=1$ if $n<0$, by Theorem \ref{thm:ST-dasr-3.2.2}, the graded pieces of $D(R_{A_3})$ are
\[ 
\hspace*{2in}
D(R_{A_3})_{\bbd} =
s^{d_1}t^{d_2} \cdot 
\left\<  \df{h_1}{ h_1(  {- \bbd} ) -1 }   \df{h_2}{ h_2(  {- \bbd} ) -1 } \right\>. 
\]

\noindent Broken down by chambers, this amounts to:
\[
\hspace*{2in}
D(R_{A_3})_{\bbd}=
\begin{cases}
s^{d_1}t^{d_2}\cdot \CC[\theta] 
&\text{if }\bbd\in\textbf{C1},\\
s^{d_1}t^{d_2} \cdot \< \df{h_2}{-3d_1+d_2-1}\> 
&\text{if }\bbd\in\textbf{C2},\\
s^{d_1}t^{d_2} \cdot \< \df{h_1}{-d_2-1} \df{h_2}{-3d_1+d_2-1}\> 
&\text{if }\bbd\in\textbf{C3},\\
s^{d_1}t^{d_2} \cdot\< \df{h_1}{-d_2-1}\> 
&\text{if }\bbd\in\textbf{C4}.\\
\end{cases}
\]  

Determining the graded pieces of $\II(J)$ and $D(R_{A_3},J)$ in $D(R_{A_3})$ is very similar to  our computations in Section \ref{sec:2D}. 
Here we include some visualizations for $\bbd\in\{(-1,0), (-1,-3), (-1,-1), (-1,-2)\}$ in Figure~\ref{fig:3Dvanish} to aid our in our description of how to obtain $\II(J)_{\bbd}$ and $D(R_{A_3},J)_{\bbd}$. 
Then, we will list the expressions for $\II(J)_{\bbd}$ and $D(R_{A_3},J)_{\bbd}$ by chamber, as we did in Section \ref{sec:2D}.

\begin{figure}[h]
\centering
\begin{subfigure}{1.65in}
  \centering
 \begin{tikzpicture}[scale=0.60]
\draw[cyan, thick] (0,0) -- (3.5,0);
\draw[gray] (0,0) -- (0,6.5);
\draw[blue, thick] (0,0) -- (2.33,7);
\draw[red] (0.33,0) -- (2.66,7);
\draw[red] (0.66,0) -- (3,7);
\draw[cyan,thick] (1,0) -- (3,6);
\foreach \y in {1,2}{ \node[draw,circle,inner sep=2pt,red,fill] at (1,\y) {};}
\foreach \x in {1,2}{ \node[draw,circle,inner sep=2pt,blue,fill] at (\x,3*\x) {};}
\foreach \x in {0,1,2,3}{ \node[draw,circle,inner sep=2pt,blue,fill] at (\x,0) {};}
\foreach \y in {1,2,...,7}{ \node[draw,circle,inner sep=2pt,red,fill] at (3,\y) {};}
\foreach \y in {1,2,...,5}{ \node[draw,circle,inner sep=2pt,red,fill] at (2,\y) {};}
\end{tikzpicture}
\captionsetup{margin=0in,width=1.5in,font=small,justification=centering}
\caption{$\bbd=(-1,0)$.}
\label{fig:3D:slant}
\end{subfigure}%
\begin{subfigure}{1.65in}
\centering
\begin{tikzpicture}[scale=0.60]
\draw[blue, thick] (0,0) -- (3.5,0);
\draw[gray] (0,0) -- (0,6.5);
\draw[cyan, thick] (0,0) -- (2.33,7);
\draw[red] (0.33,1) -- (3.5,1);
\draw[red] (0.66,2) -- (3.5,2);
\draw[cyan,thick] (1,3) -- (3.5,3);
\foreach \y in {1,2}{ \node[draw,circle,inner sep=2pt,red,fill] at (1,\y) {};}
\foreach \x in {1,2}{ \node[draw,circle,inner sep=2pt,blue,fill] at (\x,3*\x) {};}
\foreach \x in {0,1,2,3}{ \node[draw,circle,inner sep=2pt,blue,fill] at (\x,0) {};}
\foreach \y in {1,2,...,7}{ \node[draw,circle,inner sep=2pt,red,fill] at (3,\y) {};}
\foreach \y in {1,2,...,5}{ \node[draw,circle,inner sep=2pt,red,fill] at (2,\y) {};}
\end{tikzpicture}
\captionsetup{margin=0in,width=1.5in,font=small,justification=centering}
\caption{$\bbd=(-1,-3)$.}
\label{fig:3D:horizontal}
\end{subfigure}%
\begin{subfigure}{1.65in}
\centering
\begin{tikzpicture}[scale=0.60]
\draw[blue, thick] (0,0) -- (3.5,0);
\draw[gray] (0,0) -- (0,6.5);
\draw[blue, thick] (0,0) -- (2.33,7);
\draw[red] (0.33,0) -- (2.66,7);
\draw[cyan,thick] (0.66,0) -- (3,7);
\draw[cyan,thick] (0.33,1) -- (3.5,1);
\foreach \y in {1,2}{ \node[draw,circle,inner sep=2pt,red,fill] at (1,\y) {};}
\foreach \x in {1,2}{ \node[draw,circle,inner sep=2pt,blue,fill] at (\x,3*\x) {};}
\foreach \x in {0,1,2,3}{ \node[draw,circle,inner sep=2pt,blue,fill] at (\x,0) {};}
\foreach \y in {1,2,...,7}{ \node[draw,circle,inner sep=2pt,red,fill] at (3,\y) {};}
\foreach \y in {1,2,...,5}{ \node[draw,circle,inner sep=2pt,red,fill] at (2,\y) {};}
\end{tikzpicture}
\captionsetup{margin=0in,width=1.5in,font=small,justification=centering}
\caption{$\bbd=(-1,-1)$.}
\label{fig:3D:1-1}
\end{subfigure}%
\begin{subfigure}{1.65in}
\centering
\begin{tikzpicture}[scale=0.60]
\draw[blue, thick] (0,0) -- (3.5,0);
\draw[gray] (0,0) -- (0,6.5);
\draw[blue, thick] (0,0) -- (2.33,7);
\draw[red] (0.33,1) -- (3.5,1);
\draw[cyan,thick] (0.3,0) -- (2.66,7);
\draw[cyan,thick] (0.66,2) -- (3.5,2);
\foreach \y in {1,2}{ \node[draw,circle,inner sep=2pt,red,fill] at (1,\y) {};}
\foreach \x in {1,2}{ \node[draw,circle,inner sep=2pt,blue,fill] at (\x,3*\x) {};}
\foreach \x in {0,1,2,3}{ \node[draw,circle,inner sep=2pt,blue,fill] at (\x,0) {};}
\foreach \y in {1,2,...,7}{ \node[draw,circle,inner sep=2pt,red,fill] at (3,\y) {};}
\foreach \y in {1,2,...,5}{ \node[draw,circle,inner sep=2pt,red,fill] at (2,\y) {};}
\end{tikzpicture}
\captionsetup{margin=0in,width=1.5in,font=small,justification=centering}
\caption{$\bbd=(-1,-2)$.}
\label{fig:3D:1-2}
\end{subfigure}%
\captionsetup{margin=0in,width=4.5in,
justification=centering}
\caption{Vanishing at various $\bbd$}
\label{fig:3Dvanish}
\end{figure}

For each $\bbd$ in the illustrations in Figure~\ref{fig:3Dvanish}, when we apply elements of  $D(R_{A_3})_{\bbd}$ to any of the monomials represented by lattice points along the light blue lines, we obtain an integer multiple of a monomial whose exponent lies on a facet.  As in Section \ref{sec:2D}, these lines determine the linear multiples $h_i-h_i(\bbd)$ we must append to $\II(\Omega(\bbd))$ to obtain either $\II(J)_{\bbd}$ or $D(R_{A_3},J)$.  Note that the light blue lines $h_i(\bbx)-h_i(\bbd)=0$ that have contain a red dot determine the $h_i-h_i(\bbd)$ we right-multiply by to obtain $\II(J)_{\bbd}$ and that the light blue lines $h_i(\bbx)-h_i(\bbd)=0$ that have non-empty intersection with the the entire cone determine the $h_i-h_i(\bbd)$ we right-multiply by to obtain $D(R_{A_3},J)_{\bbd}$.   
Below we summarize the graded pieces of $\II(J)$ and $D(R_{A_3},J)$, as we did in Section \ref{sec:2D}:
\[
\II(J)_{\bbd} = 
\begin{cases}
s^{d_1}t^{d_2}\cdot \CC[\theta] 
&\text{if } {\bbd} \in \NN A_3,\\
s^{d_1}t^{d_2} \cdot \< \df{h_1}{-d_2} \> 
&\text{if } {\bbd} \in \textbf{C4},\\
s^{d_1}t^{d_2} \cdot \< \df{h_2}{-3d_1+d_2} \> 
&\text{if } {\bbd} \in \textbf{C2},\\
s^{d_1}t^{d_2} \cdot \< \df{h_1}{-d_2} \df{h_2}{-3d_1+d_2}\>
&\text{if } {\bbd} \in \textbf{C3},
\end{cases}
\]   
\begin{align}
\label{eq:3D:DRJ}
D(R_{A_3},J)_{\bbd} = 
\begin{cases}
s^{d_1}t^{d_2}\cdot \CC[\theta] 
&\text{if } s^{d_1}t^{d_2} \in J,\\
s^{d_1}t^{d_2}\cdot \< \df{h_1}{-d_2} \>
&\text{if } {\bbd} \in (\textbf{C4}\setminus (-\sigma_2)) \cup (\sigma_1  \setminus \{\boldzero\}), \\
s^{d_1}t^{d_2}\cdot \< \df{h_2}{-3d_1+d_2} \>  
&\text{if } {\bbd} \in (\textbf{C2}\setminus -\sigma_1) \cup (\sigma_2 \setminus \{\boldzero\}), \\
s^{d_1}t^{d_2}\cdot \< \df{h_1}{-d_2}\df{h_2}{-3d_1+d_2} \> 
&\text{if } {\bbd} \in \textbf{C3}\cup(-\sigma_1) \cup (-\sigma_2). 
\end{cases}
\end{align}

Taking the quotient, 
we obtain:
\[
\left (\displaystyle\frac{\II(J)}{D(R_{A_3},J)} \right )_{\bbd}= 
\begin{cases}
0 
&\text{if } {\bbd} \notin (\ZZ \sigma_1 \cup \ZZ \sigma_2),\\
s^{d_1}t^{d_2}\cdot \dfrac{\CC[\theta]}{\< h_1h_2\>}  
&\text{if } {\bbd=\boldzero}, \\
s^{d_1}t^{d_2}\cdot \dfrac{\CC[\theta]}{\< h_1\>  }
&\text{if } {\bbd} \in ( \sigma_1 \setminus \{\boldzero\}),\\
s^{d_1}t^{d_2}\cdot \dfrac{\CC[\theta]}{\< h_2\>  }
&\text{if } {\bbd} \in ( \sigma_2 \setminus \{\boldzero\}),\\
s^{d_1}t^{d_2}\cdot \dfrac{\< \df{h_2}{-3d_1} \>}{\< h_1\df{h_2}{-3d_1}\>  }
&\text{if } {\bbd} \in (- \sigma_1 \setminus \{\boldzero\}), \\
s^{d_1}t^{d_2}\cdot \dfrac{\< \df{h_1}{-d_2}\>}{\< h_2 \df{h_1}{-d_2} \>  }
&\text{if } {\bbd} \in (-\sigma_2 \setminus \{\boldzero\}).
\end{cases}
\]

In the remainder of the section, 
we compute the graded pieces of $JD(R_{A_3})$ to show that they are not equal to the graded pieces of $D(R_{A_3},J)$. 
This means that Proposition~\ref{prop:SmStDO-Prop1.6} does not hold for this ring. 

\begin{figure}[b]
\centering
\begin{subfigure}{1.55in}
  \centering
\begin{tikzpicture}[scale=0.45]
\draw[blue, thick] (-2.5,0) -- (3.5,0);
\draw[gray] (0,-6.5) -- (0,6.5);
\draw[blue, thick] (-2.33,-7) -- (2.33,7);
\draw[ultra thick,->,>=latex] (1,2) -- (0,0);
\draw[ultra thick,->,>=latex] (2,5) -- (1,3);
\draw[ultra thick,->,>=latex] (1,1) -- (0,0);
\draw[ultra thick,->,>=latex] (2,1) -- (1,0);
\draw[ultra thick,->,>=latex] (3,1) -- (2,0);
\foreach \y in {1,2}{ \node[draw,circle,inner sep=2pt,red,fill] at (1,\y) {};}
\foreach \x in {1,2}{ \node[draw,circle,inner sep=2pt,blue,fill] at (\x,3*\x) {};}
\foreach \x in {0,1,2,3}{ \node[draw,circle,inner sep=2pt,blue,fill] at (\x,0) {};}
\foreach \y in {1,2,...,7}{ \node[draw,circle,inner sep=2pt,red,fill] at (3,\y) {};}
\foreach \y in {1,2,...,5}{ \node[draw,circle,inner sep=2pt,red,fill] at (2,\y) {};}
\foreach \x in {0,1,2,3}{
      \foreach \y in {-7,-6,...,-1}{
        \node[draw,circle,inner sep=2pt,green,fill] at (\x,\y) {};}}
      \foreach \y in {-7,-6,...,-3}{
        \node[draw,circle,inner sep=2pt,green,fill] at (-1,\y) {};}
      \foreach \y in {-7,-6}{
        \node[draw,circle,inner sep=2pt,green,fill] at (-2,\y) {};}
\foreach \x in {-2,-1,0,1}{
      \foreach \y in {4,5,6}{
        \node[draw,circle,inner sep=2pt,yellow,fill] at (\x,\y) {};}}
\foreach \x in {-2,-1,0}{
      \foreach \y in {1,2,3}{
        \node[draw,circle,inner sep=2pt,yellow,fill] at (\x,\y) {};}}
     \foreach \x in {-2,-1,...,2}{ \node[draw,circle,inner sep=2pt,yellow,fill] at (\x,7) {};}
\foreach \x in {-1,-2}{ \node[draw,circle,inner sep=2pt,yellow,fill] at (\x,0) {};}
\foreach \y in {-5,-4,...,-1}{ \node[draw,circle,inner sep=2pt,violet,fill] at (-2,\y) {};}
\foreach \y in {-2,-1}{ \node[draw,circle,inner sep=2pt,violet,fill] at (-1,\y) {};}
\foreach \x in {-2,-1,...,3}{
      \foreach \y in {-7,-6,...,7}{
        \node[draw,circle,inner sep=2pt] at (\x,\y) {};}}
\end{tikzpicture}
\captionsetup{margin=0in,width=1.5in,font=small,justification=centering}
\caption{$(JD(R_{A_3}))_{\bbd}$ for $\bbt^{\bbd} \in J$. }
\label{fig:3D:diff-red}
\end{subfigure}%
\begin{subfigure}{1.55in}
\centering
\begin{tikzpicture}[scale=0.45]
\draw[blue, thick] (-2.5,0) -- (3.5,0);
\draw[gray] (0,-6.5) -- (0,6.5);
\draw[blue, thick] (-2.33,-7) -- (2.33,7);
\draw[ultra thick,->,>=latex] (0,2) -- (-1,0);
\draw[ultra thick,->,>=latex] (1,3) -- (0,1);
\draw[ultra thick,->,>=latex] (1,4) -- (0,2);
\draw[ultra thick,->,>=latex] (-1,2) -- (-2,0);
\draw[ultra thick,->,>=latex] (0,3) -- (-1,1);
\draw[ultra thick,->,>=latex] (-1,3) -- (-2,1);
\draw[ultra thick,->,>=latex] (0,4) -- (-1,2);
\draw[ultra thick,->,>=latex] (1,5) -- (0,3);
\draw[ultra thick,->,>=latex] (1,6) -- (0,4);
\draw[ultra thick,->,>=latex] (-1,4) -- (-2,2);
\draw[ultra thick,->,>=latex] (0,5) -- (-1,3);
\draw[ultra thick,->,>=latex] (-1,5) -- (-2,3);
\draw[ultra thick,->,>=latex] (-1,6) -- (-2,4);
\draw[ultra thick,->,>=latex] (-1,7) -- (-2,5);
\draw[ultra thick,->,>=latex] (0,6) -- (-1,4);
\draw[ultra thick,->,>=latex] (0,7) -- (-1,5);
\draw[ultra thick,->,>=latex] (2,6) -- (1,4);
\draw[ultra thick,->,>=latex] (1,7) -- (0,5);
\draw[ultra thick,->,>=latex] (2,7) -- (1,5);
\foreach \y in {1,2}{ \node[draw,circle,inner sep=2pt,red,fill] at (1,\y) {};}
\foreach \x in {1,2}{ \node[draw,circle,inner sep=2pt,blue,fill] at (\x,3*\x) {};}
\foreach \x in {0,1,2,3}{ \node[draw,circle,inner sep=2pt,blue,fill] at (\x,0) {};}
\foreach \y in {1,2,...,7}{ \node[draw,circle,inner sep=2pt,red,fill] at (3,\y) {};}
\foreach \y in {1,2,...,5}{ \node[draw,circle,inner sep=2pt,red,fill] at (2,\y) {};}
\foreach \x in {0,1,2,3}{
      \foreach \y in {-7,-6,...,-1}{
        \node[draw,circle,inner sep=2pt,green,fill] at (\x,\y) {};}}
      \foreach \y in {-7,-6,...,-3}{
        \node[draw,circle,inner sep=2pt,green,fill] at (-1,\y) {};}
      \foreach \y in {-7,-6}{
        \node[draw,circle,inner sep=2pt,green,fill] at (-2,\y) {};}
\foreach \x in {-2,-1,0,1}{
      \foreach \y in {4,5,6}{
        \node[draw,circle,inner sep=2pt,yellow,fill] at (\x,\y) {};}}
\foreach \x in {-2,-1,0}{
      \foreach \y in {1,2,3}{
        \node[draw,circle,inner sep=2pt,yellow,fill] at (\x,\y) {};}}
     \foreach \x in {-2,-1,...,2}{ \node[draw,circle,inner sep=2pt,yellow,fill] at (\x,7) {};}
\foreach \x in {-1,-2}{ \node[draw,circle,inner sep=2pt,yellow,fill] at (\x,0) {};}
\foreach \y in {-5,-4,...,-1}{ \node[draw,circle,inner sep=2pt,violet,fill] at (-2,\y) {};}
\foreach \y in {-2,-1}{ \node[draw,circle,inner sep=2pt,violet,fill] at (-1,\y) {};}
\foreach \x in {-2,-1,...,3}{
      \foreach \y in {-7,-6,...,7}{
        \node[draw,circle,inner sep=2pt] at (\x,\y) {};}}
\end{tikzpicture}
\captionsetup{margin=0in,width=1.5in,font=small,justification=centering}
\caption{$(JD(R_{A_3}))_{\bbd}$ for $\bbd \in \textbf{C2}$.}
\label{fig:3D:diff-yellow}
\end{subfigure}
\begin{subfigure}{1.55in}
\centering
\begin{tikzpicture}[scale=0.45]
\draw[blue, thick] (-2.5,0) -- (3.5,0);
\draw[gray] (0,-6.5) -- (0,6.5);
\draw[blue, thick] (-2.33,-7) -- (2.33,7);
\draw[ultra thick,->,>=latex] (1,0) -- (0,-1);
\draw[ultra thick,->,>=latex] (1,-1) -- (0,-2);
\draw[ultra thick,->,>=latex] (1,-2) -- (0,-3);
\draw[ultra thick,->,>=latex] (1,-3) -- (0,-4);
\draw[ultra thick,->,>=latex] (1,-4) -- (0,-5);
\draw[ultra thick,->,>=latex] (1,-5) -- (0,-6);
\draw[ultra thick,->,>=latex] (1,-6) -- (0,-7);
\draw[ultra thick,->,>=latex] (0,-2) -- (-1,-3);
\draw[ultra thick,->,>=latex] (0,-3) -- (-1,-4);
\draw[ultra thick,->,>=latex] (0,-4) -- (-1,-5);
\draw[ultra thick,->,>=latex] (0,-5) -- (-1,-6);
\draw[ultra thick,->,>=latex] (0,-6) -- (-1,-7);
\draw[ultra thick,->,>=latex] (-1,-5) -- (-2,-6);
\draw[ultra thick,->,>=latex] (-1,-6) -- (-2,-7);
\draw[ultra thick,->,>=latex] (2,0) -- (1,-1);
\draw[ultra thick,->,>=latex] (2,-1) -- (1,-2);
\draw[ultra thick,->,>=latex] (2,-2) -- (1,-3);
\draw[ultra thick,->,>=latex] (2,-3) -- (1,-4);
\draw[ultra thick,->,>=latex] (2,-4) -- (1,-5);
\draw[ultra thick,->,>=latex] (2,-5) -- (1,-6);
\draw[ultra thick,->,>=latex] (2,-6) -- (1,-7);
\draw[ultra thick,->,>=latex] (3,-2) -- (2,-3);
\draw[ultra thick,->,>=latex] (3,-3) -- (2,-4);
\draw[ultra thick,->,>=latex] (3,-4) -- (2,-5);
\draw[ultra thick,->,>=latex] (3,0) -- (2,-1);
\draw[ultra thick,->,>=latex] (3,-1) -- (2,-2);
\draw[ultra thick,->,>=latex] (3,-5) -- (2,-6);
\draw[ultra thick,->,>=latex] (3,-6) -- (2,-7);
\foreach \y in {1,2}{ \node[draw,circle,inner sep=2pt,red,fill] at (1,\y) {};}
\foreach \x in {1,2}{ \node[draw,circle,inner sep=2pt,blue,fill] at (\x,3*\x) {};}
\foreach \x in {0,1,2,3}{ \node[draw,circle,inner sep=2pt,blue,fill] at (\x,0) {};}
\foreach \y in {1,2,...,7}{ \node[draw,circle,inner sep=2pt,red,fill] at (3,\y) {};}
\foreach \y in {1,2,...,5}{ \node[draw,circle,inner sep=2pt,red,fill] at (2,\y) {};}
\foreach \x in {0,1,2,3}{
      \foreach \y in {-7,-6,...,-1}{
        \node[draw,circle,inner sep=2pt,green,fill] at (\x,\y) {};}}
      \foreach \y in {-7,-6,...,-3}{
        \node[draw,circle,inner sep=2pt,green,fill] at (-1,\y) {};}
      \foreach \y in {-7,-6}{
        \node[draw,circle,inner sep=2pt,green,fill] at (-2,\y) {};}
\foreach \x in {-2,-1,0,1}{
      \foreach \y in {4,5,6}{
        \node[draw,circle,inner sep=2pt,yellow,fill] at (\x,\y) {};}}
\foreach \x in {-2,-1,0}{
      \foreach \y in {1,2,3}{
        \node[draw,circle,inner sep=2pt,yellow,fill] at (\x,\y) {};}}
     \foreach \x in {-2,-1,...,2}{ \node[draw,circle,inner sep=2pt,yellow,fill] at (\x,7) {};}
\foreach \x in {-1,-2}{ \node[draw,circle,inner sep=2pt,yellow,fill] at (\x,0) {};}
\foreach \y in {-5,-4,...,-1}{ \node[draw,circle,inner sep=2pt,violet,fill] at (-2,\y) {};}
\foreach \y in {-2,-1}{ \node[draw,circle,inner sep=2pt,violet,fill] at (-1,\y) {};}
\foreach \x in {-2,-1,...,3}{
      \foreach \y in {-7,-6,...,7}{
        \node[draw,circle,inner sep=2pt] at (\x,\y) {};}}
\end{tikzpicture}
\captionsetup{margin=0in,width=1.5in,font=small,justification=centering}
\caption{$(JD(R_{A_3}))_{\bbd}$ for $\bbd \in \textbf{C4}$.}
\label{fig:3D:diff-green}
\end{subfigure}%
\begin{subfigure}{1.55in}
\centering
\begin{tikzpicture}[scale=0.45]
\draw[blue, thick] (-2.5,0) -- (3.5,0);
\draw[gray] (0,-6.5) -- (0,6.5);
\draw[blue, thick] (-2.33,-7) -- (2.33,7);
\draw[ultra thick,brown,->,>=latex] (0,0) -- (-1,-1);
\draw[ultra thick,brown,->,>=latex] (0,0) -- (-1,-2);
\draw[ultra thick,brown,->,>=latex] (-1,0) -- (-2,-1);
\draw[ultra thick,brown,->,>=latex] (-1,0) -- (-2,-3);
\draw[ultra thick,gray,->,>=latex] (-1,-1) -- (-2,-2);
\draw[ultra thick,gray,->,>=latex] (-1,-1) -- (-2,-3);
\draw[ultra thick,brown,->,>=latex] (-1,-2) -- (-2,-3);
\draw[ultra thick,brown,->,>=latex] (-1,-2) -- (-2,-4);
\draw[ultra thick,gray,->,>=latex] (-1,-3) -- (-2,-4);
\draw[ultra thick,gray,->,>=latex] (-1,-3) -- (-2,-5);
\foreach \y in {1,2}{ \node[draw,circle,inner sep=2pt,red,fill] at (1,\y) {};}
\foreach \x in {1,2}{ \node[draw,circle,inner sep=2pt,blue,fill] at (\x,3*\x) {};}
\foreach \x in {0,1,2,3}{ \node[draw,circle,inner sep=2pt,blue,fill] at (\x,0) {};}
\foreach \y in {1,2,...,7}{ \node[draw,circle,inner sep=2pt,red,fill] at (3,\y) {};}
\foreach \y in {1,2,...,5}{ \node[draw,circle,inner sep=2pt,red,fill] at (2,\y) {};}
\foreach \x in {0,1,2,3}{
      \foreach \y in {-7,-6,...,-1}{
        \node[draw,circle,inner sep=2pt,green,fill] at (\x,\y) {};}}
      \foreach \y in {-7,-6,...,-3}{
        \node[draw,circle,inner sep=2pt,green,fill] at (-1,\y) {};}
      \foreach \y in {-7,-6}{
        \node[draw,circle,inner sep=2pt,green,fill] at (-2,\y) {};}
\foreach \x in {-2,-1,0,1}{
      \foreach \y in {4,5,6}{
        \node[draw,circle,inner sep=2pt,yellow,fill] at (\x,\y) {};}}
\foreach \x in {-2,-1,0}{
      \foreach \y in {1,2,3}{
        \node[draw,circle,inner sep=2pt,yellow,fill] at (\x,\y) {};}}
     \foreach \x in {-2,-1,...,2}{ \node[draw,circle,inner sep=2pt,yellow,fill] at (\x,7) {};}
\foreach \x in {-1,-2}{ \node[draw,circle,inner sep=2pt,yellow,fill] at (\x,0) {};}
\foreach \y in {-5,-4,...,-1}{ \node[draw,circle,inner sep=2pt,violet,fill] at (-2,\y) {};}
\foreach \y in {-2,-1}{ \node[draw,circle,inner sep=2pt,violet,fill] at (-1,\y) {};}
\foreach \x in {-2,-1,...,3}{
      \foreach \y in {-7,-6,...,7}{
        \node[draw,circle,inner sep=2pt] at (\x,\y) {};}}
\end{tikzpicture}
\captionsetup{margin=0in,width=1.5in,font=small,justification=centering}
\caption{$(JD(R_{A_3}))_{\bbd}$ for $\bbd \in\textbf{C3}$.}
\label{fig:3D:diff-purple}
\end{subfigure}
\captionsetup{margin=0in,font=small,justification=centering}
\caption{Visualizing $(JD(R_{A_3}))_{\bbd}$.}
\end{figure}

To begin the computation, 
if $s^{m_1}t^{m_2} \in J$, then
\begin{align*}
D(R_{A_3})_{(d_1-m_1,d_2-m_2)}
&= s^{d_1-m_1}t^{d_2-m_2}\cdot 
\< \df{h_1}{-(d_2-m_2)-1}
\df{h_2}{-3(d_1-m_1)+d_2-m_2-1}\>.
\end{align*} 
Further, since $m_1,m_2 \geq 1$, $s^{m_1}t^{m_2}s^{d_1-m_1}t^{d_2-m_2}= s^{d_1}t^{d_2}$, and there is a containment 
\begin{align*}
&\< \df{h_1}{-(d_2-m_2)-1}\df{h_2}{-3(d_1-m_1)+d_2-m_2-1}\> \\
& \qquad\qquad \subseteq \< \df{h_1}{-d_2}\df{h_2}{-3d_1+d_2+1}, \df{h_1}{-d_2+1}\df{h_2}{-3d_1+d_2}\>. 
\end{align*}
Thus, to determine the graded piece of $(JD(R_{A_3}))_{\bbd}$ for any  $\bbd\in\ZZ^2$, 
it is enough to consider $D(R_{A_3})$ in multidegrees $(d_1-1,d_2-1)$ and $(d_1-1,d_2-2)$, where
\begin{align*}
D(R_{A_3})_{(d_1-1,d_2-1)} 
&= s^{d_1-1}t^{d_2-1} \cdot 
\< \df{h_1}{-d_2}\df{h_2}{-3d_1+d_2+1}\>
\quad\text{and}\\
D(R_{A_3})_{(d_1-1,d_2-2)}
&= s^{d_1-1}t^{d_2-2} \cdot \< \df{h_1}{-d_2+1}\df{h_2}{-3d_1+d_2}\>.
\end{align*} 

We will now break down this computation by chambers from Figure~\ref{fig:3D-lattice}, leaving off some half-lines along the way and addressing them as special cases later on. 
For \textbf{C1}, Figure~\ref{fig:3D:diff-red} helps to visualize that if $s^{d_1}t^{d_2} \in J$, then at least one of the two multidegrees  $(d_1-1,d_2-1)$ or $(d_1-1,d_2-2)$ lives in $\NN A$. 
Since $D(R_{A_3})_{\bbm}=s^{m_1}t^{m_2}\cdot \CC[\theta]$ for all $\bbm \in \NN A$, $(JD(R_{A_3}))_{\bbd}=s^{d_1}t^{d_2}\cdot \CC[\theta]$ for $s^{d_1}t^{d_2} \in J$. 
We will consider the blue multidegrees in $\textbf{C1}$ with their neighbors in \textbf{C2} and \textbf{C4}.

Now consider the graded pieces of  $(JD(R_{A_3}))_{\bbd}$ for $\bbd$ in $\sigma_2$ or \textbf{C2}, excluding the two half-lines in that chamber given by $-\sigma_1$ and $y=1$. 
Since we excluded the $x$-axis and $y=1$, both $(d_1-1,d_2-1)$ and $(d_1-1,d_2-2)$ lie in \textbf{C2} for all $\bbd$ under consideration, see Figure~\ref{fig:3D:diff-yellow}. 
Since 
\begin{align*}
D(R_{A_3})_{(d_1-1,d_2-1)} 
&= s^{d_1-1}t^{d_2-1} \cdot\< \df{h_2}{-3d_1+d_2+1}\> \quad\text{and}\\
D(R_{A_3})_{(d_1-1,d_2-2)} 
&= s^{d_1-1}t^{d_2-2} \cdot\< \df{h_2}{-3d_1+d_2} \>,
\end{align*}
and there is a containment 
$\left\< \df{h_2}{-3d_1+d_2+1} \right\>
\subseteq 
\left\< \df{h_2}{-3d_1+d_2}; \right\>$,   
for such $\bbd$, 
\[
(JD(R_{A_3}))_{\bbd} 
= s^{d_1}t^{d_2}\cdot \left\< \df{h_2}{-3d_1+d_2} \right\>.
\] 

Now consider $(JD(R_{A_3}))_{\bbd}$ for $\bbd$ in $\sigma_1$ or in \textbf{C4}, excluding those multidegrees that lie on $-\sigma_2$ or the half-line $y=3x-1$. 
Since we excluded the multidegrees on $-\sigma_2$ and $y=3x-1$, both $(d_1-1,d_2-1)$ and $(d_1-1,d_2-2)$ also lie in \textbf{C4}, as seen in  Figure~\ref{fig:3D:diff-green}. 
Since 
\begin{align*}
D(R_{A_3})_{(d_1-1,d_2-1)} 
&= s^{d_1-1}t^{d_2-1} \cdot\< \df{h_1}{-d_2}\>
\quad\text{and}\quad \\
D(R_{A_3})_{(d_1-1,d_2-2)} 
&= s^{d_1-1}t^{d_2-2} \cdot\< \df{h_1}{-d_2+1} \>,
\intertext{and there is a containment}
\< \df{h_1}{-d_2+1} \> &\subseteq \< \df{h_1}{-d_2}\>,
\intertext{it follows that for such $\bbd$,}
(JD(R_{A_3}))_{\bbd} 
&= s^{d_1}t^{d_2} \cdot \< \df{h_1}{-d_2} \>.
\end{align*}

For $(JD(R_{A_3}))_{\bbd}$ for $\bbd$ in \textbf{C3}, $-\sigma_1$, or $ -\sigma_2$, see Figure~\ref{fig:3D:diff-purple}, 
\begin{align*}
D(R_{A_3})_{(d_1-1,d_2-1)} 
&= s^{d_1-1}t^{d_2-1} \cdot 
\< \df{h_1}{-d_2} \df{h_2}{-3d_1+d_2+1} \>, \\
D(R_{A_3})_{(d_1-1,d_2-2)}
&= s^{d_1-1}t^{d_2-2} \cdot 
\< \df{h_1}{-d_2+1} \df{h_2}{-3d_1+d_2} \>.
\end{align*}
Since the ideals  
\[
\< \df{h_1}{-d_2} \df{h_2}{-3d_1+d_2+1} \> 
\quad\text{and}\quad
\< \df{h_1}{-d_2+1} \df{h_2}{-3d_1+d_2} \>
\]
are incomparable, it follows that for such $\bbd$, 
\begin{align*}
(JD&(R_{A_3}))_{\bbd} = s^{d_1}t^{d_2}\cdot 
\< \df{h_1}{-d_2} \df{h_2}{-3d_1+d_2+1}, \df{h_1}{-d_2+1} \df{h_2}{-3d_1+d_2}\>.
\end{align*}

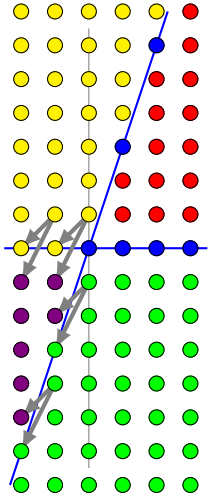
\begin{wrapfigure}[21]{l}[-5pt]{1.45in}
\vspace{-5pt}
\centering
\begin{tikzpicture}[scale=0.45]
\draw[blue, thick] (-2.5,0) -- (3.5,0);
\draw[gray] (0,-6.5) -- (0,6.5);
\draw[blue, thick] (-2.33,-7) -- (2.33,7);
\draw[ultra thick,gray,->,>=latex] (0,1) -- (-1,0);
\draw[ultra thick,gray,->,>=latex] (0,1) -- (-1,-1);
\draw[ultra thick,gray,->,>=latex] (-1,1) -- (-2,0);
\draw[ultra thick,gray,->,>=latex] (-1,1) -- (-2,-1);
\draw[ultra thick,gray,->,>=latex] (0,-1) -- (-1,-2);
\draw[ultra thick,gray,->,>=latex] (0,-1) -- (-1,-3);
\draw[ultra thick,gray,->,>=latex] (-1,-4) -- (-2,-5);
\draw[ultra thick,gray,->,>=latex] (-1,-4) -- (-2,-6);

\foreach \y in {1,2}{ \node[draw,circle,inner sep=2pt,red,fill] at (1,\y) {};}
\foreach \x in {1,2}{ \node[draw,circle,inner sep=2pt,blue,fill] at (\x,3*\x) {};}
\foreach \x in {0,1,2,3}{ \node[draw,circle,inner sep=2pt,blue,fill] at (\x,0) {};}
\foreach \y in {1,2,...,7}{ \node[draw,circle,inner sep=2pt,red,fill] at (3,\y) {};}
\foreach \y in {1,2,...,5}{ \node[draw,circle,inner sep=2pt,red,fill] at (2,\y) {};}
\foreach \x in {0,1,2,3}{
      \foreach \y in {-7,-6,...,-1}{
        \node[draw,circle,inner sep=2pt,green,fill] at (\x,\y) {};}}
      \foreach \y in {-7,-6,...,-3}{
        \node[draw,circle,inner sep=2pt,green,fill] at (-1,\y) {};}
      \foreach \y in {-7,-6}{
        \node[draw,circle,inner sep=2pt,green,fill] at (-2,\y) {};}
\foreach \x in {-2,-1,0,1}{
      \foreach \y in {4,5,6}{
        \node[draw,circle,inner sep=2pt,yellow,fill] at (\x,\y) {};}}
\foreach \x in {-2,-1,0}{
      \foreach \y in {1,2,3}{
        \node[draw,circle,inner sep=2pt,yellow,fill] at (\x,\y) {};}}
     \foreach \x in {-2,-1,...,2}{ \node[draw,circle,inner sep=2pt,yellow,fill] at (\x,7) {};}
\foreach \x in {-1,-2}{ \node[draw,circle,inner sep=2pt,yellow,fill] at (\x,0) {};}
\foreach \y in {-5,-4,...,-1}{ \node[draw,circle,inner sep=2pt,violet,fill] at (-2,\y) {};}
\foreach \y in {-2,-1}{ \node[draw,circle,inner sep=2pt,violet,fill] at (-1,\y) {};}
\foreach \x in {-2,-1,...,3}{
      \foreach \y in {-7,-6,...,7}{
        \node[draw,circle,inner sep=2pt] at (\x,\y) {};}}
\end{tikzpicture}
\captionsetup{margin=0in,width=1.25in,font=small,labelsep=newline,justification=centering}
\caption{($JD(R_{A_3}))_{\bbd}$ for $\bbd$ on exceptional lines in \textbf{C2} and \textbf{C4}}
\label{fig:3D:diff-cross-lines}
\end{wrapfigure}

The multidegrees $\bbd$ that we have not yet considered for $(JD(R_{A_3}))_{\bbd}$ are those in \textbf{C2} along $y=1$ and in \textbf{C4} along $y=3x-1$. 
For such ${\bbd}$, one of $(d_1-1,d_2-1)$ or $(d_1-1,d_2-2)$ belongs to \textbf{C3}, 
see Figure~\ref{fig:3D:diff-cross-lines}. 
First, for $\bbd$ in \textbf{C2} along the line $y=1$, $(d_1-1,-1)\in\textbf{C3}$, and
\[
\hspace*{2in}
D(R_{A_3})_{(d_1-1,-1)} = 
s^{d_1-1}t^{-1}\cdot 
\< h_1 \df{h_2}{-3d_1+1} \>.
\] 
Combining the fact that 
$D(R_{A_3})_{(d_1-1,0)} = 
s^{d_1-1} \cdot \< \df{h_2}{-3d_1+2} \>$,  
we compute that 
\[
\hspace*{2in}
(JD(R_{A_3}))_{\bbd} 
= s^{d_1}t \cdot 
\< \df{h_2}{-3d_1+2}, h_1\df{h_2}{-3d_1+1} \>.
\] 
Second, for $\bbd\in\textbf{C4}$ along the line $y=3x-1$, 
\[
\hspace*{2in}
D(R_{A_3})_{(d_1-1,3d_1-2)} = 
s^{d_1-1}t^{3d_1-2}\cdot\< \df{h_1}{-3d_1+1} h_2\>. 
\]
Connecting with the fact that 
\[
\hspace*{2in}
D(R_{A_3})_{(d_1-1,3d_1-3)} 
= s^{d_1-1}t^{3d_1-3} \cdot \< \df{h_1}{-3d_1+2} \>,
\]
we determine that 
\[
\hspace*{2in}
(JD(R_{A_3}))_{\bbd} 
= s^{d_1}t^{3d_1-1}\cdot 
\< \df{h_1}{-3d_1+1} h_2, \df {h_1}{-3d_1+2} \>.
\]
Having now computed $JD(R_{A_3})$ in all multidegrees, we have that 
\[
(JD(R_{A_3}))_{\bbd} = 
\begin{cases}
s^{d_1}t^{d_2}\cdot \CC[\theta] 
&\text{if } s^{d_1}t^{d_2} \in J,\\
s^{d_1}t^{d_2}\cdot\< \df{h_2}{-3d_1+d_2} \> 
&\text{if }{\bbd} \in \left(\textbf{C2} \setminus (\sigma_1 \cup \{y =1\})\right) \cup ( \sigma_2 \setminus \{\boldzero\}),\\
s^{d_1}t^{d_2}\cdot \\
\ \ \  
\< \df{h_1}{-d_2} \df{h_2}{-3d_1 
+d_2+1}, &\\
\ \ \quad \df{h_1}{-d_2+1}\df{h_2}{-3d_1+d_2}\> 
&\text{if }{\bbd} \in \textbf{C3} \cup (-\sigma_1 \cup -\sigma_2)  \cup ( \NN A_3), 
\\
s^{d_1}t^{d_2}\cdot \< \df{h_1}{-d_2}\> 
&\text{if }{\bbd} \in \left(\textbf{C4} \setminus (\sigma_2 \cup \{ y=3x-1 \})\right) \cup ( \sigma_1 \setminus \{\boldzero\}).
\end{cases}
\]
Now compare this with $D(R_{A_3},J)$ from \eqref{eq:3D:DRJ}; it becomes clear that 
$(JD(R_{A_3}))_{\bbd} \neq D(R_{A_3},J)_{\bbd}$ whenever $\bbd$ belongs to any of the following: 
$-\sigma_1$, $-\sigma_2$, \textbf{C3}, \textbf{C2} along the half-line $\{y=1\}$, or \textbf{C4} along the half-line $\{y=3x-1\}$.
In particular, $\II(J)/D(R_{A_3},J) \neq \II(J)/JD(R_{A_3})$, in contrast to Proposition~\ref{prop:SmStDO-Prop1.6}.

\section{Differential operators on a rational normal curve}
\label{sec:higherD}

The techniques used in Sections~\ref{sec:2D} and~\ref{sec:3D} can also be applied to compute $\II(J)_{\bbd}$, $D(R_{A_n},J)_{\bbd}$, $(JD(R_{A_n}))_{\bbd}$ and $\left(\II(J)/D(R_{A_n},J) \right)_{\bbd}$, 
where the radical ideal $J=\<st,st^2,\ldots,st^{n-1}\>$ is again the intersection of the primes defined by the facets of $A_n$.  
We denote the facets of $A_n$ by 
\begin{align*}
\sigma_1 = 
\{(x,y) \in \NN^2 \mid x \geq 0, y=0\}
\quad\text{and}\quad
\sigma_2 = 
\{(x,y) \in \NN^2 \mid x, y \geq 0, y=3x\},
\end{align*} 
which have primitive integral support functions 
\[
h_1=\theta_2 \quad \text{and}\quad h_2=3\theta_1-\theta_2.
\] 
Again, we divide $\ZZ^2$ into chambers, analogous to those used in Sections~\ref{sec:2D} and~\ref{sec:3D}, so  
\begin{itemize}
\item[\textbf{C1}] $= \{\bbd\in\ZZ^2\mid \bbd\in\NN A\}$, 
\item[\textbf{C2}] $= \{\bbd\in\ZZ^2\mid h_1({\bbd})\geq 0, h_2 ({\bbd})< 0\}$, 
\item[\textbf{C3}] $= \{\bbd\in\ZZ^2\mid h_1({\bbd})<0, h_2({\bbd})<0\}$, and 
\item[\textbf{C4}] $= \{\bbd\in\ZZ^2\mid  h_1({\bbd})<0, h_2 ({\bbd})\geq 0\}$. 
\end{itemize}

These computations yield the following formulas:
\[
\II(J)_{\bbd} =
\begin{cases}
s^{d_1}t^{d_2}\cdot\CC[\theta] 
&\text{if } {\bbd} \in \NN A_n,\\
s^{d_1}t^{d_2}\cdot \< \df{h_1}{-d_2} \> 
&\text{if } {\bbd} \in \textbf{C4},\\
s^{d_1}t^{d_2} \cdot\< \df{h_2}{-n d_1+d_2} \> 
&\text{if } {\bbd} \in \textbf{C2},\\
s^{d_1}t^{d_2} \cdot\< \df{h_1}{-d_2} \df{h_2}{-n d_1+d_2}\> 
&\text{if } {\bbd} \in \textbf{C3},
\end{cases}
\]  
\[
D(R_{A_n},J)_{\bbd} = 
\begin{cases}
s^{d_1}t^{d_2}\cdot\CC[\theta] 
&\text{if } s^{d_1}t^{d_2} \in J,\\
s^{d_1}t^{d_2}\cdot\< \df{h_1}{-d_2} \> 
&\text{if } {\bbd} \in (\textbf{C4}\setminus (-\sigma_2)) \cup (\sigma_1 \setminus \{\boldzero\}), \\
s^{d_1}t^{d_2}\cdot\< \df{h_2}{-n d_1+d_2} \>  
&\text{if } {\bbd} \in (\textbf{C2}\setminus (-\sigma_1)) \cup (\sigma_2 \setminus \{\boldzero\}), \\
s^{d_1}t^{d_2}\cdot\< \df{h_1}{-d_2}\df{h_2}{-n d_1+d_2} \> 
&\text{if } {\bbd} \in \textbf{C3} \cup(-\sigma_1 \cup -\sigma_2), \\
\end{cases}
\]
\[
(JD(R_{A_n}))_{\bbd} = 
\begin{cases}
s^{d_1}t^{d_2}\cdot\CC[\theta] 
&\text{if } s^{d_1}t^{d_2} \in J,\\
s^{d_1}t^{d_2}\cdot\left\< \df{h_1}{-d_2}\right\> 
&\text{if }{\bbd} \in \textbf{C4} \setminus
\left(\displaystyle\bigcup_{j=0}^{n-2}\{ y=n x-j \}\right),\\
s^{d_1}t^{d_2}\cdot\left\< \df{h_2}{-n x+y} \right\> 
&\text{if }{\bbd} \in \textbf{C2} \setminus 
\left(
    \displaystyle\bigcup_{j=0}^{n-2}\{y=j\}
\right),\\
s^{d_1}t^{d_2}\cdot
\left\< 
    \df{h_1}{-d_2+j}\, \cdot \right.&\\
\ \ 
\left.
    \df{h_2}{-n d_1+d_2+n-2-j} 
\right\>_{j=0}^{n-2} 
&\text{if }{\bbd} \text{ in }  
\textbf{C3}
\cup 
\left[
    \textbf{C2} \cap 
    \left(\displaystyle\bigcup_{j=0}^{n-2}\{y=j\}\right)
\right]
\cup \\
& \qquad\qquad\qquad\ \ \ 
\left[
    \textbf{C4} \cap 
    \left(\displaystyle\bigcup_{j=0}^{n-2}\{ y=n x-j \}\right)
\right]. 
\end{cases}
\]

Recall from Section \ref{sec:3D} that $JD(R_{A_3})\neq D(R_{A_3},J)$; we see that this is true for all rings of rational normal curves  $R_{A_n}=\CC[s,st,\dots,st^n]$ with $n \geq 3$.  
Specifically, comparing $D(R_{A_n},J)_{\bbd}$ and $(JD(R_{A_n}))_{\bbd}$ for various $\bbd\in\ZZ^2$, 
the graded pieces differ whenever $\bbd$ belongs to 
$\sigma_1$, 
$\sigma_2$,  
\textbf{C3}, 
the half-lines in \textbf{C2} inside 
\[
\{(x,y)\in\RR^2\mid x<0, y=i, 0\leq i\leq n-2\},
\]
or the half-lines in \textbf{C4} inside 
\[
\{(x,y)\in\RR^2\mid x<0, y=nx-i, 0\leq i\leq n-2\}.
\]
Figure~\ref{fig:higherD:7-drj-jdr} has these multidegrees $\bbd$ in blue for the case $n=7$.

\begin{wrapfigure}[16]{r}[5pt]{1.4in}
\centering
\begin{tikzpicture}[scale=0.65]
\draw[blue, thick] (-2.5,0) -- (2.5,0);
\draw[gray] (0,-6.5) -- (0,6.5);
\draw[blue, thick] (-1,-7) -- (1,7);
\draw[thick, blue,->,>=latex] (0,1) -- (-2.5,1);
\draw[thick, blue,->,>=latex] (0,2) -- (-2.5,2);
\draw[thick, blue,->,>=latex] (0,3) -- (-2.5,3);
\draw[thick, blue,->,>=latex] (0,4) -- (-2.5,4);
\draw[thick, blue,->,>=latex] (0,5) -- (-2.5,5);
\draw[thick, blue,->,>=latex] (0,-1) -- (-.86,-7);
\draw[thick, blue,->,>=latex] (0,-2) -- (-.71,-7);
\draw[thick, blue,->,>=latex] (0,-3) -- (-.57,-7);
\draw[thick, blue,->,>=latex] (0,-4) -- (-.43,-7);
\draw[thick, blue,->,>=latex] (0,-5) -- (-.29,-7);

\foreach \x in {-2,-1,...,2}{
      \foreach \y in {-7,-6,...,7}{
        \node[draw,circle,inner sep=2pt] at (\x,\y) {};}}
\foreach \x in{-2,-1,...,2}{ \node[draw,circle,blue,inner sep=2pt,fill] at (\x,0) {};}
\foreach \x in{-1,0,1}{ \node[draw,circle,blue,inner sep=2pt,fill] at (\x,7*\x) {};}
\foreach \y in{-1,-2,...,-5}{ \node[draw,circle,blue,inner sep=2pt,fill] at (0,\y) {};}
\foreach \x in {-2,-1,0}{
        \foreach \y in {1,2,...,5}{ \node[draw,circle,blue,inner sep=2pt,fill] at (\x,\y) {};}}
        
\foreach \x in {-2,-1}{
      \foreach \y in {-7,-6,...,-1}{
        \node[draw,circle,inner sep=2pt,blue,fill] at (\x,\y) {};}}
\end{tikzpicture}
\captionsetup{margin=0in,width=1in,font=small,labelsep=newline,justification=centering}
\caption{
$(JD(R_{A_7}))_{\bbd}$ and  $D(R_{A_7},J)_{\bbd}$ differ at blue $\bbd$
}
\label{fig:higherD:7-drj-jdr}
\end{wrapfigure}

As with the rational normal cones in degrees $2$ and $3$ from Sections~\ref{sec:2D} and~\ref{sec:3D}, 
for $\bbd\in\ZZ^2$, 
\begin{equation}
\left(\dfrac{\II(J)}{D(R_{A_n},J)}\right)_{\bbd} = 
\begin{cases}
0 
&\text{if } {\bbd} \in \ZZ^2\setminus (\ZZ \sigma_1 \cup \ZZ \sigma_2), \\
\dfrac{\CC[\theta]}{\< h_1h_2\>} 
&\text{if } {\bbd=\boldzero}, \\
s^{d_1}\cdot\dfrac{\< \df{h_2}{-n d_1} \>}{\< h_1\df{h_2}{-n d_1}\>  }
&\text{if } {\bbd} \in \ZZ \sigma_1 \setminus \{\boldzero\},\\
s^{d_1}t^{d_2}\cdot\dfrac{\< \df{h_1}{-d_2}\>}{\< h_2 \df{h_1}{-d_2} \>  }
&\text{if } {\bbd} \in \ZZ \sigma_2 \setminus \{\boldzero\}.
\end{cases}
\hspace*{1.5in}
\label{eq:diamond}
\end{equation}

Example~\ref{ex:Traves} considered the ordinary double point $\CC[x,y]/\<x y\>$, which is isomorphic to $R_{A_n}/J$ for all $n \geq 1$. 
Comparing 
$D(\CC[x,y]/\<xy\>))_{\bbd}$ from  \eqref{eq:bigtriangleup} and $\left(\II(J)/D(R_{A_n},J)\right)_{\bbd}$ from
\eqref{eq:diamond}, 
we see that there is an isomorphism between the graded components, viewed as $\CC$-vector spaces, given by $\varphi\colon D\left({\CC[x,y]}/{\<x y\>}\right) \rightarrow D\left({R_{A_n}}/{J}\right)$ with, for $m\in\ZZ$, 
\begin{align*}
\varphi\left(\frac{\CC[\theta_x,\theta_y]}{\< \theta_x \theta_y\> }\right ) 
&=\displaystyle\frac{\CC[\theta]}{\< h_1 h_2\> },
\hspace*{1.5in}
\\
\varphi\left (x^m \displaystyle\frac{\< \df{\theta_x}{-m} \>}{\< \df{\theta_x}{-m}\theta_y\> }\right ) 
&= s^m \cdot \displaystyle\frac{\< \df{h_2}{-nm} \>}{\< \df{h_2}{-nm}h_1\> }, 
\quad \text{and}
\hspace*{1.5in}
\\
\varphi\left (y^m \cdot \displaystyle\frac{\< \df{\theta_y}{-m} \>}{\< \df{\theta_y}{-m}\theta_x\> }\right )
&= s^mt^{n m} \cdot\displaystyle\frac{\< \df{h_1}{-n m} \>}{\< \df{h_1}{-n m}h_2\> }.
\hspace*{1.5in}
\end{align*}

However, this isomorphism is not an isomorphism of rings. Again as in Example~\ref{ex:A_2:P1}, the degree of the polynomial generator in $\theta_1$ and $\theta_2$ in multidegree $\bbd$ is $n$ times the degree of the polynomial generator of in $\theta_x$ and $\theta_y$. In fact, noting that along $\ZZ\sigma_2$, $d_2=nd_1$, we have
\begin{equation}
D\left(R_{A_n}/J\right)_{\bbd} = 
\begin{cases}
0 
&\text{if } {\bbd} \in \ZZ^2\setminus (\ZZ \sigma_1 \cup \ZZ \sigma_2), \\
\dfrac{\CC[\theta]}{\left\< \dfrac{h_1}{n}\dfrac{h_2}{n}\right\>} 
&\text{if } {\bbd=\boldzero}, \\
s^{d_1}\cdot\dfrac{\left\< \df{\dfrac{h_2}{n}}{- d_1} \right\>}{\left\< \dfrac{h_1}{n}\df{\dfrac{h_2}{n}}{- d_1}\right\>  }
&\text{if } {\bbd} \in \ZZ \sigma_1 \setminus \{\boldzero\},\\
s^{d_1}t^{d_2}\cdot\dfrac{\left\< \df{\dfrac{h_1}{n}}{-d_1}\right\>}{\left\< \dfrac{h_2}{n} \df{\dfrac{h_1}{n}}{-d_1} \right\>  }
&\text{if } {\bbd} \in \ZZ \sigma_2 \setminus \{\boldzero\}.
\end{cases}
\label{eq:dopsratnormn}
\end{equation}

Now comparing \eqref{eq:bigtriangleup} and \eqref{eq:dopsratnormn}, we see there is an isomorphism of the ring of differential operators between the two rings given by, for each $m\in\ZZ$, 
\[
\psi(\df{\theta_x}{-m}+\<\df{\theta_x}{-m}\theta_y \>) = s^{-m}\cdot \df{\frac{h_2}{n}}{-m}+\left\<\df{\frac{h_2}{n}}{-m}h_1 \right\>, 
\]
\[
\psi(\df{\theta_y}{-m}+\<\df{\theta_y}{-m}\theta_x \>) = s^{-m}t^{-nm} \cdot \df{\frac{h_1}{n}}{-m}+\left\<\df{\frac{h_1}{n}}{-m}h_2 \right\>,
\]
\[
\psi(x^{m}) = s^m, 
\quad \text{and} \quad 
\psi(y^{m}) = s^mt^{nm}.
\]

\section{Higher dimensional examples}
\label{sec:hdim}

Thus far, we have considered 
differential operators determined by radical ideals in two-dimensional semigroup rings. 
In this section, we turn to looking at some of the differential operators in the three-dimensional semigroup ring given by 
\[
A=\begin{bmatrix} 1 & 1 & 1 & 1\\
0 & 1 & 0 & 1\\
0 & 0 & 1 & 1\\ 
\end{bmatrix}
\]
and describe some subsets of differential operators of $R_A$ given by two different radical ideals $J$. 
For both choices of $J$, we will compute the graded pieces $\II(J)_{\bbd}$, $D(R_A,J)_{\bbd}$ and $\left(\II(J)/D(R_A,J)\right)_{\bbd}$. 
For the given matrix $A$, 
\[
R_A 
=\CC[\NN A]
=\CC[t_1,t_1t_2,t_1t_3,t_1t_2t_3].
\]  
The facets of the cone $\RR_{\geq 0} A$ are 
\[
\sigma_1=\NN\{ e_1, e_1+e_2\}, \ \ \ 
\sigma_2=\NN\{ e_1, e_1+e_3\},\ \ \  
\sigma_3=\NN\{ e_1+e_3, e_1+e_2+e_3\}, \ \ \   
\sigma_4=\NN\{ e_1+e_2, e_1+e_2+e_3\}, 
\]
and the corresponding primitive integral support functions are 
\[
h_1 = h_1(\theta) =  \theta_3,
\quad h_2 = h_2(\theta) = \theta_2,
\quad 
h_3 =h_3(\theta) =  \theta_1-\theta_3,   \quad
h_4 = h_4(\theta) = \theta_1-\theta_2.
\]
The prime ideals associated to the facets are 
\[
P_{\sigma_1}=\<t_1t_3,t_1t_2t_3\>, \quad
P_{\sigma_2}=\<t_1t_2,t_1t_2t_3\>,  \quad 
P_{\sigma_3}=\<t_1, t_1t_2\>, \quad 
P_{\sigma_4}=\<t_1,t_1t_3\>, 
\]
and the prime ideals associated to the rays 
(or $1$-dimensional faces) of the cone are  
\begin{alignat*}{2}
&P_{\sigma_1 \cap \sigma_2} = \<t_1t_2,t_1t_3,t_1t_2t_3\>, \qquad
&&P_{\sigma_2 \cap \sigma_3} = \<t_1,t_1t_2,t_1t_2t_3\>,
\\
&P_{\sigma_3 \cap \sigma_4} = \<t_1,t_1t_2,t_1t_3\>, 
&&P_{\sigma_1 \cap \sigma_4} = \<t_1,t_1t_3,t_1t_2t_3\>.
\end{alignat*}

\begin{wrapfigure}[13]{l}[5pt]{2.7in}
    \centering
    \includegraphics[scale=.12]{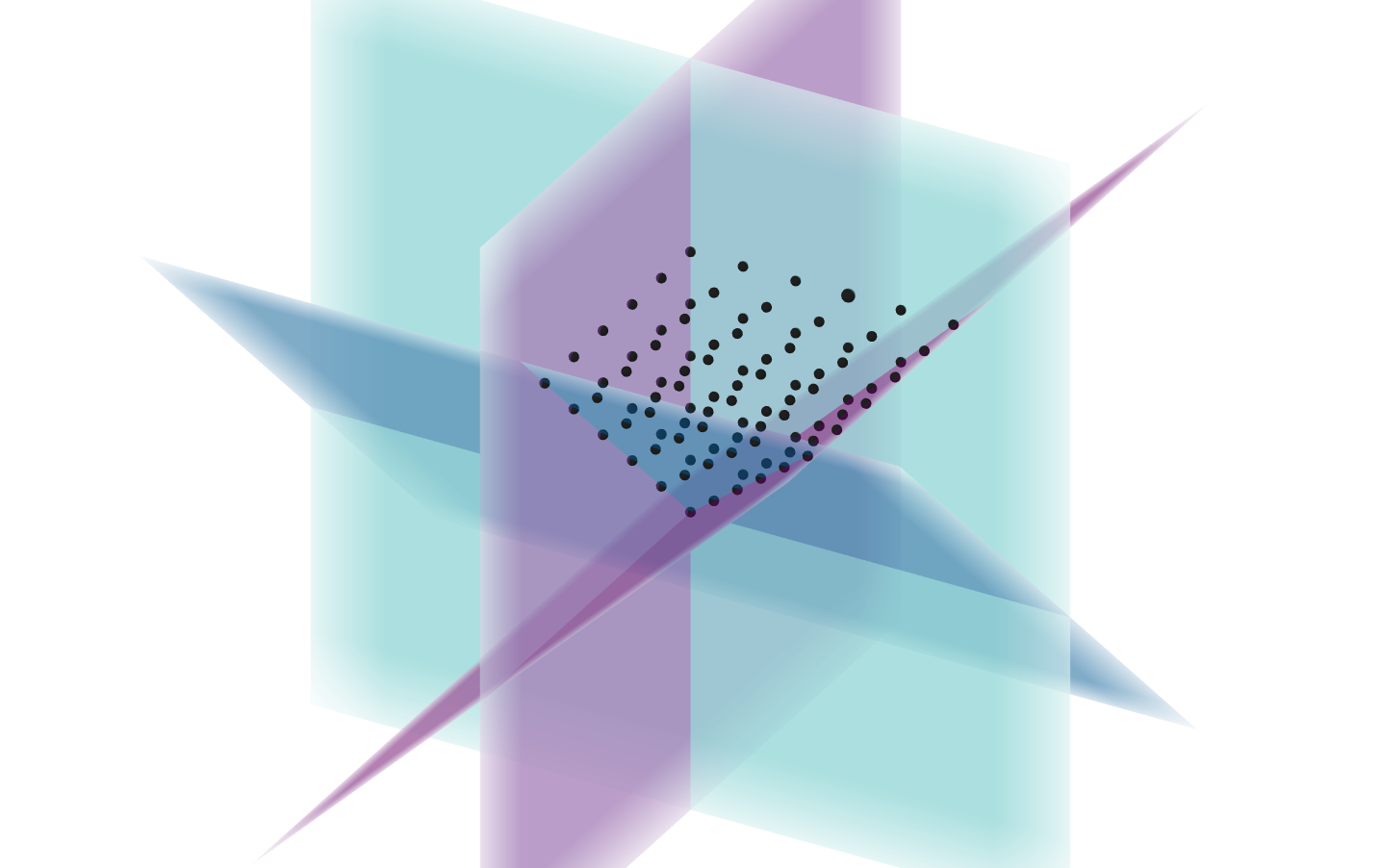}
\captionsetup{margin=0in,width=2.5in,font=small,labelsep=newline,justification=centering}
\caption{Chambers for $\CC[t_1,t_1t_2,t_1t_3,t_1t_2t_3]$}
    \label{fig:14Chambers}
\end{wrapfigure}
In contrast with the two dimensional cases considered in 
Sections~\ref{sec:2D},~\ref{sec:3D}, and~\ref{sec:higherD}, 
the increased dimension of the semigroup $\NN A$ allows for many more choices for radical monomial ideal $J$ in $R_A$. 
We will compute $\II(J)/D(R_A,J)$ for two such choices, 
\begin{align*}
\hspace*{2.4in}
J 
= P_{\sigma_1} \cap P_{\sigma_2} \cap P_{\sigma_3} \cap P_{\sigma_4}
= \<t_1^2 t_2 t_3\> 
\quad\text{and}\\
\hspace*{2.4in}
J 
= P_{\sigma_1} \cap P_{\sigma_2} \cap P_{\sigma_3 \cap \sigma_4} 
= \< t_1^2t_2t_3, t_1^2t_2^2t_3, t_1^2t_2t_3^2\>.
\end{align*}
For both choices of $J$, to compute the graded pieces of $\II(J)$, we will divide up $\ZZ^3$ into 14 chambers, depending on various combinations of signs of $h_i({\bbd})$. 
Table~\ref{tab:chambers} describes the 14 chambers in a list, while 
Figure~\ref{fig:14Chambers} illustrates them, with chamber \textbf{C1} (given by $\NN A$) in the top right with some of the lattice points of $\NN A$ shown. Note the $x$-axis is the vertical axis in this picture.

We must modify the chambers slightly to compute the graded pieces of $D(R_A,J)$, as we now consider which differential operators, when applied to an element in $R_A$, yield an element in $J$.  
Although all of the operators in the graded pieces of $\II(J)$ lying on the facets send $J$ into $J$, it is not necessarily the case that these operators applied to an element of $R_A$ will output an element in $J$.  
So we will switch the lattice points on the facets that correspond to monomials that do not lie in $J$ to lie in the adjacent chambers by interchanging $\geq$ with $>$.  
We will describe these new \emph{regions} within the examples.

\begin{ex}
\label{ex:3dinterior}
First consider the ideal 
\[
J= P_{\sigma_1} \cap P_{\sigma_2} \cap P_{\sigma_3} \cap P_{\sigma_4}=\<t_1^2 t_2 t_3\>
\quad 
\text{in}\quad 
R_A = \CC[t_1,t_1t_2,t_1t_3,t_1t_2t_3].
\]
\end{ex}
%
\begin{wrapfigure}[15]{l}[-10pt]{2.4in}
\centering
    \includegraphics[scale=.2]{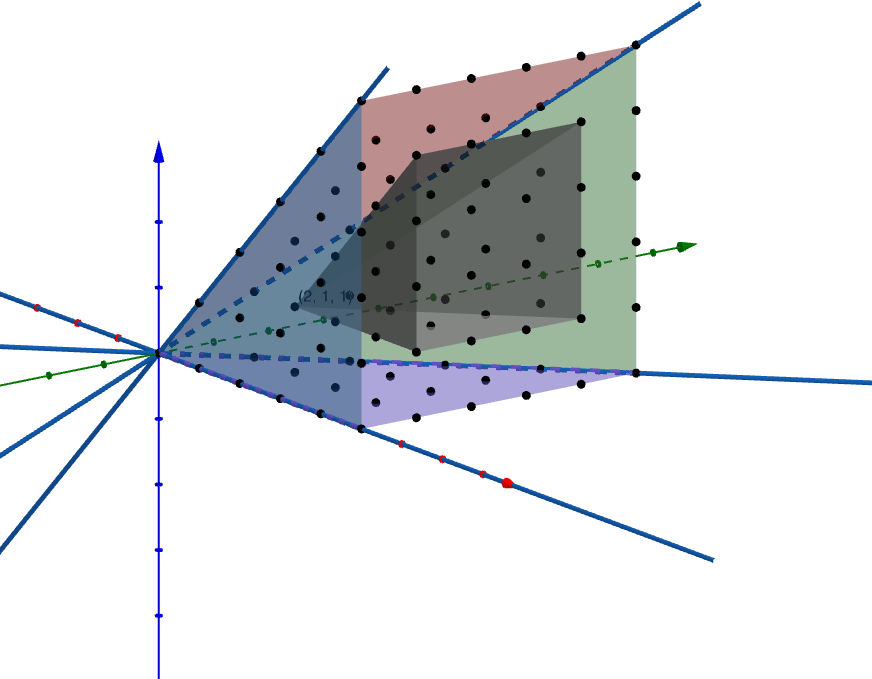}
\captionsetup{margin=0in,width=2.4in,font=small,labelsep=newline,justification=centering}    
    \caption{Cone  of $R_A$ 
    with  ideal $J = \<t_1^2 t_2 t_3\>$}
    \label{fig:x2yz}
\end{wrapfigure}

\renewcommand{\arraystretch}{1.5}
\begin{table}[b]
$\begin{array}{@{\extracolsep{\fill}} |c|l|l|}
\hline 
\text{Chamber} & \text{Halfspace inequalities} &  \text{Lattice point inequalities}
\\ \hline 
\textbf{C1} 
&\{{\bbd} \in \ZZ^3 \mid h_1({\bbd}), h_2({\bbd}), h_3({\bbd}), h_4({\bbd}) \geq 0\} 
& \NN A 
\\ \hline 
\textbf{C2} 
& \{{\bbd} \in \ZZ^3 \mid h_1({\bbd}), h_3({\bbd}) \geq 0, h_2({\bbd})>0, h_4({\bbd})< 0\} 
& \{\bbd\in\ZZ^3\mid d_2> d_1 \geq d_3 \geq 0\}
\\ \hline 
\textbf{C3} 
&\{{\bbd} \in \ZZ^3 \mid h_2({\bbd}), h_4({\bbd}) \geq 0, h_1({\bbd}) >0, h_3({\bbd}) < 0\} 
& \{\bbd\in\ZZ^3\mid d_3> d_1 \geq d_2 \geq 0\}
\\ \hline 
\textbf{C4} 
&\{{\bbd} \in \ZZ^3 \mid h_1({\bbd}), h_3({\bbd}) \geq 0, h_4({\bbd}) > 0, h_2({\bbd})<0\} 
& \{\bbd\in\ZZ^3\mid d_1 \geq d_3\geq 0 >d_2\}
\\ \hline 
\textbf{C5} 
&\{{\bbd} \in \ZZ^3 \mid h_2({\bbd}), h_4({\bbd}) \geq 0 , h_3({\bbd})>0, h_1({\bbd})<0\} 
& \{\bbd\in\ZZ^3\mid d_1 \geq d_2\geq 0 >d_3\}
\\ \hline 
\textbf{C6} 
&\{{\bbd} \in \ZZ^3 \mid h_1({\bbd}), h_2({\bbd}) \geq 0, h_3({\bbd}), h_4({\bbd})<0\} 
& \{\bbd\in\ZZ^3\mid d_3\geq 0, d_2\geq 0, d_2>d_1, d_3 > d_1\}
\\ \hline 
\textbf{C7} 
&\{{\bbd} \in \ZZ^3 \mid h_2({\bbd}), h_3({\bbd}) \geq 0, h_1({\bbd}), h_4({\bbd})<0\} 
& \{\bbd\in\ZZ^3\mid d_2\geq 0> d_3, d_2> d_1\geq d_3\}
\\ \hline 
\textbf{C8} 
&\{{\bbd} \in \ZZ^3 \mid h_1({\bbd}), h_4({\bbd}) \geq 0, h_2({\bbd}), h_3({\bbd})<0\} 
& \{\bbd\in\ZZ^3\mid d_3\geq 0> d_2, d_3> d_1\geq d_2\}
\\ \hline 
\textbf{C9} 
&\{{\bbd} \in \ZZ^3 \mid h_3({\bbd}), h_4({\bbd}) \geq 0, h_1({\bbd}), h_2({\bbd})<0\} 
& \{\bbd\in\ZZ^3\mid 0 > d_2, 0>d_3, d_1\geq d_2, d_1\geq d_3\}
\\ \hline 
\textbf{C10} 
&\{{\bbd} \in \ZZ^3 \mid h_1({\bbd}), h_3({\bbd}), h_4({\bbd}) < 0, h_2({\bbd})\geq 0\} 
& \{\bbd\in\ZZ^3\mid d_2\geq 0  > d_3>d_1\}
\\ \hline 
\textbf{C11} 
&\{{\bbd} \in \ZZ^3 \mid h_2({\bbd}), h_3({\bbd}), h_4({\bbd}) < 0, h_1({\bbd})\geq 0\} 
& \{\bbd\in\ZZ^3\mid d_3\geq 0  > d_2>d_1\}
\\ \hline 
\textbf{C12} 
&\{{\bbd} \in \ZZ^3 \mid h_1({\bbd}), h_2({\bbd}), h_3({\bbd}) < 0, h_4({\bbd})\geq 0\} 
& \{\bbd\in\ZZ^3\mid 0 >d_3 > d_1\geq d_2\}
\\ \hline 
\textbf{C13} 
&\{{\bbd} \in \ZZ^3 \mid h_1({\bbd}), h_2({\bbd}), h_4({\bbd}) < 0, h_3({\bbd})\geq 0\}
& \{\bbd\in\ZZ^3\mid 0 >d_2 > d_1\geq d_3\}
\\ \hline 
\textbf{C14} 
&\{{\bbd} \in \ZZ^3 \mid h_1({\bbd}), h_2({\bbd}), h_3({\bbd}), h_4({\bbd}) < 0\} 
& -\text{Int}(\NN A) 
\\ \hline
\end{array}$
\caption{The chambers used to compute $\II(J)$ in Examples \ref{ex:3dinterior} and \ref{ex:3dbookwtail}.}
\label{tab:chambers}
\end{table}

In Figure~\ref{fig:x2yz}, the multidegrees in the gray cone whose vertex lies at $(2,1,1)$ correspond to the monomials that lie in $J$, and the monomials in the outer cone lie in $R_A\setminus J$.
Note also that the view  of the cone that we see in Figure~\ref{fig:x2yz} is from the side of the $xz$-plane. 

Since none of the points on the facets of the cone have monomials that lie in $J$, the 14 regions that we consider in determining $D(R_A,J)$ are listed in Table~\ref{tab:regions1}. 

If ${\bbd}$ is in \textbf{C1} (which are the lattice points corresponding to points in the semigroup) when considering the idealizer or \textbf{R1} (which are the lattice points corresponding to the monomials in $J$) when considering $D(R_A,J)$, then $D(R_A)_{\bbd}$ is generated by multiplication by ${\bbt^{\bbd}}$.  
Multiplying any element of $J$ by ${\bbt^{\bbd}}$ remains in $J$ making the graded pieces of $\II(J)_{\bbd}={\bbt^{\bbd}}\cdot\CC[\theta]$ and $D(R_A,J)_{\bbd}={\bbt^{\bbd}}\cdot\CC[\theta]$.

\begin{figure}[h]
\begin{subfigure}{1.90in}
    \centering
    \includegraphics[scale=.35]{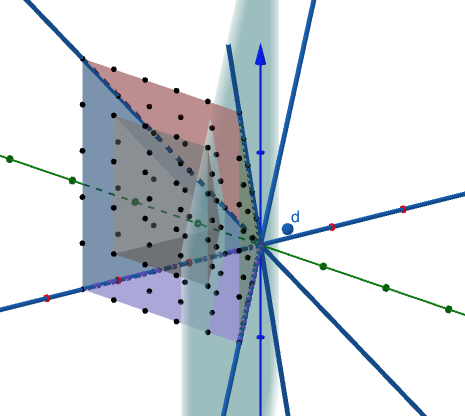}
 \captionsetup{margin=0in,width=2in,font=small,justification=centering,labelsep=newline}
   \caption{${\bbd} = (1,2,1)\in \textbf{C2}$ 
    and $\{x-y+h_4({\bbd}) = 0\}$}
    \label{fig:C2nodots}
\end{subfigure}
\begin{subfigure}{2.10in}
    \centering
    \includegraphics[scale=.3]{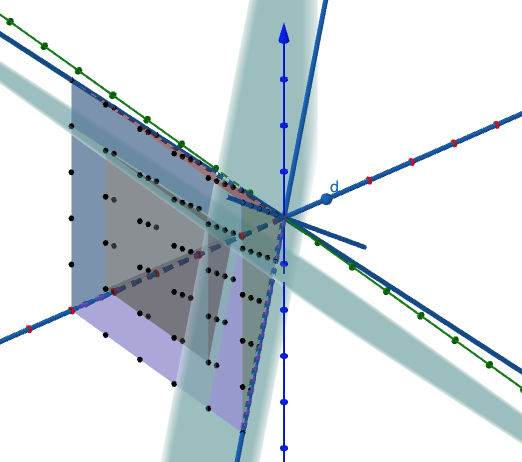}
\captionsetup{margin=0in,width=1.5in,font=small,justification=centering,labelsep=newline}
    \caption{${\bbd} = (-1,0,0)\in\textbf{C6}$, 
    $\{x-y+h_4({\bbd}) = 0\}$,  and $\{x - z +h_4({\bbd})=0\}$
    }
    \label{fig:x2yz:C6}
\end{subfigure}
\begin{subfigure}{2.40in}
    \centering
    \includegraphics[scale=.25]{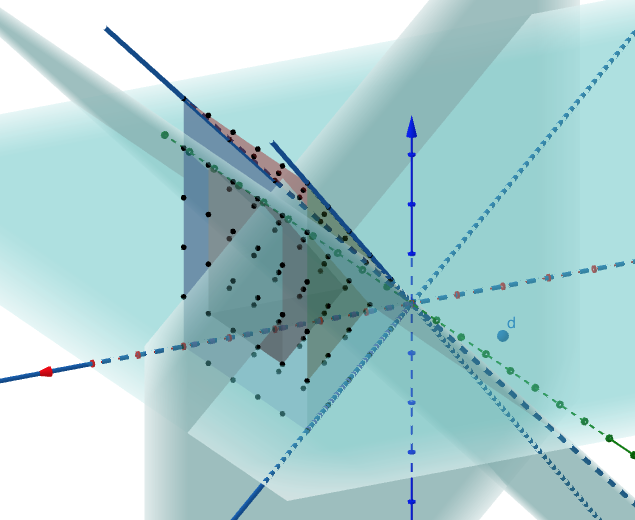}
\captionsetup{margin=0in,width=1.5in,font=small,labelsep=newline,justification=centering}
    \caption{${\bbd} = (-2,0,-1)\in\textbf{C10}$, 
    $\{z+h_1({\bbd}) = 0\}$,  $\{x-z+ h_3({\bbd})= 0\}$,  and $\{x - y +h_4({\bbd})=0\}$}
    \label{fig:C10nodots}
\end{subfigure}
\captionsetup{margin=0in,width=2.5in,font=small,justification=centering}
    \caption{Vanishing planes}
\end{figure}

\renewcommand{\arraystretch}{2}
\begin{table}
\scalebox{1.0}{
$\begin{array}{@{\extracolsep{\fill}} |c|l|l|}
\hline 
\text{Region} & \text{Halfspace inequalities} &  \text{Lattice point inequalities}
\\ \hline 
\textbf{R1} 
&\{ {\bbd} \in \ZZ^3 \mid h_1({\bbd}), h_2({\bbd}), h_3({\bbd}), h_4({\bbd}) > 0\}
& \text{Int}(\NN A) 
\\ \hline 
\textbf{R2} 
&\{ {\bbd} \in \ZZ^3 \mid h_1({\bbd}), h_2({\bbd}), h_3({\bbd}) > 0, h_4({\bbd})\leq 0\}
& \{\bbd\in\ZZ^3\mid d_2\geq d_1 > d_3 > 0\}
\\ \hline 
\textbf{R3} 
&\{ {\bbd} \in \ZZ^3 \mid h_1({\bbd}), h_2({\bbd}), h_4({\bbd}) > 0, h_3({\bbd})\leq 0\}
& \{\bbd\in\ZZ^3\mid d_3\geq d_1 > d_2 > 0\}
\\ \hline 
\textbf{R4} 
&\{ {\bbd} \in \ZZ^3 \mid h_1({\bbd}), h_3({\bbd}), h_4({\bbd}) > 0, h_2({\bbd})\leq 0\}
& \{\bbd\in\ZZ^3\mid d_1 > d_3> 0 \geq d_2\}
\\ \hline 
\textbf{R5} 
&\{ {\bbd} \in \ZZ^3 \mid h_2({\bbd}), h_3({\bbd}), h_4({\bbd}) > 0, h_1({\bbd}) \leq 0\}
& \{\bbd\in\ZZ^3\mid d_1 > d_2 > 0 \geq d_3\}
\\ \hline 
\textbf{R6} 
&\{ {\bbd} \in \ZZ^3 \mid h_1({\bbd}), h_2({\bbd}) > 0, h_3({\bbd}), h_4({\bbd})\leq 0\}
& \{\bbd\in\ZZ^3\mid d_2>0, d_3> 0, d_2\geq d_1, d_3 \geq d_1\}
\\ \hline 
\textbf{R7} 
&\{ {\bbd} \in \ZZ^3 \mid h_2({\bbd}), h_3({\bbd}) > 0, h_1({\bbd}), h_4({\bbd})\leq 0\}
& \{\bbd\in\ZZ^3\mid d_2> 0 \geq d_3, d_2 \geq d_1>d_3\}
\\ \hline 
\textbf{R8} 
&\{ {\bbd} \in \ZZ^3 \mid h_1({\bbd}), h_4({\bbd}) >0, h_2({\bbd}), h_3({\bbd}) \leq  0\}
& \{\bbd\in\ZZ^3\mid d_3> 0 \geq d_2,  d_3\geq  d_1> d_2\}
\\ \hline 
\textbf{R9} 
&\{ {\bbd} \in \ZZ^3 \mid h_3({\bbd}), h_4({\bbd}) > 0, h_1({\bbd}), h_2({\bbd})\leq 0\}
& \{\bbd\in\ZZ^3\mid d_1> d_2, d_1>d_3, 0 \geq d_2, 0\geq d_3\}
\\ \hline 
\textbf{R10} 
&\{ {\bbd} \in \ZZ^3 \mid h_1({\bbd}), h_3({\bbd}) \leq 0, h_4({\bbd}) < 0, h_2({\bbd})> 0\}
& \{\bbd\in\ZZ^3\mid d_2> 0  \geq d_3\geq d_1\}
\\ \hline 
\textbf{R11} 
&\{ {\bbd} \in \ZZ^3 \mid h_2({\bbd}), h_4({\bbd}) \leq 0, h_3({\bbd}) <  0, h_1({\bbd}) >0\}
& \{\bbd\in\ZZ^3\mid d_3> 0  \geq d_2\geq d_1\}
\\ \hline 
\textbf{R12} 
&\{ {\bbd} \in \ZZ^3 \mid h_1({\bbd}),  h_3({\bbd}) \leq 0, h_2({\bbd})<0, h_4({\bbd})> 0\}
& \{\bbd\in\ZZ^3\mid 0 \geq d_3 \geq d_1> d_2\}
\\ \hline 
\textbf{R13} 
&\{ {\bbd} \in \ZZ^3 \mid h_2({\bbd}), h_4({\bbd}) \leq  0, h_1({\bbd})<0, h_3({\bbd}) > 0\}
& \{\bbd\in\ZZ^3\mid 0 \geq d_2 \geq  d_1> d_3\}
\\ \hline 
\textbf{R14} 
&\{ {\bbd} \in \ZZ^3 \mid h_1({\bbd}), h_2({\bbd}), h_3({\bbd}), h_4({\bbd}) \leq  0\}
& -\NN A 
\\ \hline
\end{array}$}
\caption{Regions to compute $D(R_A,J)$ in Example \ref{ex:3dinterior}.}
\label{tab:regions1}
\end{table}

Since the multidegrees in the chambers and regions corresponding to $\textbf{C2}-\textbf{C5}$ and $\textbf{R2}-\textbf{R5}$ 
will all need operators adjusted only for monomials in $J$ whose exponents are parallel to a single facet in $A$, 
we will only describe the process to determine the graded pieces of $\II (J)$ and $D(R_A,J)$ for chamber \textbf{C2} and region \textbf{R2}, respectively, as precisely the same type of argument holds for the other chambers and regions.  
If ${\bbd}$ is in \textbf{C2} in the case of the idealizer or in \textbf{R2} in the case of $D(R_A,J)$, consider the lattice points corresponding to the monomials in $J$ on the plane $x-y+h_4({\bbd})=0$. 
When an element in $\II(\Omega({\bbd}))$ is applied to such a monomial in $J$, the result is a constant times a monomial with exponent on the facet $\sigma_4$, which is not in $J$. 
Hence we need to multiply $\II(\Omega({\bbd}))$ by $h_4+h_4({\bbd})$, so for ${\bbd}$ in \textbf{C2},  
\[
\II(J)_{\bbd}={\bbt^{\bbd}}\cdot\<\df{h_4}{h_4({-\bbd})}\>
\]
and  for ${\bbd}$ in \textbf{R2}, 
\[
D(R_A,J)_{\bbd}={\bbt^{\bbd}}\cdot\<\df{h_4}{h_4(-{\bbd})}\>.
\]
Figure~\ref{fig:C2nodots} illustrates the  plane in \textbf{C2} that determines the linear form that we multiply by $\II(\Omega(\bbd))$ to obtain $\II(J)_{\bbd}$.

Since the multidegrees in the chambers and regions corresponding to $\textbf{C6}-\textbf{C9}$ and $\textbf{R6}-\textbf{R6}$ 
will all need operators adjusted only for monomials in $J$ whose exponents are parallel to two facets of $A$, 
we will only describe the process to determine the graded pieces of $\II (J)$ and $D(R_A,J)$ for chamber \textbf{C6} and region \textbf{R6}, respectively, as precisely the same type of argument holds for the other chambers and regions. 
If ${\bbd}$ is in \textbf{C6} in the case of the idealizer or \textbf{R6} in the case of $D(R_A,J)$, consider the lattice points corresponding to the monomials in $J$ on the planes $x-y+h_4({\bbd})=0$ and $x-z+h_3({\bbd})=0$. 
When an element in $\II(\Omega({\bbd}))$ is applied to such monomial in $J$, the result is a constant times a monomial with exponent on $x-y=0$ or $x-z=0$,  which is not in $J$. 
Hence we need to multiply $\II(\Omega({\bbd}))$ by $(h_3+h_3({\bbd}))(h_4+h_4({\bbd}))$, so for ${\bbd}$ in \textbf{C6},   
\[
\II(J)_{\bbd}={\bbt^{\bbd}}\cdot\<\df{h_3}{h_3({-\bbd})}\df{h_4}{h_4({-\bbd})}\>,
\]
and for ${\bbd}$ in \textbf{R6},  
\[
D(R_A,J)_{\bbd}={\bbt^{\bbd}}\cdot\<\df{h_3}{h_3(-{\bbd})}\df{h_4}{h_4({-\bbd})}\>.
\]
Figure~\ref{fig:x2yz:C6} illustrates the two planes in \textbf{C6} that determine the linear forms that we multiply $\II(\Omega(\bbd))$ by to obtain $\II(J)_{\bbd}$.

Since the multidegrees in the chambers and regions corresponding to $\textbf{C10}-\textbf{C13}$ and $\textbf{R10}-\textbf{R13}$ will all need operators adjusted only for monomials in $J$ whose exponents are parallel to three facets of $A$, we will only describe the process to determine the graded piece of $\II (J)$ and $D(R_A,J)$ for chamber \textbf{C10} and region \textbf{R10}, respectively, as precisely the same type of argument holds for the other chambers and regions.  
If ${\bbd}$ is in \textbf{C10} in the case of the idealizer or \textbf{R10} in the case of $D(R_A,J)$, consider the lattice points corresponding to the monomials in $J$ on the planes $z+h_1({\bbd})=0$, $x-z+h_3({\bbd})=0$ and $x-y+h_4({\bbd})=0$. 
When an element of $\II(\Omega({\bbd}))$ is applied to such a monomial in $J$, 
the result is a constant times a monomial whose exponent is in $\sigma_1$, $\sigma_3$, or $\sigma_4$, 
which is not in $J$. 
Hence we need to multiply $\II(\Omega({\bbd}))$ by $(h_1+h_1({\bbd}))(h_3+h_3({\bbd}))(h_4+h_4({\bbd}))$, so for ${\bbd}$ in \textbf{C10}, 
\[
\II(J)_{\bbd}={\bbt^{\bbd}}\cdot\<\df{h_1}{h_1({-\bbd})}\df{h_3}{h_3(-{\bbd})}\df{h_4}{h_4(-{\bbd})}\>,
\]
and for ${\bbd}$ in \textbf{R10}, 
\[
D(R_A,J)_{\bbd}={\bbt^{\bbd}}\cdot\<\df{h_1}{h_1({-\bbd})}\df{h_3}{h_3({-\bbd})}\df{h_4}{h_4(-{\bbd})}\>.
\]
Figure~\ref{fig:C10nodots} illustrates the three planes in \textbf{C10} that determine the linear forms that we multiply by $\II(\Omega(\bbd))$ to obtain $\II(J)_{\bbd}$.

If ${\bbd}$ is in \textbf{C14} in the case of the idealizer or \textbf{R14} in the case of $D(R_A,J)$, consider the lattice points corresponding to the monomials in $J$ on the planes $z+h_1({\bbd})=0$, $y+h_2({\bbd})=0$, $x-z+h_3({\bbd})=0$ and $x-y+h_4({\bbd})=0$. 
When an element of $\II(\Omega({\bbd}))$ is applied to such a monomial in $J$, the result is a constant times a monomial with exponent on one of the facets $\sigma_i$ of $A$, which is not in $J$.  
Hence we need to multiply $\II(\Omega({\bbd}))$ by $(h_1+h_1({\bbd}))(h_2+h_2({\bbd}))(h_3+h_3({\bbd}))(h_4+h_4({\bbd}))$, so for ${\bbd}$ in \textbf{C14},  
\[
\II(J)_{\bbd}={\bbt^{\bbd}}\cdot\<\df{h_1}{h_1({-\bbd})}\df{h_2}{h_2({-\bbd})}\df{h_3}{h_3({-\bbd})}\df{h_4}{h_4({-\bbd})}\>,
\]
and for ${\bbd}$ in \textbf{R14}, 
\[
D(R_A,J)_{\bbd}={\bbt^{\bbd}}\cdot\<\df{h_1}{h_1({-\bbd})}\df{h_2}{h_2({-\bbd})}\df{h_3}{h_3({-\bbd})}\df{h_4}{h_4({-\bbd})}\>.
\]

Combining the information from all $14$ chambers \textbf{C1-C14} and regions  \textbf{R1-R14}, the general formula for the graded piece of the idealizer of $J$ at ${\bbd}\in\ZZ^3$ is 
\[
\II(J)_{\bbd} 
= {\bbt^{\bbd}} \cdot 
\left\<\prod\limits_{h_i({\bbd})<0}\df{h_i}{h_i(-{\bbd})}\right\>,
\] 
and the general formula for the graded piece of the $D(R_A,J)$ at ${\bbd}\in\ZZ^3$ is 
\[
D(R_A,J)_{\bbd} 
= {\bbt^{\bbd}} \cdot 
\left\<\prod\limits_{h_i({\bbd})\leq 0}\df{h_i}{h_i(-{\bbd})}\right\>.
\]  
For any $\bbd\in\ZZ^3$, 
\[
\left(\dfrac{\II(J)}{D(R_A,J)}\right)_{\bbd} 
= {\bbt^{\bbd}}\cdot 
\displaystyle\frac{\left\<\prod\limits_{h_i({\bbd})<0}\df{h_i}{h_i(-{\bbd})}\right\>}{\left\<\prod\limits_{h_i({\bbd})\leq 0}\df{h_i}{h_i(-{\bbd})}\right\>}.
\]

\begin{figure}[h]
    \centering
    \includegraphics[scale=.4]{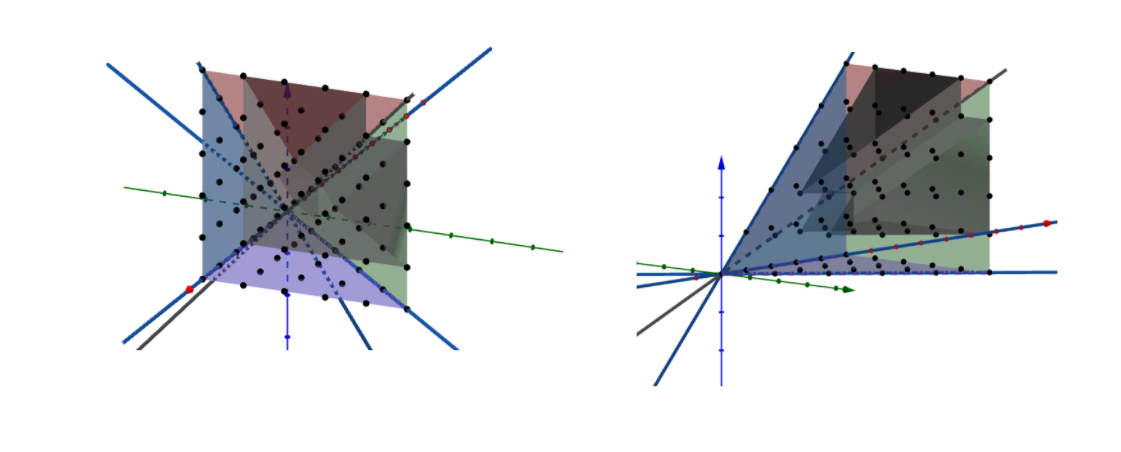}
 \captionsetup{margin=0in,width=2.5in,font=small,justification=centering,labelsep=newline}
     \caption{Cone of $R_A$ with $J=P_{\sigma_1}\cap P_{\sigma_2} \cap P_{\sigma_3 \cap \sigma_4}$}
    \label{fig:higherD:second-cone}
\end{figure}

\begin{ex}
\label{ex:3dbookwtail}
Now consider the ideal 
\[
J= P_{\sigma_1} \cap P_{\sigma_2} \cap P_{\sigma_3 \cap \sigma_4}=\<t_1^2t_2t_3, t_1^2t_2^2t_3,t_1^2t_2t_3^2\>
\] 
in $R_A = \CC[t_1,t_1t_2,t_1t_3,t_1t_2t_3]$.
In Figure~\ref{fig:higherD:second-cone}, we give two views of the cone with the lattice points in $J$ shaded in grey.
Note that multidegrees of the monomials in $J$ lie in the interior of the cone and on the interiors of the two faces $\sigma_3$ and $\sigma_4$, so this is the grey portion of the figure. 
\end{ex}
 
Since the only points on the facets of the cone that have monomials that lie in $J$ lie in the interiors of $\sigma_3$ and $\sigma_4$, the 14 \emph{regions} that we consider in determining $D(R_A,J)$ appear in Table~\ref{tab:regions2}. 

\renewcommand{\arraystretch}{2}
\begin{table}[ht]
\hspace{-.1in}
\centering
\scalebox{0.88}{
$\begin{array}{@{\extracolsep{\fill}} |c|l|l|}
\hline 
\text{Region} & \text{Halfspace inequalities} &  \text{Lattice point inequalities}
\\ \hline 
\textbf{R1} 
& 
\begin{array}{l}
    \{ {\bbd} \in \ZZ^3 \mid h_1({\bbd}), h_2({\bbd}), h_3({\bbd})>0, h_4({\bbd}) \geq 0\}\ \cup \vspace{-6pt} \\
    \qquad \{{\bbd} \in \ZZ^3 \mid h_1({\bbd}), h_2({\bbd}), h_4({\bbd})>0, h_3({\bbd}) \geq 0\}
\end{array}
& 
\begin{array}{l}
    \{\bbd\in\ZZ^3\mid {\bbx}^{\bbd} \in J \}
\end{array}
\\ \hline 
\textbf{R2} 
&\{ {\bbd} \in \ZZ^3 \mid  h_1({\bbd}), h_2({\bbd})>0, h_3({\bbd}) \geq  0, h_4({\bbd}) < 0\}
& \{\bbd\in\ZZ^3\mid d_2>  d_1 \geq  d_3 > 0\}
\\ \hline 
\textbf{R3} 
&\{ {\bbd} \in \ZZ^3 \mid  h_1({\bbd}), h_2({\bbd})> 0, h_4({\bbd}) \geq 0, h_3({\bbd}) < 0\}
& \{\bbd\in\ZZ^3\mid d_3 >  d_1 \geq d_2 > 0\}
\\ \hline 
\textbf{R4} 
&
\begin{array}{l}
    \{ {\bbd} \in \ZZ^3 \mid h_1({\bbd}), h_3({\bbd}), h_4({\bbd})>0,  h_2({\bbd})\leq 0\}\ \vspace{-6pt} \cup \\
    \qquad \{{\bbd} \in \ZZ^3 \mid h_3({\bbd})=0, h_1({\bbd}),h_4({\bbd})>0, h_2({\bbd})<0\}
\end{array}
& 
\begin{array}{l}
    \{\bbd\in\ZZ^3\mid d_1 > d_3> 0 \geq d_2\}\ \vspace{-6pt} \cup \\
    \qquad \{\bbd\in\ZZ^3\mid d_1=d_3>0>d_2\}
\end{array}
\\ \hline 
\textbf{R5} 
&
\begin{array}{l}
    \{ {\bbd} \in \ZZ^3 \mid  h_2({\bbd}), h_3({\bbd}), h_4({\bbd})>0, h_1({\bbd}) \leq 0\}\ \vspace{-6pt} \cup \\
    \qquad \{{\bbd} \in \ZZ^3 \mid h_4({\bbd})=0, h_2({\bbd}),h_3({\bbd})>0, h_1({\bbd})<0\}
\end{array}
& 
\begin{array}{l}
    \{\bbd\in\ZZ^3\mid d_1 \geq d_2 > 0 \geq d_3\}\ \vspace{-6pt} \cup\\
    \qquad \{\bbd\in\CC^3\mid d_1=d_2>0>d_3\}
\end{array}
\\ \hline 
\textbf{R6} 
&
\begin{array}{l}
    \{ {\bbd} \in \ZZ^3 \mid h_1({\bbd}), h_2({\bbd}) > 0\text{ and } h_3({\bbd}), h_4({\bbd})<0\}\ \vspace{-6pt} \cup \\
    \qquad  \{{\bbd} \in \ZZ^3 \mid h_1({\bbd})>0, h_2({\bbd}) > 0, h_3({\bbd})=h_4({\bbd})=0\}
\end{array}
& 
\begin{array}{l}
    \{\bbd\in\ZZ^3\mid d_3>0, d_2> 0, d_2>d_1, d_3 >d_1\}\ \vspace{-6pt} \cup\\
    \qquad \qquad \qquad \qquad \{\bbd\in\ZZ^3\mid d_1=d_2=d_3>0\}
\end{array}
\\ \hline 
\textbf{R7} 
&
\begin{array}{l}
    \{ {\bbd} \in \ZZ^3 \mid h_2({\bbd}), h_3({\bbd}) > 0, h_1({\bbd}) \leq 0, h_4({\bbd})< 0\}\ \vspace{-6pt}\cup \\
    \qquad \qquad \quad \{{\bbd} \in \ZZ^3 \mid h_2({\bbd})>0,h_3({\bbd})=0, h_1({\bbd})<0\}
\end{array}
& 
\begin{array}{l}
    \{\bbd\in\ZZ^3\mid d_2> 0 \geq d_3, d_2 >d_1>d_3\}\ \vspace{-6pt} \cup\\
    \qquad \qquad \qquad \{\bbd\in\ZZ^3\mid d_2>0 >d_1=d_3\}
\end{array}
\\ \hline 
\textbf{R8} 
&
\begin{array}{l}
    \{ {\bbd} \in \ZZ^3 \mid h_1({\bbd}), h_4({\bbd}) > 0, h_2({\bbd}) \leq 0, h_3({\bbd})<0\}\ \vspace{-6pt}\cup \\
    \qquad \qquad \quad \{{\bbd} \in \ZZ^3 \mid  h_1({\bbd})>0,h_4({\bbd})=0, h_2({\bbd})<0\}
\end{array}
& 
\begin{array}{l}
    \{\bbd\in\ZZ^3\mid d_3> 0 \geq d_2, d_3 >d_1>d_2\}\ \vspace{-6pt} \cup\\ 
    \qquad \qquad \qquad \{\bbd\in\ZZ^3\mid d_3>0 >d_1=d_2\}
\end{array}
\\ \hline 
\textbf{R9} 
&
\begin{array}{l}
    \{ {\bbd} \in \ZZ^3 \mid h_3({\bbd}), h_4({\bbd}) > 0, h_1({\bbd}), h_2({\bbd})\leq 0\}\ \vspace{-6pt}\cup \\
    \quad \{ {\bbd} \in \ZZ^3 \mid h_1({\bbd}), h_2({\bbd})<0,h_4({\bbd})>0,h_3({\bbd})=0\}\ \vspace{-6pt} \cup \\
    \qquad \quad \{\bbd\in\ZZ^3\mid h_1({\bbd}), h_2({\bbd}<0, h_3({\bbd})>0,h_4({\bbd})=0\}
\end{array}
& 
\begin{array}{l}
    \{\bbd\in\ZZ^3\mid d_1 > d_2,d_3, 0 \geq d_2,d_3\}\ \vspace{-6pt}\cup \\
    \qquad \{\bbd\in\ZZ^3\mid 0>d_1=d_2, 0>d_1>d_3\}\ \vspace{-6pt} \cup \\
    \qquad \qquad \{\bbd\in\ZZ^3\mid 0>d_1=d_3, 0>d_1>d_2\}
\end{array}
\\ \hline 
\textbf{R10} 
&\{ {\bbd} \in \ZZ^3 \mid  h_1({\bbd})\leq 0, h_3({\bbd}), h_4({\bbd}) < 0, h_2({\bbd})> 0\}
& \{\bbd\in\ZZ^3\mid d_2> 0  \geq d_3> d_1\}
\\ \hline 
\textbf{R11} 
&\{ {\bbd} \in \ZZ^3 \mid h_2({\bbd}) \leq 0, h_4({\bbd}), h_3({\bbd}) <  0, h_1({\bbd}) >0\}
& \{\bbd\in\ZZ^3\mid d_3> 0  \geq d_2>d_1\}
\\ \hline 
\textbf{R12} 
&\{ {\bbd} \in \ZZ^3 \mid  h_1({\bbd}), h_3({\bbd})<0, h_2({\bbd}) \leq  0, h_4({\bbd}) \geq  0 \}
& \{\bbd\in\ZZ^3\mid 0 \geq d_2 > d_1\geq  d_3\}
\\ \hline 
\textbf{R13} 
&\{ {\bbd} \in \ZZ^3 \mid h_2({\bbd}), h_4({\bbd})<0, h_1({\bbd}) \leq  0, h_3({\bbd}) \geq  0\}
& \{\bbd\in\ZZ^3\mid 0 \geq d_3 > d_1\geq  d_2\}
\\ \hline 
\textbf{R13} 
&\{ {\bbd} \in \ZZ^3 \mid  h_1({\bbd}), h_2({\bbd}), h_4({\bbd}) \leq  0, h_3({\bbd}) > 0\}
& \{\bbd\in\ZZ^3\mid 0 \geq d_2 \geq  d_1> d_3\}
\\ \hline
\textbf{R14} 
&\begin{array}{l} \{ {\bbd} \in \ZZ^3 \mid h_1({\bbd}), h_2({\bbd}) h_3({\bbd})\leq 0, h_4({\bbd}) <  0\} \ \vspace{-6pt} \cup \\
   \qquad \{ {\bbd} \in \ZZ^3 \mid h_1({\bbd}), h_2({\bbd}) h_4({\bbd})\leq 0, h_3({\bbd}) <  0\}
\end{array}
& \begin{array}{l}\{\bbd\in\ZZ^3\mid d_1< d_2\leq 0, d_1\leq  d_3\leq 0\}\ \vspace{-6pt} \cup \\
    \qquad  \{\bbd\in\ZZ^3\mid d_1\leq  d_2\leq 0, d_1 <  d_3\leq 0\}
\end{array}
\\ \hline
\end{array}$}
\caption{Regions for $D(R_A,J)$ in Example~\ref{ex:3dbookwtail}. 
}
\label{tab:regions2}
\end{table}

If ${\bbd}$ is in $\textbf{C1}=\NN A$, when considering the idealizer or \textbf{R1}, which are the exponents of the monomials in $J$, when considering $D(R_A,J)$, then $D(R_A)_{\bbd}$ is generated by ${\bbt^{\bbd}}$. 
The product of any element of $J$ and ${\bbt^{\bbd}}$ is also an element of $J$, which means that the graded pieces are $\II(J)_{\bbd}={\bbt^{\bbd}}\cdot\CC[\theta]$ and $D(R_A,J)_{\bbd}={\bbt^{\bbd}}\cdot\CC[\theta]$.

Since the multidegrees in $\textbf{C2}-\textbf{C5}$ and $\textbf{R2}-\textbf{R5}$ potentially affect monomials in $J$ whose exponents are parallel to a single facet of $A$, we will only describe the process to determine the graded pieces of $\II (J)$ and $D(R_A,J)$ for chamber \textbf{C2} and region \textbf{R2}, respectively, as precisely the same type of argument holds for the other chambers and regions. 
The reader may refer to Figure~\ref{fig:C2nodots} to help visualize the argument below. 
If ${\bbd}$ is in \textbf{C2} in the case of the idealizer or \textbf{R2} in the case of $D(R_A,J)$, consider exponents of monomials in $J$ on the plane $x-y+h_4({\bbd})=0$. 
When an element of $\II(\Omega({\bbd}))$ is applied to such a monomial in $J$, the result is a constant times a monomial on the interior of the facet $\sigma_4$, which is in $J$.  
Hence, for ${\bbd}$ in \textbf{C2}, 
\[
\II(J)_{\bbd}={\bbt^{\bbd}}\cdot\<\df{h_4}{h_4({-\bbd})-1}\>,
\]
and for ${\bbd}$ in \textbf{R2}, 
\[
D(R_A,J)_{\bbd}={\bbt^{\bbd}}\cdot\<\df{h_4}{h_4({-\bbd})-1}\>.
\]

Since the multidegrees in the chambers and regions corresponding to $\textbf{C6}-\textbf{C9}$ and $\textbf{R6}-\textbf{R9}$ can 
will all need operators adjusted only for monomials in $J$ whose exponents are parallel to two facets of $A$, we will only describe the process to determine the graded piece of $\II (J)$ and $D(R_A,J)$ for chamber \textbf{C6} and region \textbf{R6}, respectively, as precisely the same type of argument holds for the other chambers and regions.  
The reader may refer to Figure~\ref{fig:x2yz:C6} to help visualize the argument below.  
If ${\bbd}$ is in \textbf{C6} in the case of the idealizer or \textbf{R6} in the case of $D(R_A,J)$, consider the multidegrees of monomials in $J$ on the planes $x-z+h_3({\bbd})=0$ or $x-y+h_4({\bbd})=0$. 
When an element of $\II(\Omega({\bbd}))$ is applied to such a monomial in $J$, the result is a constant times a monomial on $\sigma_3$ or $\sigma_4$. 
This is not in $J$ only if this monomial's exponent is in the intersection of the two planes. 
Hence, we need to multiply $\II(\Omega({\bbd}))$ by either $(h_3+h_3({\bbd}))$ or $(h_4+h_4({\bbd}))$,  so that  for ${\bbd}$ in \textbf{C6}, 
\[
\II(J)_{\bbd} 
= {\bbt^{\bbd}}\cdot 
\df{h_4}{h_4({-\bbd})-1}\df{h_3}{h_3({-\bbd})-1} \cdot 
\< (h_3+h_3({\bbd})), (h_4+h_4({\bbd})) \>,
\] 
and for ${\bbd}$ in \textbf{R6}, 
\[
D(R_A,J)_{\bbd} = 
{\bbt^{\bbd}}\cdot 
\df{h_4}{h_4({-\bbd})-1}\df{h_3}{h_3({-\bbd})-1}\cdot 
\< (h_3+h_3({\bbd})), (h_4+h_4({\bbd})) \>.
\] 

Since the multidegrees in the chambers and regions corresponding to $\textbf{C10}-\textbf{C13}$ and $\textbf{R10}-\textbf{R13}$ can will all need operators adjusted only for monomials in $J$ whose exponents are parallel to a  three facets of $A$, we will only describe the process to determine the graded piece of $\II (J)$ and $D(R_A,J)$ for chamber \textbf{C10} and region \textbf{R10}, respectively, as precisely the same type of argument holds for the other chambers and regions. 
The reader may refer to Figure~\ref{fig:C10nodots} to help visualize the argument below. 
If ${\bbd}$ is in \textbf{C10} in the case of the idealizer or \textbf{R10} in the case of $D(R_A,J)$, consider the multidegrees of monomials in $J$ on the planes $z+h_1({\bbd})=0$, $x-z+h_3({\bbd})=0$, and $x-y+h_4({\bbd})=0$. 
When an element in $\II(\Omega({\bbd}))$ is applied to such a monomial in $J$, the result is 
a constant times a monomial on $\sigma_1$, which is not in $J$, 
or 
a constant times a monomial in $\sigma_3$ or $\sigma_4$, 
which is not in $J$ 
only if the monomial is in the intersection of $\sigma_3$ and $\sigma_4$.  
Hence, we need to multiply $\II(\Omega({\bbd}))$ by either $(h_1+h_1({\bbd}))(h_3+h_3({\bbd}))$ or $(h_1+h_1({\bbd}))(h_4+h_4({\bbd}))$, so that for ${\bbd}$ in \textbf{C10},  
\[
\II(J)_{\bbd} 
= {\bbt^{\bbd}}\cdot 
\df{h_1}{h_1({-\bbd})}\df{h_3}{h_3({-\bbd})-1}\df{h_4}{h_4({-\bbd})-1}\cdot 
\< (h_3+h_3({\bbd})), (h_4+h_4({\bbd}))\>, 
\]
and for ${\bbd}$ in \textbf{R10}, 
\[
D(R_A,J)_{\bbd} 
= {\bbt^{\bbd}}\cdot 
\df{h_1}{h_1({-\bbd})}\df{h_3}{h_3({-\bbd})-1}\df{h_4}{h_4({-\bbd})-1}\cdot 
\< (h_3+h_3({\bbd})), (h_4+h_4({\bbd}))\>. 
\] 

If ${\bbd}$ is in \textbf{C14} in the case of the idealizer or \textbf{R14} in the case of $D(R_A,J)$, consider the multidegrees of the monomials in $J$ on the planes $z+h_1({\bbd})=0$, $y+h_2({\bbd})=0$, $x-z+h_3({\bbd})=0$, and $x-y+h_4({\bbd})=0$. 
When an element of $\II(\Omega({\bbd}))$ is applied to such a monomial in $J$, the result is 
a constant times a monomial in $\sigma_1$ or $\sigma_2$, 
which is not in $J$, 
or 
a constant times a monomial in $\sigma_3$ or $\sigma_4$, 
which is not in $J$ if this monomial lies in both planes parallel to these faces.
Hence we need to multiply $\II(\Omega({\bbd}))$ by $(h_1+h_1({\bbd}))(h_2+h_2({\bbd}))(h_3+h_3({\bbd}))(h_4+h_4({\bbd})-1)$ or $(h_1+h_1({\bbd}))(h_2+h_2({\bbd}))(h_3+h_3({\bbd})-1)(h_4+h_4({\bbd}))$, so that  for ${\bbd}$ in \textbf{C14}, 
\[
\II(J)_{\bbd} 
= {\bbt^{\bbd}} \cdot 
\prod\limits_{i=1}^4\df{h_i}{h_i({-\bbd})-1} \cdot 
\left\<\left( \prod\limits_{i=1}^3(h_i+h_i({\bbd}))\right), (h_4+h_4({\bbd}))\left(\prod\limits_{i=1}^2(h_i+h_i({\bbd}))\right)\right\>
\] 
and for ${\bbd}$ in \textbf{R14}, 
\[
D(R_A,J)_{\bbd} 
= {\bbt^{\bbd}} \cdot 
\prod\limits_{i=1}^4 \df{h_i}{h_i({-\bbd})-1}\cdot 
\left\<\left( \prod\limits_{i=1}^3(h_i+h_i({\bbd}))\right), 
(h_4+h_4({\bbd}))\left(\prod\limits_{i=1}^2(h_i+h_i({\bbd}))\right)\right\>. 
\] 

Putting all the information together from all $14$ chambers \textbf{C1-C14} and regions \textbf{R1-R14}, the general formula for the graded piece of the idealizer of $J$ at ${\bbd}\in\ZZ^3$ is 
\[
\II(J)_{\bbd} 
= {\bbt^{\bbd}}\cdot 
\prod\limits_{i=1}^4\df{h_i}{h_i({-\bbd})-1}\cdot 
\left\<\left( \prod\limits_{\substack{h_i({\bbd})<0\\ \text{for } i \neq 4}}(h_i+h_i({\bbd}))\right), 
\left(\prod\limits_{\substack{h_i({\bbd})<0\\ \text{for } i \neq 3}}(h_i+h_i({\bbd}))\right)\right\>, 
\] 
and the general formula for the graded piece of $D(R_A,J)$ at ${\bbd}\in\ZZ^3$ is 
\[
D(R_A,J)_{\bbd} 
= {\bbt^{\bbd}}\cdot 
\prod\limits_{i=1}^4\df{h_i}{h_i({-\bbd})-1}\cdot 
\left\<\left( \prod\limits_{\substack{h_i({\bbd})\leq 0\\ 
\text{for } i \neq 4}}(h_i+h_i({\bbd}))\right), 
\left(\prod\limits_{\substack{h_i({\bbd}) \leq 0\\ 
\text{for } i \neq 3}}(h_i+h_i({\bbd}))\right)\right\>.
\] 
We also have for any $\bbd\in\ZZ^3$, 
\[
\left(\dfrac{\II(J)}{D(R_A,J)}\right)_{\bbd} 
= {\bbt^{\bbd}}\cdot 
\displaystyle\frac{\ 
    \prod\limits_{i=1}^4\df{h_i}{h_i({-\bbd})-1}\cdot 
    \left\<\left( \prod\limits_{\substack{h_i({\bbd})<0\\ 
    \text{for } i \neq 4}}(h_i+h_i({\bbd}))\right), 
    \left(\prod\limits_{\substack{h_i({\bbd})<0\\ 
    \text{for } i \neq 3}}(h_i+h_i({\bbd}))\right)\right\>
\ }{\ 
    \prod\limits_{i=1}^4\df{h_i}{h_i({-\bbd})-1}\cdot 
    \left\<\left( \prod\limits_{\substack{h_i({\bbd})\leq 0\\ 
    \text{for } i \neq 4}}(h_i+h_i({\bbd}))\right), 
    \left(\prod\limits_{\substack{h_i({\bbd}) \leq 0\\ 
    \text{for } i \neq 3}}(h_i+h_i({\bbd}))\right)\right\>
\ }. 
\]

\section{Non-normal examples}
\label{sec:nonnormal}

Thus far, we have restricted our attention to normal semigroup rings. 
However, Saito and Traves determined the ring of differential operators for all saturated affine semigroup rings 
in \cite[Theorem 3.3.1]{Sai-Tr-DASR} (see also \cite[Theorem 2.1]{Sai-Tr-fgdo}),
and we can use this work to compute $\II(J)$, $D(R_A,J)$ and $\II(J)/D(R_A,J)$ for non-normal $R_A$ and a graded radical $R_A$-ideal $J$. 
We broaden slightly from normal to \emph{scored} semigroup rings, see Definition~\ref{def:normal-scored}.   
We include two examples; one is scored and the other is not.  
However, the two quotient rings modulo the interior ideals are isomorphic and we will see that the expressions for $\II(J)/D(R_A,J)$ in both rings are isomorphic.  
\begin{ex}
\label{ex:one-line-gone}
Consider the matrix  
$A=\begin{bmatrix} 
    2 & 3 & 0\\ 
    0 & 0 & 1
\end{bmatrix}$, 
which is associated to the semigroup ring $R_A=\CC[\NN A]=\CC[s^2,s^3,t]$. 
Since $R_A$ is not normal, note that integer points in the faces of $A$ and $\NN A$ differ from those in $\RR_{\geq 0}A$. 
We let 
\[
\sigma_1=\NN\{ e_2\}
\quad\text{and}\quad
\sigma_2=\NN\{ e_1\}
\]
be the integer points in the facets of the cone $\RR_{\geq 0} A$. 
The primitive integral support functions of $\NN A$ are 
\[
h_1=h_1(\theta) = \theta_1
\quad\text{and}\quad
h_2= h_2(\theta) =\theta_2.
\]
The prime ideals associated to the facets of $A$ are 
\[
P_1=\<s^2,s^3\>
\quad\text{and}\quad
P_2=\<t\>.
\]
Set 
$J=P_1 \cap P_2=\<s^2t, s^3t\>$. 
Note that $R_A$ is a Gorenstein ring that is not normal, and, in this case, the ideal $J$, which is the interior, is not the canonical module of $R_A$.  
\end{ex}

\begin{wrapfigure}[13]{l}[10pt]{2in}
\centering
\begin{tikzpicture}[scale=0.55]
\draw[blue, thick] (-3.5,0) -- (3.5,0);
\draw[blue, thick] (0,-3.5) -- (0,3.5);
\draw[orange, thick] (1,0) -- (1,3.5);
\draw[orange, thick] (-1,0) -- (-1,3.5);
\draw[teal, thick] (1,0) -- (1,-3.5);
\draw[teal, thick] (-1,0) -- (-1,-3.5);
\foreach \y in {1,2,3}{ \node[draw,circle,inner sep=2pt,red,fill] at (2,\y) {};}
\foreach \y in {1,2,3}{ \node[draw,circle,inner sep=2pt,red,fill] at (3,\y) {};}
\foreach \y in {0,1,2,3}{ \node[draw,circle,inner sep=2pt,orange,fill] at (1,\y) {};}
\foreach \y in {0,1,2,3}{ \node[draw,circle,inner sep=2pt,orange,fill] at (-1,\y) {};}
\foreach \x in {0,2,3}{ \node[draw,circle,inner sep=2pt,blue,fill] at (\x,0) {};}
\foreach \y in {1,2,3}{ \node[draw,circle,inner sep=2pt,blue,fill] at (0,\y) {};}
\foreach \x in {-1,1}{
     \foreach \y in {-3,-2,-1}{
       \node[draw,circle,inner sep=2pt,teal,fill] at (\x,\y) {};}}

\foreach \x in {-3,-2,...,3}{
     \foreach \y in {-3,-2,...,3}{
       \node[draw,circle,inner sep=2pt] at (\x,\y) {};}}
\foreach \x in {-3,-2,...,-2}{
     \foreach \y in {0,1,2,3}{
       \node[draw,circle,inner sep=2pt,yellow,fill] at (\x,\y) {};}}
\foreach \x in {-3,-2}{
     \foreach \y in {-3,-2,-1}{
       \node[draw,circle,inner sep=2pt,violet,fill] at (\x,\y) {};}}     
\foreach \x in {0,2,3}{
     \foreach \y in {-3,-2,-1}{
       \node[draw,circle,inner sep=2pt,green,fill] at (\x,\y) {};}}
\end{tikzpicture}
\captionsetup{margin=0in,width=1.5in,font=small,justification=centering}
\caption{\\Chambers for $\CC[s^2,s^3,t]$}
\label{fig:nonnormal:lattice}
\end{wrapfigure}

We will again use \emph{chambers} of $\ZZ^2$ to aid in our computations of the graded pieces of 
$\II(J)/D(R_A,J)$;
in particular, we show that 
$\II(J)/D(R_A,J)$ is only nonzero in multidegree ${\bbd}$ if ${\bbd} \in \sigma_1 \cup \sigma_2$, as occurred in the normal cases computed thus far.
In Figure~\ref{fig:nonnormal:lattice}, the red multidegrees correspond to the monomials in $J$, and the blue multidegrees are exponents for monomials in the semigroup ring $R_A$ that do not lie in $J$. 
The monomials at the orange multidegrees behave somewhat more like the red and blue, 
in that, for all such $\bbd$, 
$\II(J)_{\bbd} = D(R_A)_{\bbd}$. 
Thus, together the red, blue, and orange multidegrees form chamber \textbf{C1}. 
The yellow multidegrees constitute \textbf{C2}, the violet points form \textbf{C3}, and both the green and dark green make \textbf{C4}. 
The vertical lines $\{x=-1\}$ and $\{x=1\}$ contain multidegrees $\bbd$ for which extra care is needed to determine $\II(J)_{\bbd}$ and $D(R_A,J)_{\bbd}$.

Eriksen showed in \cite[Proposition 2.1]{Er} that the ring of differential operators of a numerical semigroup ring $R=\CC[t^{\Gamma}]$ is 
\[
D(R) 
= \bigoplus\limits_{d \in \ZZ} 
s^{d_1}t^{d_2} \cdot \left\< \prod\limits_{\gamma \in \Omega(d)} (\theta-\gamma) \right\>, 
\]
where $\Omega(d)=\{\gamma \in \Gamma \mid \gamma +d \notin \Gamma\}$. 
We encode these polynomial generators using  $G_{\bbd}(\theta)$, where 
\[
G_{\bbd}(\theta) := 
\begin{cases}
1 
&\text{ if } d_1 \in \NN\< 2,3\>,\\
h_1 
&\text{ if } d_1=1,\\
h_1-1 
&\text {if } d_1=-1,\\
\displaystyle\frac{\df{h_1}{1-d_1}}{(h_1-1)(h_1+d_1)} 
&\text{ if } d_1\leq -2.
\end{cases} 
\]
By Theorem~\ref{thm:ST-dasr-3.2.2}, for any $\bbd\in\ZZ^2$, 
\[
D(R_A)_{\bbd} = 
s^{d_1}t^{d_2}\cdot
\left\< G_{\bbd}(\theta)\df{h_2}{-d_2-1}\right\>.
\]
To determine $\II(J)$ or $D(R_A,J)$, we must determine the lines of multidegrees that contain operators that, when applied to monomials in $J$ or $R$, produce elements in  $R\setminus J$. 
Note that given $\bbm\in\NN A$, 
\begin{align}
\label{eq:apply-diff-op-023}
s^{d_1}t^{d_2}G_{\bbd}(\theta)\df{h_2}{-d_2-1} \ast s^{m_1}t^{m_2} = 
G_{\bbd}(\bbm)(h_2(\bbm), -d_2-1)!\cdot  s^{d_1+m_1}t^{d_2+m_2}. 
\end{align}
If $\bbd\in\textbf{C1}$ and $s^{m_1}t^{m_2}\in J$, then the result of \eqref{eq:apply-diff-op-023} will remain in $J$. 
Thus, for all $\bbd\in\textbf{C1}$,   $D(R_A,J)_{\bbd}=s^{d_1}t^{d_2}\cdot\<G_{\bbd}(\theta)\>$.

Next, for ${\bbd}\in\textbf{C2}$,  
when $\bbm\in\NN A\cap \{x=-d_1\}$, then the result of \eqref{eq:apply-diff-op-023} is a constant times a monomial in $\sigma_1$, which does not belong to $J$.
Hence, for all $\bbd\in\textbf{C2}$, 
\[
\II(J)_{\bbd} = s^{d_1}t^{d_2}\cdot\< G_{\bbd}(\theta)(h_1+d_1)\>.
\] 
Further, if $\bbd\in \textbf{C2}\cup (\sigma_1\setminus\{\boldzero\})$, then 
\[
D(R_A,J)_{\bbd} = 
s^{d_1}t^{d_2}\cdot\< G_{\bbd}(\theta)(h_1+d_1)\>.
\]

For ${\bbd}\in\textbf{C4}$,  
when $\bbm\in\NN A\cap \{y=-d_2\}$, then the result of \eqref{eq:apply-diff-op-023} is a constant times a monomial in $\sigma_2$, which does not belong to $J$.
Hence, for all $\bbd\in\textbf{C4}$, 
\[
\II(J)_{\bbd} = s^{d_1}t^{d_2}\cdot\< G_{\bbd}(\theta) \df{h_2}{-d_2}\>.
\]
Further, if ${\bbd}\in\textbf{C4}\cup(\sigma_2\setminus\{\boldzero\})$, then  
\[
D(R_A,J)_{\bbd} = s^{d_1}t^{d_2}\cdot\<  \df{h_2}{-d_2}\>.
\]

Finally, for ${\bbd}\in\textbf{C3}$, 
if $\bbm\in \{x=-d_1\}\cap\{\bbc\mid s^{c_1}t^{c_1}\in J\}$ or if $\bbm\in\{y=-d_2\}\cap\NN A$, respectively, 
then the result of \eqref{eq:apply-diff-op-023} 
is either a constant times a monomial in $\sigma_1$ or $\sigma_2$, respectively.  
Hence, 
\[
\II(J)_{\bbd}=s^{d_1}t^{d_2}\cdot\< G_{\bbd}(\theta)(h_1+d_1)\df{h_2}{-d_2}\>.
\]
Finally, if ${\bbd}\in\textbf{C3}\cup (-\sigma_1)\cup(-\sigma_2)$, 
then 
\[
D(R_A,J)_{\bbd} = 
s^{d_1}t^{d_2}\cdot\< G_{\bbd}(\theta)(h_1+d_1)\df{h_2}{-d_2}\>.
\]

Gathering these computations, 
\[
\II(J)_{\bbd} = 
\begin{cases}
s^{d_1}t^{d_2} \cdot \< G_{\bbd}(\theta)\> 
&\text{ if } {\bbd} \in \textbf{C1},\\
s^{d_1}t^{d_2}\cdot \< G_{\bbd}(\theta)(h_1+d_1)\> 
&\text{ if } {\bbd} \in \textbf{C2},\\
s^{d_1}t^{d_2}\cdot \< G_{\bbd}(\theta)(h_1+d_1)\df{h_2}{-d_2}\> 
&\text{ if } {\bbd} \in \textbf{C3},\\
s^{d_1}t^{d_2}\cdot \< G_{\bbd}(\theta) \df{h_2}{-d_2} \> 
&\text{ if } {\bbd} \in \textbf{C4},
\end{cases}
\]
and 
\[
D(R_A,J)_{\bbd} = 
\begin{cases}
s^{d_1}t^{d_2} \cdot
\< G_{\bbd}(\theta)\> 
&\text{ if } 
\bbd\in\textbf{C1}\setminus(\sigma_1\cup\sigma_2),\\
s^{d_1}t^{d_2}\cdot
\< G_{\bbd}(\theta)(h_1+d_1)\> 
&\text{ if } {\bbd} \in (\textbf{C2}\setminus \sigma_2) \cup (\sigma_1 \cap \{d_2>0\}),\\
s^{d_1}t^{d_2}\cdot
\< G_{\bbd}(\theta)(h_1+d_1)\df{h_2}{-d_2}\> 
&\text{ if } {\bbd} \in \textbf{C3} \cup \{\boldzero\} \cup (-\sigma_1) \cup (-\sigma_2),\\
s^{d_1}t^{d_2}\cdot
\< G_{\bbd}(\theta) \df{h_2}{-d_2}\> 
&\text{ if } {\bbd} \in (\textbf{C4}\setminus \sigma_1) \cup  (\sigma_2 \cap \{d_1>0\}).
\end{cases}
\]
Taking the quotient of the above multigraded modules, we then obtain a description of the graded pieces of $\II(J)/D(R_A,J)$ as follows: 

\[
\left(\dfrac{\II(J)}{D(R_A,J)}\right)_{\bbd} = 
\begin{cases}
0 
&\text{ if } {\bbd} \notin \ZZ{\sigma_1} \cup \ZZ{\sigma_2},\\
\displaystyle\frac{\CC[\theta]}{\< h_1h_2 \>} 
&\text{ if } {\bbd}=\boldzero,\\
s^{d_1}t^{d_2} \cdot
\displaystyle\frac{\< G_{\bbd}(\theta) \>}{\< G_{\bbd}(\theta)h_2\>}, 
&\text{ if } d_1 \neq 0, d_2=0,\\
s^{d_1}t^{d_2}\cdot
\displaystyle\frac{\< G_{\bbd}(\theta) \>}{\< G_{\bbd}(\theta)h_1\>} 
&\text{ if } d_1=0, d_2 >0,\\
s^{d_1}t^{d_2}\cdot
\displaystyle\frac{ \< G_{\bbd}(\theta)\df{h_2}{-d_2}\>}{\< h_1G_{\bbd}(\theta)\df{h_2}{-d_2} \>} 
&\text{ if } d_1=0, d_2<0.
\end{cases}
\]

In this case, $JD(R_A)=D(R_A,J)$, so \eqref{eq:JDDJequal} happens to hold in this example. 
To see this, compute $(JD(R_A))_{\bbd}$ by looking at $\II(\Omega(\bbd-(2,1))$ and $\II(\Omega(\bbd-(3,1))$ for any $\bbd\in\ZZ^3$. 
For an example of how showing this equality works, when $(d_1-2,d_2-1)$ and $(d_1-3,d_2-1)$ are both in $\textbf{C2} \setminus \sigma_2$, then for $d_1 <0$,
\[
G_{(d_1-2,d_2-1)}(\theta)=\displaystyle\frac{\df{h_1}{3-d_1}}{(h_1-1)(h_1+d_1-2)}
\quad\text{and}\quad
G_{(d_1-3,d_2-1)}(\theta)=\displaystyle\frac{\df{h_1}{4-d_1}}{(h_1-1)(h_1+d_1-3)}.  
\]
Both of these polynomials are divisible by $G_{\bbd}(\theta)(h_1+d_1)$:
\[
G_{(d_1-2,d_2-1)}(\theta) 
= \displaystyle\frac{\df{h_1}{3-d_1}}{(h_1-1)(h_1+d_1-2)} 
= G_{\bbd}(\theta)(h_1+d_1)(h_1+d_1-3)
\] 
and
\[
G_{(d_1-2,d_2-1)}(\theta) 
= \displaystyle\frac{\df{h_1}{4-d_1}}{(h_1-1)(h_1+d_1-2)}
= G_{\bbd}(\theta)(h_1+d_1)(h_1+d_1-2)(h_1+d_1-4).
\] 
Further,  
\[
1=(h_1+d_1-3)^2-(h_1+d_1-2)(h_1+d_1-4),
\] 
which implies that 
\[
\left\< G_{\bbd}(\theta)(h_1+d_1)\right\> =\left\< G_{(d_1-2,d_2-1)}(\theta), G_{(d_1-3,d_2-1)}(\theta)
\right\>.
\]
Hence for $\bbd\in\textbf{C2}$, there is an equality
$D(R_A,J)_{\bbd} = s^{d_1}t^{d_2}\cdot \< G_{\bbd}(\theta)(h_1+d_1)\>=(JD(R_A))_{\bbd}$. 
\begin{ex}
\label{ex:nonnormal-nonGor}
Consider the matrix $\widetilde{A}=
\begin{bmatrix}
1&2&3&0\\
1&0&0&1
\end{bmatrix}$, which is associated to the semigroup ring $R_{\widetilde{A}}=\CC[\NN \widetilde{A}]=\CC[st,s^2,s^3,t]$, which is not normal, scored, or Gorenstein. 
As in Example~\ref{ex:one-line-gone}, we let 
\[
\sigma_1=\NN\{e_2\}
\quad\text{and}\quad
\sigma_2=\NN\{e_1\}
\]
be the integer points in the facets of the cone $\RR_{\geq 0} \widetilde{A}$. 
The primitive integral support functions of $\NN \widetilde{A}$ are 
\[
h_1=h_1(\theta)= \theta_1
\quad\text{and}\quad
h_2=h_2(\theta)=\theta_2.
\]
The prime ideals associated to the facets of $\widetilde{A}$ are 
\[
P_1=\<st,s^2,s^3\>
\quad\text{and}\quad
P_2=\<st,t\>.
\]
Finally, set $\widetilde{J}=P_1 \cap P_2=\<st, s^2t\>$. 
\end{ex}

\begin{wrapfigure}[14]{r}[10pt]{2in}
\centering
\begin{tikzpicture}[scale=0.6]
\draw[blue, thick] (-3.5,0) -- (3.5,0);
\draw[blue, thick] (0,-3.5) -- (0,3.5);
\filldraw[cyan] (1,0) circle (2pt);
\foreach \y in {1,2,3}{ \node[draw,circle,inner sep=2pt,red,fill] at (2,\y) {};}
\foreach \y in {1,2,3}{ \node[draw,circle,inner sep=2pt,red,fill] at (3,\y) {};}
\foreach \y in {1,2,3}{ \node[draw,circle,inner sep=2pt,red,fill] at (1,\y) {};}
\foreach \x in {0,2,3}{ \node[draw,circle,inner sep=2pt,blue,fill] at (\x,0) {};}
\foreach \y in {0,1,2,3}{ \node[draw,circle,inner sep=2pt,blue,fill] at (0,\y) {};}
\foreach \x in {-3,-2,...,3}{
     \foreach \y in {-3,-2,...,3}{
       \node[draw,circle,inner sep=2pt] at (\x,\y) {};}}
\foreach \x in {-3,-2,...,-1}{
     \foreach \y in {0,1,2,3}{
       \node[draw,circle,inner sep=2pt,yellow,fill] at (\x,\y) {};}}
\foreach \x in {-3,-2,...,-1}{
     \foreach \y in {-3,-2,-1}{
       \node[draw,circle,inner sep=2pt,violet,fill] at (\x,\y) {};}}     
\foreach \x in {0,1,2,3}{
     \foreach \y in {-3,-2,-1}{
       \node[draw,circle,inner sep=2pt,green,fill] at (\x,\y) {};}}
\end{tikzpicture}
\captionsetup{margin=0in,width=1.25in,font=small,labelsep=newline,justification=centering}
\caption{\ Chambers for $\CC[st,s^2,s^3,t]$}
\label{fig:nonnormal:2nd-lattice}
\end{wrapfigure}
Note that $R_A/J$ from Example \ref{ex:one-line-gone} is isomorphic to $R_{\widetilde{A}}/\widetilde{J}$.  
Hence, 
we would also expect that $\II(J)/D(R_A,J) \cong \II(J)/D(R_{\widetilde{A}},\widetilde{J})$. 
We will determine the graded pieces of $\II(\widetilde{J})$ and $D\left(R_{\widetilde{A}},\widetilde{J}\right)$, as well as ${\II(\widetilde{J})}/{D(R_{\widetilde{A}},\widetilde{J})}$; 
then, we will exhibit an isomorphism between  the graded pieces of ${\II(\widetilde{J})}/{D(R_{\widetilde{A}},\widetilde{J})}$ and $\II(J)/D(R_A,J)$.

We will again use \emph{chambers} of $\ZZ^2$ to aid in our computations of the graded pieces of $\II(J)/D(R_A,J)$. 
In Figure~\ref{fig:nonnormal:2nd-lattice}, the red multidegrees correspond to monomials in $\widetilde{J}$, while the blue ones (both light and dark) are the remaining $\bbd\in\ZZ^2$ for which $D(R_{\widetilde{A}})_{\bbd}=\II(\widetilde{J})_{\bbd}$; together, these red and blue multidegrees form chamber \textbf{C1}.  
The remaining multidegrees are split into chambers as in the previous examples: 
yellow points form \textbf{C2}, 
violet points make \textbf{C3}, 
and green points give \textbf{C4}.

In \cite[Examples 3.2.7 and 3.3.4]{Sai-Tr-DASR}, 
Saito and Traves found that for any $\bbd\in\ZZ^2$, 
$D(R_{\widetilde{A}})_{\bbd}= 
s^{d_1}t^{d_2}\cdot  \II(\Omega({\bbd}))$, where 
\begin{align}
\label{eq:IIomega-nonGor}
\II(\Omega({\bbd})) = 
\begin{cases}
\left\< \prod\limits_{i=1}^2 \df{h_i}{h_i({-\bbd})-1}\right\> 
&\text{ if } {\bbd} \notin e_1 -\NN \widetilde{A},\\
\prod\limits_{i=1}^2 \df{h_i}{h_i(-{\bbd})-1}\cdot &\\
\qquad
\left\< h_1+d_1-1,h_2+d_2 \right\> 
&\text{ if } {\bbd} \in e_1 -\NN \widetilde{A}.
\end{cases}
\end{align}
To compute $\II(\widetilde{J})$, we first determine which elements of $\II(\Omega(d))$ are already in $\II(\widetilde{J})$. 
To do this, we first break up \eqref{eq:IIomega-nonGor} into cases by chamber: 
\[
\II(\widetilde{J})_{\bbd} = 
\begin{cases}
s^{d_1}t^{d_2}\cdot 
\II(\Omega({\bbd})) 
&\text{ if } {\bbd} \in\textbf{C1},\\
s^{d_1}t^{d_2}\cdot
\< \df{h_1}{h_1({-\bbd})}\> 
&\text{ if } 
\bbd\in\textbf{C2}\setminus(-\sigma_2),\\
s^{d_1}t^{d_2}\cdot 
\df{h_1}{h_1({-\bbd})}\cdot \< h_1+d_1-1,h_2 \> 
&\text { if } 
\bbd\in(-\sigma_2),\\
s^{d_1}t^{d_2}\cdot 
\left\< \df{h_1}{h_1({-\bbd})}\df{h_2}{h_2({-\bbd})}\right\> 
&\text{ if } 
\bbd\in\textbf{C3},\\
s^{d_1}t^{d_2}\cdot 
\< \df{h_2}{h_2({-\bbd})}\> 
&\text{ if } 
\bbd\in\textbf{C4}.
\end{cases}
\]
Similarly, 
\[
D(R_{\widetilde{A}},\widetilde{J})_{\bbd} = 
\begin{cases}
s^{d_1}t^{d_2}\cdot \CC[\theta] 
&\text{ if } s^{d_1}t^{d_2} \in \widetilde{J},\\
s^{d_1}t^{d_2}\cdot \< \df{h_1}{h_1(-{\bbd})}\> 
&\text{ if } 
\bbd\in(\textbf{C2}\setminus(-\sigma_2))\cup(\sigma_1\setminus\{\boldzero\}),\\
s^{d_1}t^{d_2}\cdot \< \df{h_2}{h_2({-\bbd})}\> 
&\text{ if } 
\bbd\in\textbf{C4}\setminus(-\sigma_1),\\
s^{d_1}t^{d_2}\cdot \< \df{h_1}{h_1({-\bbd})}\df{h_2}{h_2({-\bbd})}\> 
&\text{ if } 
\bbd\in\textbf{C3}\cup(-\sigma_1)\cup(-\sigma_2).
\end{cases}
\]
Putting these together, the graded piece of $\II(\widetilde{J})/D(R_{\widetilde{A}},\widetilde{J})$ at ${\bbd} \in \ZZ^2$ is
\[
\displaystyle\frac{\II(\widetilde{J})}{D(R_{\widetilde{A}},\widetilde{J})}_{\bbd} = 
\begin{cases}
0 
&\text{ if } {\bbd} \notin \ZZ\sigma_1 \cup \ZZ\sigma_2,
\\
s^{d_1}t^{d_2}\cdot\dfrac{\CC[\theta]}{\< h_1\> }
&\text{ if } 
\bbd\in\sigma_1\setminus\{\boldzero\},
\\
s^{d_1}t^{d_2}\cdot\dfrac{\CC[\theta]}{\< h_2\>}
&\text{ if } 
\bbd\in\sigma_2\setminus\{\boldzero,(1,0)\},
\\
s^{d_1}t^{d_2}\cdot\dfrac{\CC[\theta]}{\<\theta_1\theta_2\>}
&\text{ if }\bbd=\boldzero,
\\
s^{d_1}t^{d_2}\cdot\dfrac{\< h_1,h_2\>}{\< h_2\> }
&\text{ if }  {\bbd}=(1,0),
\\
s^{d_1}t^{d_2}\cdot\dfrac{\df{h_1}{h_1({-\bbd})}\< h_1+d_1-1,h_2 \>}{\< \df{h_1}{h_1({-\bbd})}h_2\> }
&\text{ if } 
\bbd\in(-\sigma_2)\setminus\{\boldzero\}, 
\\
s^{d_1}t^{d_2}\cdot\dfrac{\< \df{h_2}{h_2({-\bbd})}\>}{\< h_1 \df{h_2}{h_2({-\bbd})}\> }
&\text{ if } 
\bbd\in(-\sigma_1)\setminus\{\boldzero\}.
\end{cases}
\]

Looking at the graded pieces of $\II(J)/D(R_A,J)$ from Example \ref{ex:one-line-gone} and noting that $G_{\bbd}(\theta)=\CC[\theta]$ for $\bbd \in (\sigma_1) \cup (\sigma_2 \setminus\{ \textbf{0},(1,0)\}$, it follows that $\left(\II(J)/D(R_A,J)\right)_{\bbd}$ is identical to  $\left({\II(\widetilde{J})}/{D(R_{\widetilde{A}},\widetilde{J})}\right)_{\bbd}$ for all $\bbd$ except $\bbd=(1,0)$ and $\bbd \in (-\sigma_2) \setminus \{\textbf{0}\}$. 
So we only need to exhibit that the graded pieces are isomorphic for $\bbd \in \{(1,0)\} \cup (-\sigma_2) \setminus \{\textbf{0}\}$.  Note that for  $d_1 <0$ and $d_1 \neq -1$, 
\[ 
\displaystyle\frac{\left\<\displaystyle\frac{\df{h_1}{1-d_1}}{(h_1-1)(h_1+d_1)}\right\>}{\left\< \displaystyle\frac{\df{h_1}{1-d_1}\, h_2}{(h_1-1)(h_1+d_1)}\right\>} \cong \displaystyle\frac{\<\df{h_1}{1-d_1}\>}{\<\df{h_1}{1-d_1}\, h_2\>}
\]
by the Third Isomorphism Theorem, and  
\[
\dfrac{\df{h_1}{h_1({-\bbd})}\< h_1+d_1-1,h_2 \>}{\< \df{h_1}{h_1({-\bbd})}h_2\> }\cong \displaystyle\frac{\<\df{h_1}{1-d_1}\>}{\<\df{h_1}{1-d_1}\, h_2\>} 
\]
by the Second Isomorphism Theorem.
A similar argument gives an isomorphism in 
the remaining multidegrees. 
Combining the information on all graded pieces, 
we have explicitly shown that 
\[
\dfrac{\II(\widetilde{J})}{D\left(R_{\widetilde{A}},\widetilde{J}\right)}
\cong 
\dfrac{\II(J)}{D\left(R_{A},J\right)}.
\]

\bibliographystyle{amsalpha}
\bibliography{refs}
\end{document}